\documentclass[11pt,reqno]{amsart}

\usepackage{fixltx2e}                   
\usepackage{tgpagella} 

\usepackage[utf8]{inputenc}            
\usepackage{amsmath}                     
\usepackage{amssymb, latexsym, stmaryrd, amsthm, dsfont, amsfonts, amsbsy }            

\usepackage{mathtools}                  
\usepackage{bm}                          
\usepackage{enumerate}                   
\usepackage{verbatim}                    
\usepackage{url}   
\usepackage{lscape}                      
                     
\usepackage{soul}

\usepackage[margin=1.09in]{geometry}

\usepackage[all]{xy}
\usepackage{microtype,tikz,tikz-cd}  
\makeatletter                            
\def\MT@register@subst@font{\MT@exp@one@n\MT@in@clist\font@name\MT@font@list
 \ifMT@inlist@\else\xdef\MT@font@list{\MT@font@list\font@name,}\fi}
\makeatother

\usepackage[capitalize]{cleveref}

\newcommand{\bit}{\begin{itemize}}    
\newcommand{\eit}{\end{itemize}}
\newcommand{\ben}{\begin{enumerate}}
\newcommand{\een}{\end{enumerate}}

\newcommand{\benroman}{\ben[\normalfont (i)]}  
\let\eroman\een

\newcommand{\bde}{\begin{description}}
\newcommand{\ede}{\end{description}}

\theoremstyle{theorem}
\newtheorem{Theorem}{Theorem}[section]
\newtheorem{Theorem-n}{Theorem}

\newtheorem{Proposition}[Theorem]{Proposition}
\newtheorem{Antidichotomy Theorem}[Theorem]{Antidichotomy Theorem}
\newtheorem{Jankov Theorem}[Theorem]{Jankov Theorem}
\newtheorem{Jankov Lemma}[Theorem]{Jankov Lemma}
\newtheorem{Dual Jankov Lemma}[Theorem]{Dual Jankov Lemma}
\newtheorem{Finite Esakia Duality}[Theorem]{Finite Esakia Duality}
\newtheorem{Fine's Completeness Theorem}[Theorem]{Fine's Completeness Theorem}
\newtheorem{Esakia Duality}[Theorem]{Esakia Duality}
\newtheorem{Width Dichotomy Theorem}[Theorem]{Width Dichotomy Theorem}
\newtheorem{FMP Dichotomy Theorem}[Theorem]{FMP Dichotomy Theorem}
\newtheorem{Modal Antidichotomy Theorem}[Theorem]{Modal Antidichotomy Theorem}
\newtheorem{Blok-Esakia Theorem}[Theorem]{Blok-Esakia Theorem}

\newtheorem{Lemma}[Theorem]{Lemma}
\newtheorem{Corollary}[Theorem]{Corollary}
\newtheorem{Claim}[Theorem]{Claim}

\theoremstyle{definition}
\newtheorem{Definition}[Theorem]{Definition}
\newtheorem{exa}[Theorem]{Example}

\theoremstyle{remark}

\newtheorem{Remark}[Theorem]{Remark}

\let\leq=\leqslant
\let\nleq=\nleqslant
\let\geq=\geqslant

\bmdefine{\A}{A} 
\bmdefine{\C}{C}                                
                              
\bmdefine{\2}{2}
\bmdefine{\B}{B}
\bmdefine{\D}{D}

\usepackage{scalerel}

\subjclass[2010]{03B55, 03B45, 06D20, 06E15}

\keywords{Superintuitionistic logics, modal logics, Kripke completeness, the finite model property}

\begin{document}
		
\title[Degrees of the finite model property: The Antidichotomy Theorem]{Degrees of the finite model property: \\ The Antidichotomy Theorem}

\author{Guram Bezhanishvili, Nick Bezhanishvili, and Tommaso Moraschini}

\address{Guram Bezhanishvili: Department of Mathematical Sciences, New Mexico State University, Las Cruces NM $88003$, USA}\email{guram@nmsu.edu}

\address{Nick Bezhanishvili: Institute for Logic, Language and Computation, University of Amsterdam, Postbus $94242$, $1090$GE Amsterdam, The Netherlands}\email{N.Bezhanishvili@uva.nl}

\address{Tommaso Moraschini: Departament de Filosofia, Facultat de Filosofia, Universitat de Barcelona (UB), Carrer Montalegre, $6$, $08001$ Barcelona, Spain}\email{tommaso.moraschini@ub.edu}

\date{}

\begin{abstract}
A classic result in modal logic, known as the Blok Dichotomy Theorem,  states that the degree of incompleteness of a normal extension of the basic modal logic $\sf K$ is $1$ or $2^{\aleph_0}$. It is a long-standing open problem whether Blok Dichotomy holds for normal extensions of other prominent modal logics (such as $\sf S4$ or  $\sf K4$) or for extensions of the intuitionistic propositional calculus $\mathsf{IPC}$ (see \cite[Prob.~10.5]{ChZa97}).  In this paper,  we introduce the notion of the degree of finite model property (fmp), which is a natural variation of the degree of incompleteness.  It is a consequence of the Blok Dichotomy Theorem that the degree of fmp of a normal extension of $\sf K$ remains $1$ or $2^{\aleph_0}$. In contrast, our main result establishes the following Antidichotomy Theorem for the degree of fmp for extensions of $\mathsf{IPC}$: each nonzero cardinal $\kappa$  such that $\kappa \leq \aleph_0$ or $\kappa = 2^{\aleph_0}$  is realized as the degree of fmp of some extension of $\mathsf{IPC}$. 
We then use the Blok-Esakia theorem to establish the same Antidichotomy Theorem for normal extensions of $\sf S4$ and $\sf K4$. This provides a solution of the reformulation of \cite[Prob.~10.5]{ChZa97} for the degree of fmp.  
\end{abstract}

\maketitle

\tableofcontents

\section{Introduction}

 Since its inception in the late 1950s/early 1960s, Kripke semantics became the most popular tool to study modal and intuitionistic logics.  
However, examples of Kripke incomplete logics began to emerge already in the 1970s (see, e.g., \cite[Ch.~6]{ChZa97}). 
In order to shed light on the phenomenon of Kripke incompleteness, 
Fine \cite{Fin74} associated with each normal modal logic $\mathsf{L}$ a cardinal that measures the degree of incompleteness of  $\mathsf{L}$. More precisely, let $\sf Fr(L)$ be the class of Kripke frames validating $\sf L$. We say that the {\em degree of incompleteness} of $\sf L$ is the cardinal $\kappa$ if there are exactly $\kappa$ logics $\sf L'$ such that $\sf Fr(L')=Fr(L)$.  
 Notice  that all but one of these $\mathsf{L}'$
 are Kripke incomplete. 
 
 Blok \cite{B78,Blok-dichotomy78} gave a very unexpected  characterization of degrees of incompleteness,  \color{black} which became known as the {\em Blok Dichotomy Theorem}. It states that a normal modal logic $\sf L$ has the degree of incompleteness either 1 or $2^{\aleph_0}$; 
it is 1 iff $\sf L$ is a join-splitting logic (see Section~\ref{sec: jankov} for the definition); otherwise 
it is $2^{\aleph_0}$. Chagrova \cite{MR1688513} proved that the Blok Dichotomy Theorem also holds for the neighbourhood semantics. 
We refer to  \cite{RZW06} and \cite{Lit08} 
 for a detailed discussion of Blok Dichotomy and its importance in modal logic.

Blok's result implies that some of the most studied normal modal logics, such as $\sf K4$ (the logic of transitive Kripke frames) and $\sf S4$ (the logic of reflexive and transitive Kripke frames), have the degree of incompleteness $2^{\aleph_0}$. However, the logics sharing the Kripke frames with $\sf K4$ and $\sf S4$ are not necessarily normal extensions of $\sf K4$ or $\sf S4$. Thus, Blok's result does not automatically transfer to normal extensions of $\sf K4$ or $\sf S4$ (or, more generally, to normal extensions of a given normal modal logic). There have been several attempts to investigate Blok Dichotomy for normal extensions of $\sf K4$ and $\sf S4$. However, this remains an outstanding open problem in modal logic \cite[Prob.~10.5]{ChZa97}. 

For a logic $\sf L$, let $\sf Fin(L)$ be the class of finite Kripke frames validating $\sf L$. We recall that $\sf L$ has the {\em finite model property} ({\em fmp} for short) if $\sf L$ is complete with respect to $\sf Fin(L)$. Clearly each logic with the fmp is Kripke complete. Taking inspiration from degrees of incompleteness, it is natural to introduce a similar concept for the fmp. 
We say that the {\em degree of fmp} of a logic $\sf L$ is $\kappa$ provided there exist exactly  $\kappa$ logics $\sf L'$ such that $\sf Fin(L')=Fin(L)$. 
 As with the degree of incompleteness, all but one of such $\mathsf{L}'$ lack the fmp. 
Our main result establishes a complete opposite of Blok Dichotomy theorem for superintuitionistic logics and  transitive (normal) modal logics.
Namely, we prove that if $\kappa$ is a nonzero cardinal such that $\kappa \leq \aleph_0$ or $\kappa = 2^{\aleph_0}$, then there exists a  superintuitionistic logic (or a transitive modal logic) $\sf L$ 
such that the degree of fmp of $\sf L$ is $\kappa$. Under the Continuum Hypothesis (CH) this implies that each  nonzero \color{black} $\kappa \leq 2^{\aleph_0}$ is realized as the degree of fmp of some  superintuitionistic \color{black} logic 
(or some transitive modal logic).
For this reason, we refer to these results as the {\em Antidichotomy Theorems for degrees of fmp} (see Theorems \ref{Thm:dichotomy:IPC:main} and \ref{Thm:modal:dichotomy:main}).

In \cite[p.~409]{Lit08} Litak asks ``if there is any nontrivial completeness notion for which the Blok Dichotomy does not hold." Our main result provides such a nontrivial and, in our opinion, very natural notion for superintuitionistic logics and transitive modal logics.
 It also provides a solution of a variant of \cite[Prob.~10.5]{ChZa97} when the degree of incompleteness is replaced by the degree of fmp. \color{black}

To give more context, we recall that {\em superintuitionistic logics} are (axiomatic) extensions of the intuitionistic propositional calculus $\sf IPC$. They have been studied extensively in the literature (see, e.g., \cite{ChZa97}). In particular, there is a close connection between superintuitionistic  logics and normal extensions of $\sf S4$. The {\em G\"odel translation} embeds $\sf IPC$ into $\sf S4$ 
faithfully \cite{McKT48}. Thus, each superintuitionistic logic $\sf L$ is embedded 
into a normal extension of $\sf S4$, called a {\em modal companion} of $\sf L$ \cite[Sec.~9.6]{ChZa97}.
Each $\sf L$ has many modal companions, but remarkably each $\sf L$ possesses a \color{black}  largest modal companion.
By Esakia's theorem \cite{Esakia76,MR579150}, the largest modal companion of $\sf IPC$ is the well-known Grzegorczyk logic $\sf Grz$. Consequently, the largest modal companion of each superintuitionistic logic is a normal extension of $\sf Grz$, and there  exists an isomorphism between the lattice of superintuitionistic logics and the lattice of normal extensions of $\sf Grz$ \color{black} (the Blok-Esakia theorem) \cite{Blok76,Esakia76}. 

Notice that it is a consequence of the Blok Dichotomy Theorem that the degree of fmp of a normal extension of the basic modal logic $\sf K$ remains $1$ or $2^{\aleph_0}$. Thus, in the lattice of all normal modal logics the dichotomy holds also for the degrees of fmp (see Theorem \ref{Thm:strict:dichotomy:modal}).  In contrast, it is a consequence of our Modal  Antidichotomy Theorem that the situation is drastically different for transitive modal logics (see Corollary \ref{Cor:antidic:S4:K4}).

We conclude the introduction by discussing how we establish our main results.
We first prove the Antidichotomy Theorem for degrees of fmp of superintuitionistic logics. 
We heavily rely on Esakia duality for Heyting algebras 
\cite{Esakia-book85}, as well as on Fine's completeness theorem for logics of bounded width \cite{Fin74b} and the theory of splittings \cite[Sec.~10.5]{ChZa97}. Our proof is broken into two parts, depending on whether $\kappa\leq\aleph_0$ or $\kappa=2^{\aleph_0}$. 

When $\kappa\leq\aleph_0$ 
we work with extensions of the superintuitionistic logic $\sf KG$, which was introduced by Kuznetsov and Ger\v{c}iu \cite{GeKuz70,KuzGer70a} and bears  their name. The logic $\sf KG$ is the logic of  sums of one-generated Heyting algebras, the combinatorics of which allows to  
 construct \color{black}extensions of $\sf KG$ that lack the fmp 
\cite{KuzGer70a,Kr93b,BBdeJ08}. 
First, we use Fine’s completeness theorem to prove that $\mathsf{KG}$ is a join-splitting logic over $\mathsf{IPC}$ (for a similar result see \cite{Kr93b}). Then we develop a method, utilizing a technique of \cite{BBdeJ08}, that produces an extension $\mathsf{L}$ of $\mathsf{KG}$ whose degree of fmp is $\kappa$ for every nonzero cardinal $\kappa \leq \aleph_0$.

To show that there exist superintuitionistic  logics whose degree of fmp is $2^{\aleph_0}$ we work with superintuitionistic logics of finite width. Transitive modal logics of finite width were introduced by Fine \cite{Fin74b} who showed that each transitive modal logic of finite width  is Kripke complete.  
The concept was adapted to superintuitionistic logics by Sobolev \cite{So77}.  For every positive integer $n$, let $\mathsf{BW}_n$ be the least superintuitionistic logic of width $n$. These logics are known to have the fmp (see \cite{Fine85}). 
We prove that if $n>2$, then the degree of fmp of $\mathsf{BW}_n$ is $2^{\aleph_0}$. This is done by a careful analysis of the combinatorics of posets of bounded width.

Under CH our results show that  for every nonzero cardinal $\kappa\leq 2^{\aleph_0}$ there exists a superintuitionistic logic $\sf L$ whose degree of fmp is $\kappa$, thus yielding the Antidichotomy Theorem for degrees of fmp of superintuitionistic logics. Nonetheless, determining the degree of fmp of a given superintuitionistic logic remains an interesting open problem.

Finally, we transfer our results to the setting of modal logics. Following the notation of \cite{ChZa97}, for a normal modal logic $\mathsf{L}$, let $\textup{Next}\,\mathsf{L}$ be the lattice of normal extensions of $\mathsf{L}$. We first use the Blok-Esakia theorem to prove our Antidichotomy Theorem for $\textup{Next}\,\mathsf{Grz}$. We next show that for each normal modal logic $\mathsf{L}\subseteq\mathsf{Grz}$ with the fmp, the Antidichotomy Theorem holds for $\textup{Next}\,\mathsf{L}$ provided $\mathsf{Grz}$ is a join-splitting logic above $\mathsf{L}$. Since $\mathsf{S4}$ and $\mathsf{K4}$ have the fmp and $\mathsf{Grz}$ is a join-splitting logic above both, it follows that the Antidichotomy Theorem holds for $\textup{Next}\,\mathsf{S4}$ and $\textup{Next}\,\mathsf{K4}$. We conclude the paper by listing several open problems and possible future research directions.

\section{Superintuitionistic logics}

We recall that a {\em superintuitionistic logic}, or a {\em si-logic} for short, is a set of formulas $\sf L$ containing $\sf IPC$ and closed under the inference rules of modus ponens and substitution. It is well known 
(see, e.g., \cite[Thm.~4.1]{ChZa97}) that consistent si-logics are exactly the logics situated between $\sf IPC$ and the classical propositional calculus $\sf CPC$. Thus, consistent si-logics are often referred to as {\em intermediate logics}. Given a set of formulas $\Sigma$, we denote by $\mathsf{IPC} + \Sigma$ the si-logic \textit{axiomatized} by $\Sigma$; that is, the least si-logic containing $\Sigma$.

 When ordered by set inclusion, the set of si-logics forms a complete lattice, denoted by $\textup{Ext}\,\mathsf{IPC}$, whose bottom and top are $\sf IPC$ and the inconsistent logic, respectively. The meet and join operations in $\textup{Ext}\,\mathsf{IPC}$ are defined as
\[
\bigwedge_{i \in I}\mathsf{L}_i = \bigcap_{i \in I}\mathsf{L}_i \quad \text{and} \quad
\bigvee_{i \in I}\mathsf{L}_i = \text{the si-logic axiomatized by }\bigcup_{i \in I}\mathsf{L}_i.
\]
It is a well-known result of Jankov \cite{Jankov68} that the cardinality of $\textup{Ext}\,\mathsf{IPC}$ is $2^{\aleph_0}$.

Kripke semantics for si-logics is given by partially ordered sets (posets for short). For a poset $X$, we call $U\subseteq X$ an {\em upset} (upward closed set) if 
\[
x\in U \mbox{ and } x\leq y \mbox{ imply } y\in U. 
\]
A {\em valuation} $\nu$ on $X$ assigns to each propositional letter $p$ an upset of $X$. For $x\in X$ and a formula $\varphi$ we write  $x\Vdash_\nu\varphi$ \color{black} when $x$ satisfies $\varphi$ under $\nu$. As usual, the satisfaction relation $\Vdash$ is defined by recursion on the construction of formulas:
\[
\begin{array}{lll}
x\nVdash_\nu\bot & & \\
x\Vdash_\nu p & \mbox{iff} & x\in\nu(p) \\
x\Vdash_\nu \varphi\wedge\psi & \mbox{iff} & x\Vdash_\nu\varphi \mbox{ and } x\Vdash_\nu\psi \\
x\Vdash_\nu \varphi\vee\psi & \mbox{iff} & x\Vdash_\nu\varphi \mbox{ or } x\Vdash_\nu\psi \\
x\Vdash_\nu \varphi\to\psi & \mbox{iff} & \forall y(x\leq y \mbox{ and } y\Vdash_\nu\varphi \mbox{ imply } y\Vdash_\nu\psi).
\end{array}
\]
A formula $\varphi$ is said to be \textit{true} in $X$ under $\nu$ if $x \Vdash_\nu \varphi$ for every $x \in X$ and it is said to be \textit{valid} in $X$ if it is true under each valuation, in which case we write $X \vDash \varphi$. \color{black}

Algebraic semantics for si-logics is given by Heyting algebras. We recall that a {\em Heyting algebra}  $\A = \langle A; \land, \lor, \to, 0, 1 \rangle$ \color{black} is a bounded distributive lattice such that $\land$ has a residual $\to$ given by
\[
a \land b \leq c \Longleftrightarrow a \leq b \to c
\]
for all $a, b, c \in A$. 

A valuation $\nu$ in a Heyting algebra $\A$ assigns to each propositional letter an element of $\A$. The logical connectives are then interpreted as the corresponding operations in $\A$. A formula $\varphi$ is {\em true} in $\A$ under $\nu$ if $\nu(\varphi)=1$ and it is {\em valid} in $\A$ if it is true under each valuation, in which case we write $\A \vDash \varphi$. 

There is a close connection between Kripke and algebraic semantics for si-logics. For a poset $X$ and $U \subseteq X$, let
\begin{align*}
{\uparrow} U &= \{x\in X :  \exists u\in U \mbox{ with } u\leq x\} \\
{\downarrow} U &= \{x\in X :  \exists u\in U \mbox{ with } x\leq u\}.
\end{align*}
If $U = \{ x \}$, we simply write ${\uparrow} x$ and ${\downarrow} x$ instead of ${\uparrow} \{ x \}$ and ${\downarrow} \{ x\}$. 
Let ${\sf Up}(X)$ be the set of upsets of $X$. Then ${\sf Up}(X)$ is a Heyting algebra where join and meet are set-theoretic union and intersection, bottom and top are $\varnothing$ and $X$, and $\to$ is defined by
\[
U \to V = X \smallsetminus {\downarrow} (U \smallsetminus V) = \{ x\in X : {\uparrow}x\cap U\subseteq V \}.
\]

Conversely, for a Heyting algebra $\A$, let $X_{\A}$ be the poset of prime filters of $\A$ ordered by inclusion. Define $\gamma_{\A} \colon \A\to{\sf Up}(X_{\A})$ by
\[
\gamma_{\A}(a)=\{ x\in X_{\A} : a\in x \}.
\]
Then $\gamma_\A$ is a Heyting algebra embedding. To recognize the image of $\A$ in ${\sf Up}(X_{\A})$, we introduce the topology $\tau$ on $X_{\A}$ given by the subbasis
\[
\{ \gamma_\A(a) : a \in A \} \cup \{ X_\A \smallsetminus \gamma_\A(a) : a \in A \}.
\]
It is well known that $\tau$ is a Stone topology on $X_{\A}$ (that is, it is compact, Hausdorff, and zero-dimensional). The triple $\A_{\ast} = \langle X_{\A}, \tau, \subseteq \rangle$ is known as the {\em Esakia space} of $\A$. The map $\gamma_\A$ is an isomorphism from $\A$ onto the Heyting algebra of clopen upsets of $\A_{\ast}$. Thus, each Heyting algebra is represented as the algebra of clopen upsets of an Esakia space. 

Esakia spaces are characterized abstractly as triples $X=\langle X, \tau, \leq \rangle$ where $\tau$ is a Stone topology and $\leq$ is a partial order on $X$ that, moreover, is {\em continuous} in the sense that
\begin{enumerate}
\item ${\uparrow} x$ is closed for all $x \in X$;
\item $U\subseteq X$ is clopen implies ${\downarrow} U$ is clopen.
\end{enumerate}

We point out that the partial order $\leq$ is continuous iff the corresponding map $\rho:X \to \mathcal V X$ from $X$ to the Vietoris space $\mathcal V X$, given by $\rho(x)={\uparrow}x$, is a well-defined continuous map \cite{Es74,Ab05b,KupKuVe04}. 

We thus obtain the object level of Esakia duality, namely that there is a one-to-one  correspondence between Heyting algebras and Esakia spaces. To extend this correspondence to full duality, we recall that a {\em p-morphism} (or {\em bounded morphism}) between two posets $X$ and $Y$ is a map $\alpha \colon X\to Y$ such that ${\uparrow}\alpha(x)=\alpha({\uparrow}x)$ for each $x\in X$. 

Let $\mathsf{ES}$ be the category of Esakia spaces and continuous p-morphisms between them. Let also $\mathsf{HA}$ be the category of Heyting algebras and Heyting homomorphisms between them. The two categories are related as follows \cite{Es74,Esakia-book85}: 

\begin{Theorem}[\textbf{Esakia Duality}]
$\mathsf{HA}$ is dually equivalent to $\mathsf{ES}$.
\end{Theorem}

We denote the contravariant functors establishing Esakia duality by $(-)_*:\mathsf{HA}\to\mathsf{ES}$ and $(-)^*:\mathsf{ES}\to\mathsf{HA}$. The functor $(-)_*$ assigns to each Heyting algebra $\A$ the Esakia space $\A_*$.
If $f \colon \A \to \B$ is a Heyting homomorphism, define $f_{\ast} \colon \B_{\ast} \to \A_{\ast}$ by $f_{\ast}(x) = f^{-1}(x)$ for all $x \in  B_{\ast}$. Then $f_{\ast}$ is a continuous p-morphism and $(-)_*$ assigns $f_*$ to $f$. 

The functor $(-)^*$ assigns to an Esakia space $X$ the Heyting algebra $X^{\ast}$ of clopen upsets of $X$.
If $\alpha \colon X \to Y$ is a continuous p-morphism, define $\alpha^{\ast} \colon Y^{\ast} \to X^{\ast}$ by $\alpha^{\ast}(U) = \alpha^{-1}(U)$ for all $U \in Y^{\ast}$. Then $\alpha^\ast$ is a Heyting homomorphism 
and $(-)^*$ assigns $\alpha^*$ to $\alpha$.

The topology of a finite Esakia space is discrete (since it is Hausdorff). Therefore, the full subcategory of $\mathsf{ES}$ consisting of finite Esakia spaces is isomorphic to the category of finite posets and p-morphisms between them. Consequently, in the finite case, Esakia duality restricts to the following \cite{Es74,Esakia-book85}:

\begin{Theorem}[\textbf{Finite Esakia Duality}]
The category of finite Heyting algebras and Heyting homomorphisms is dually equivalent to the category of finite posets and p-morphisms between them.
\label{thm: FED}
\end{Theorem}

In view of Esakia duality, we can define the notion of validity for Esakia spaces as follows. We say that a formula $\varphi$ is \textit{valid} in an Esakia space $X$, and write $X \vDash \varphi$, when it is valid in the Heyting algebra $X^\ast$. This allows us to associate an si-logic with each class of Esakia spaces (resp.\ Heyting algebras or posets) as follows.

\begin{Definition}
Let $K$ be a class of Esakia spaces (resp.\ Heyting algebras or posets). The \textit{logic of} $K$, in symbols $\mathsf{Log}(K)$, is the set of formulas valid in each member of $K$.
\end{Definition}

Notice that $\mathsf{Log}(K)$ is always an si-logic. While every si-logic has the form  $\mathsf{Log}(K)$ for some class of Esakia spaces (resp.\ Heyting algebras), the logics of the form $\mathsf{Log}(K)$ for a class $K$ of posets are precisely the Kripke complete ones.

We conclude this preliminary section by a brief dual description of homomorphic images and subalgebras of Heyting algebras. Henceforth, we will freely use these results. To this end, we recall that if $\alpha \colon X \to Y$ is a p-morphism between posets, the map $\alpha^{-1} \colon \mathsf{Up}(Y) \to \mathsf{Up}(X)$ is a Heyting homomorphism that, moreover, is \emph{complete} (i.e., it preserves arbitrary meets and joins).
For part (1) of the next result see \cite[Lem.~3.3.13(3)]{Esakia-book85}, and for part (2) see \cite[Thms.~3.4, 3.5, 4.6]{MR0197372}. 

\begin{Theorem} \label{thm: 1-1 and onto}
The following conditions hold.
\begin{enumerate}
\item Let $X$ and $Y$ be Esakia spaces, $\alpha \colon X\to Y$ a continuous p-morphism, and $\alpha^{-1} \colon Y^\ast\to X^\ast$ the corresponding Heyting homomorphism. Then $\alpha^{-1}$ is one-to-one iff $\alpha$ is onto, and $\alpha^{-1}$ is onto iff $\alpha$ is one-to-one. 
\item Let $X$ and $Y$ be posets, $\alpha \colon X\to Y$ a p-morphism, and $\alpha^{-1} \colon {\sf Up}(Y)\to{\sf Up}(X)$ the corresponding complete Heyting homomorphism. Then $\alpha^{-1}$ is one-to-one iff $\alpha$ is onto, and $\alpha^{-1}$ is onto iff $\alpha$ is one-to-one. 
\end{enumerate}
\end{Theorem}

A closed upset of an Esakia space is an Esakia space (see, e.g., \cite[Lem.~3.4.11]{Esakia-book85}). Since one-to-one (continuous) p-morphisms correspond to (closed) upsets, we obtain the following characterization of quotients. For part (1) see \cite[Thm.~3.4.16]{Esakia-book85}, and for part (2) see \cite[Thms.~3.4, 3.5]{MR0197372}. 

\begin{Corollary} \label{cor: hom images}
The following conditions hold.
\begin{enumerate}
\item For an Esakia space $X$, the map $U \mapsto U^\ast$ is a bijection between the closed upsets of $X$ and the quotients of $X^\ast$.
\item For a poset $X$, the map $U \mapsto \mathsf{Up}(U)$ is a bijection between the upsets of $X$ and the complete quotients of ${\sf Up}(X)$.
\end{enumerate}
\end{Corollary} 

We next describe the kernels of onto (continuous) p-morphisms. To this end, given a binary relation $R$ on a set $X$ and $U \subseteq X$, we let
\[
R(U) = \{ x \in X : \langle y, x \rangle \in R  \text{ for some }y \in U \}.
\]
If $R$ is an equivalence relation, then $R(U)=U$ iff $U$ is a union of equivalence classes of $R$. In such a case, we say that $U$ is $R$-\emph{saturated}. 

\begin{Definition} \label{def: E partition}
\
\begin{enumerate}
\item Let $X$ be an Esakia space. An \emph{Esakia partition} (or {\em E-partition} for short) of $X$ is an equivalence relation $R$ on $X$ satisfying the following conditions:
\begin{enumerate}
\item\label{item:E-partition a} If $\langle x, y \rangle \in R$ and $x \leq z$, then there is $u \in X$ such that $y \leq u$ and $\langle z, u \rangle \in R$;
\item\label{item:E-partition b} If $\langle x, y \rangle \notin R$, then there is an $R$-saturated clopen upset $U$ such that $x\in U$ and $y\notin U$. 
\end{enumerate}
\item Let $X$ be a poset. An {\em E-partition} of $X$ is an equivalence relation $R$ on $X$ satisfying Condition~(\ref{item:E-partition a}) and the following version of Condition~(\ref{item:E-partition b}):
\begin{enumerate}
\item[(b')]\label{item:E-partition2} If $\langle x, y \rangle \notin R$, then there is an $R$-saturated upset $U$ such that $x \in U$ and $y \notin U$.
\end{enumerate}
\end{enumerate}
\end{Definition}

Let $X$ be an Esakia space or a poset. If $R$ is an E-partition of $X$, we define a partial order $\leq_R$ on $X / R$ as follows for every $x, y \in X$:
\[
[ x ] \leq_R [y] \Longleftrightarrow \text{ there are } x'\in [x] \text{ and } y' \in [y] \text{ such that } x' \leq y'.
\]
Since $R$ is an E-partition, the partial order $\leq_R$ is well defined and the map $x \mapsto [x]$ is a p-morphism from $X$ to $X /R$. Furthermore, when $X$ is an Esakia space, the poset $X/R$ endowed with the quotient topology (i.e., the open sets of $X/R$ are the $R$-saturated open sets of $X$) is an Esakia space and the map $x \mapsto [x]$ is a continuous p-morphism.

A subalgebra $\A$ of a complete Heyting algebra $\B$ is called \emph{complete} when $\A$ is also a complete sublattice of $\B$. Since E-partitions are exactly the kernels of (continuous) p-morphisms, from Theorem~\ref{thm: 1-1 and onto} we deduce:

\begin{Corollary} \label{cor: subs}
\
\begin{enumerate}
\item For an Esakia space $X$, the map $R \mapsto X / R$ is a bijection between the E-partitions of $X$ and the subalgebras of $X^*$.
\item For a poset $X$, the map $R \mapsto X / R$ is a bejection between the E-partitions of $X$ and the complete subalgebras of $\mathsf{Up}(X)$.
\end{enumerate}
\end{Corollary}

\noindent For part (1) of the above result see \cite[Cor.\ 2.3.1]{Bez-PhD}, and for part (2) see \cite[Thm.~4.6]{MR0197372}.

\section{Degrees of the finite model property} \label{sec: jankov}

We denote the set of posets validating an si-logic $\mathsf{L}$ by $\mathsf{Fr}(\mathsf{L})$. The \textit{degree of incompleteness} of $\sf L$ is the number of si-logics $\sf L'$ such that $\mathsf{Fr}(\mathsf{L}) = \mathsf{Fr}(\mathsf{L}')$. In this paper we are concerned with the degree of fmp. Thus, we restrict our attention to finite posets and let $\mathsf{Fin}(\mathsf{L})$ be the set of finite members of $\mathsf{Fr}(\mathsf{L})$.

\begin{Definition}
Let $\sf L$ be an si-logic.
\begin{enumerate}
\item The \textit{fmp span} $\textup{fmp}({\sf L})$ of $\sf L$ is the set of si-logics $\sf L'$ such that $\sf Fin(L')=Fin(L)$.
\item The \textit{degree of fmp} $\textup{deg}({\sf L})$ of $\sf L$ is the cardinality of $\textup{fmp}({\sf L})$.
\end{enumerate}
\end{Definition}

We call a poset $X$ {\em rooted} if there is $x\in X$ such that $X={\uparrow}x$. Such an $x$ is clearly unique and we call it the {\em root} of $X$. Given an si-logic $\mathsf{L}$, we denote the class of the rooted members of $\mathsf{Fin(L)}$ by $\mathsf{RFin(L)}$. Notice that, for each pair $\mathsf{L}$ and $\mathsf{L}'$ of si-logics, we have
\[
\mathsf{Fin}(\mathsf{L}) = \mathsf{Fin}(\mathsf{L}')\, \, \text{ iff } \, \, \mathsf{RFin}(\mathsf{L}) = \mathsf{RFin}(\mathsf{L}').
\]
To see this, suppose that $\mathsf{Fin}(\mathsf{L}) \neq  \mathsf{Fin}(\mathsf{L}')$. By symmetry, we may assume that there is a finite poset $X$ validating $\mathsf{L}$ and refuting $\mathsf{L}'$. Then there is $x \in X$ such that ${\uparrow}x$ validates $\mathsf{L}$ and refutes 
$\mathsf{L}'$. Consequently, ${\uparrow}x$ is a member of $\mathsf{RFin}(\mathsf{L}) \smallsetminus  \mathsf{RFin}(\mathsf{L}')$, as desired. In view of this,
the fmp span of an si-logic $\mathsf{L}$ is the set of si-logics $\mathsf{L}'$ such that $\mathsf{RFin}(\mathsf{L}) = \mathsf{RFin}(\mathsf{L}')$. We will use this fact without further notice.

Since each si-logic $\sf L$ belongs to its own fmp span and there are exactly $2^{\aleph_0}$ si-logics, the obvious lower and upper bounds for $\textup{deg}({\sf L})$ are $1$ and $2^{\aleph_{0}}$. The main result of this paper is the Antidichotomy Theorem stating that these restrictions are indeed optimal in that each nonzero  cardinal $\kappa$ such that $\kappa \leq \aleph_0$ or $\kappa = 2^{\aleph_0}$ occurs as the degree of fmp of some si-logic. Thus, under CH, every   cardinal  $1\leq\kappa\leq 2^{\aleph_0}$ occurs as the degree of fmp of some si-logic.\footnote{It is not known whether it is consistent with ZFC that there are si-logics with the degree of fmp $\kappa$ for $\aleph_0<\kappa < 2^{\aleph_0}$ (see Problem 1 in the Conclusions).}\ 
More precisely, we will prove the following:

\begin{Theorem}[\textbf{Antidichotomy Theorem}]\label{Thm:dichotomy:IPC:main}
For each nonzero cardinal $\kappa$ such that $\kappa \leq \aleph_0$ or $\kappa = 2^{\aleph_0}$ there is an si-logic $\sf L$ such that $\textup{deg}({\sf L}) = \kappa$.
\end{Theorem}

As we pointed out in the introduction, one of the techniques required to prove this theorem is that of splittings and Jankov formulas. We recall that a pair of elements $(a,b)$ of a lattice $L$ {\em splits} $L$ if $L$ is the disjoint union of ${\uparrow}a$ and ${\downarrow}b$ \cite[Sec.~9.4]{ChZa97}. An si-logic $\sf L$ is a {\em splitting logic} if there is an si-logic $\sf M$ such that the pair $\sf (L,M)$ splits the lattice  $\textup{Ext}\,\mathsf{IPC}$. An si-logic is {\em join-splitting} if it is the join in $\textup{Ext}\,\mathsf{IPC}$ of a set of splitting si-logics.

Jankov \cite{Jankov63for} provided an axiomatization of the join-splitting si-logics. We recall that a Heyting algebra $\A$ is \textit{subdirectly irreducible} (SI for short) if it has the second largest element (equivalently, the filter $\{1\}$ is completely meet-irreducible in the lattice of filters of $\A$). 
By the Jankov Theorem \cite{Jankov63for}, with each finite SI Heyting algebra $\A$ we can associate a formula $\mathcal{J}(\A)$ (referred to as the \textit{Jankov formula} of $\A$) that axiomatizes the least si-logic $\sf L$ such that $\A \nvDash \sf L$: 

\begin{Theorem}[\textbf{Jankov Theorem}]
An si-logic $\sf L$ is a splitting logic iff there exists a finite SI Heyting algebra $\A$ such that ${\sf L}={\sf IPC}+\mathcal{J}(\A)$. Consequently, $\sf L$ is a  join-splitting logic iff $\sf L$ is axiomatizable by Jankov formulas.
\end{Theorem}

The following lemma governs the behavior of Jankov formulas \cite{Jankov69}: 

\begin{Lemma}[\textbf{Jankov Lemma}]\label{lem: Jankov lemma}
Let $\A$ and $\B$ be Heyting algebras with $\A$ finite and SI. Then $\B\nvDash\mathcal{J}(\A)$ iff $\A$  is a subalgebra of a homomorphic image of $\B$.
\end{Lemma}

It is well known that a Heyting algebra $\A$ is SI iff 
$\A_{\ast}$ has a root which, moreover, is isolated (see, e.g.,  \cite[Appendix 1.1]{Esakia-book85}). Therefore, the Finite Esakia Duality implies that the finite SI Heyting algebras are those of the form ${\sf Up}(X)$ where $X$ is a finite rooted poset. Because of this, given a finite rooted poset $X$, we denote by $\mathcal{J}(X)$ the Jankov formula of the finite SI Heyting algebra ${\sf Up}(X)$. 
 Thus, in view of Theorem~\ref{thm: 1-1 and onto} and Corollary~\ref{cor: hom images}, the Jankov Lemma can be formulated dually as follows:
 
\begin{Lemma}[\textbf{Dual Jankov Lemma}]\label{lem: dual Jankov lemma}
Let $X$ be a finite rooted poset. For every Esakia space $Y$ we have $Y \nvDash\mathcal{J}(X)$ iff $X$ is a continuous p-morphic image of a closed upset of $Y$.
\end{Lemma}

\begin{Remark}
In \cite{ChZa97}, (continuous) p-morphisms are referred to as \emph{reductions}. Using this terminology, the Dual Jankov Lemma can be formulated as follows: $Y \nvDash\mathcal{J}(X)$ iff a closed upset of $Y$ is reducible to $X$. 
\end{Remark}

 Notably, the following variant of the Dual Jankov Lemma for posets holds too \cite{Fin74b}:

\begin{Lemma}[\textbf{Fine Lemma}]\label{lem: Fine lemma}
Let $X$ be a finite rooted poset. For every poset $Y$ we have $Y \nvDash\mathcal{J}(X)$ iff $X$ is a p-morphic image of an upset of $Y$.
\end{Lemma}

 The next immediate consequence of the Dual Jankov Lemma governs the interaction between Jankov formulas and si-logics.

\begin{Corollary}\label{Cor:Jankov-easy}
For every finite rooted poset $X$ and si-logic $\mathsf{L}$ we have $X \vDash \mathsf{L}$ iff $\mathcal{J}(X) \notin \mathsf{L}$.
\end{Corollary}

We rely on the following folklore result. We provide a full proof of part (2) since we were not able to find one in the literature.\color{black}

\begin{Lemma}\label{lem: auxiliary} The following conditions holds.
\begin{enumerate}
\item\label{item: auxiliary} Let $X$ be a finite rooted poset and $K$ a class of Esakia spaces. Then $X\vDash \mathsf{Log}(K)$ iff there is $Y\in K$ such that $X$ is a continuous p-morphic image of a closed upset of $Y$.
\item\label{item: splittings} Two si-logics $\mathsf{L}$ and $\mathsf{L}'$ contain the same Jankov formulas iff $\mathsf{Fin}(\mathsf{L}) = \mathsf{Fin}(\mathsf{L}')$.
\end{enumerate}
\end{Lemma}

\begin{proof}
(\ref{item: auxiliary}) Immediate from the Dual Jankov Lemma.

(\ref{item: splittings}) First suppose that $\mathsf{Fin}(\mathsf{L}) \ne \mathsf{Fin}(\mathsf{L}')$. Since a poset validates a formula iff each of its principal upsets does, 
without loss of generality we may assume that there is a finite rooted $X \in \mathsf{Fin}(\mathsf{L}) \smallsetminus \mathsf{Fin}(\mathsf{L}')$. By Corollary \ref{Cor:Jankov-easy} we have $\mathcal{J}(X) \in \mathsf{L}' \smallsetminus \mathsf{L}$. Conversely, suppose that $\mathsf{L}$ and $\mathsf{L}'$ do not contain the same Jankov formulas. We may assume without loss of generality that $\mathcal{J}(X) \in \mathsf{L} \smallsetminus \mathsf{L}'$ for a finite rooted poset $X$. From Corollary \ref{Cor:Jankov-easy} it follows that $X \in \mathsf{Fin}(\mathsf{L}') \smallsetminus \mathsf{Fin}(\mathsf{L})$. 
\end{proof}

In order to describe fmp spans, it is convenient to introduce the following concept.

\begin{Definition}
For an si-logic $\mathsf{L}$, define 
\begin{enumerate}
\item $\mathsf{L}^{+} = \mathsf{Log}({\sf Fin}(\mathsf{L}))$; 
\item $\mathsf{L}^- = {\sf IPC}+\{\mathcal{J}(X) : X \notin{\sf Fin}(\mathsf{L})\}$.
\end{enumerate}
\end{Definition}

\noindent Let $[\mathsf{L}^-, \mathsf{L}^+]$ be the interval in the lattice $\textup{Ext}\,\mathsf{IPC}$.

\begin{Theorem}\label{Thm:interval}
For an si-logic $\mathsf{L}$ we have:
\begin{enumerate}
\item\label{item:interval:1} $\textup{fmp}(\mathsf{L})=[\mathsf{L}^-, \mathsf{L}^+]$.
\item $\mathsf{L}^+$ is the only member of $\textup{fmp}(\mathsf{L})$ that has the fmp.
\item $\mathsf{L}^-$ is the only member of $\textup{fmp}(\mathsf{L})$ that is axiomatizable by Jankov formulas.
\end{enumerate}
\end{Theorem}

\begin{proof}
(1) We begin by proving that $\mathsf{L}^- \in \textup{fmp}(\mathsf{L})$. By Lemma~\ref{lem: auxiliary}(\ref{item: splittings}), it suffices to show that $\mathsf{L}$ and $\mathsf{L}^-$ contain the same Jankov formulas. In view of Corollary \ref{Cor:Jankov-easy}, $\{\mathcal{J}(X) : X \notin{\sf Fin}(\mathsf{L})\}$ is the set of Jankov formulas in $\mathsf{L}$. Since $\mathsf{L}^- = \mathsf{IPC} + \{\mathcal{J}(X) : X \notin{\sf Fin}(\mathsf{L})\}$, every Jankov formula in $\mathsf{L}$ belongs to $\mathsf{L}^-$ and $\mathsf{L}^- \subseteq \mathsf{L}$. The latter implies that every Jankov formula in $\mathsf{L}^-$ belongs to $\mathsf{L}$. Thus, $\mathsf{L}^- \in \textup{fmp}(\mathsf{L})$, as desired. Since $\mathsf{L}^-$ is axiomatized by Jankov formulas, this implies that it is the least element of $\textup{fmp}(\mathsf{L})$.

We next prove that $\mathsf{L}^+$ is the greatest logic in $\textup{fmp}(\mathsf{L})$.  Clearly ${\sf Fin}(\mathsf{L})\subseteq{\sf Fin}(\mathsf{L}^+)$ by the definition of $\mathsf{L}^+$. The other inclusion follows from Lemma~\ref{lem: auxiliary}(\ref{item: auxiliary}) and the fact that ${\sf Fin}(\mathsf{L})$ is closed under the formation of upsets and p-morphic images. Thus, ${\sf Fin}(\mathsf{L}^+)={\sf Fin}(\mathsf{L})$, and so  
$\mathsf{L}^{+} \in \textup{fmp}(\mathsf{L})$. 
Let $\mathsf{L}' \in \textup{fmp}(\mathsf{L})$. Then ${\sf Fin}(\mathsf{L}')={\sf Fin}(\mathsf{L})$. Since $\mathsf{L}^+$ is the logic of ${\sf Fin}(\mathsf{L})$, we conclude that $\mathsf{L}'\subseteq\mathsf{L}^+$. Thus,
$\mathsf{L}^+$ is the greatest element of $\textup{fmp}(\mathsf{L})$.

It follows from the definition of $\textup{fmp}(\mathsf{L})$ that $\textup{fmp}(\mathsf{L})$ is an interval in the lattice of si-logics. Together with the fact that $\mathsf{L}^-$ and $\mathsf{L}^+$ are
the least and greatest elements of $\textup{fmp}(\mathsf{L})$, this implies that $\textup{fmp}(\mathsf{L})=[\mathsf{L}^-, \mathsf{L}^+]$.

(2) By definition, 
$\mathsf{L}^+$ has the fmp.
If $\mathsf{L}' \in \textup{fmp}(\mathsf{L})$ has the fmp, then $\mathsf{L}'$ is the logic of ${\sf Fin}(\mathsf{L}')$. But ${\sf Fin}(\mathsf{L}')={\sf Fin}(\mathsf{L})$, so $\mathsf{L}=\mathsf{L}^+$. Thus,
$\mathsf{L}^+$ is the only member of $\textup{fmp}(\mathsf{L})$ with the fmp.

(3) By definition, $\mathsf{L}^-$ is axiomatized by Jankov formulas. Let $\mathsf{L}' \in \textup{fmp}(\mathsf{L})$ be also axiomatized by Jankov formulas. Since $\mathsf{Fin}(\mathsf{L}') = \mathsf{Fin}(\mathsf{L}) = \mathsf{Fin}(\mathsf{L}^-)$, we can apply Lemma \ref{lem: auxiliary}(\ref{item: splittings}) to obtain that $\mathsf{L}'$  and $\mathsf{L}^-$ contain the same Jankov formulas. As both $\mathsf{L}^-$ and $\mathsf{L}'$ are axiomatized by Jankov formulas, we conclude that $\mathsf{L}^- = \mathsf{L}'$. Thus, $\mathsf{L}^-$ is the only member of $\textup{fmp}(\mathsf{L})$ axiomatizable by Jankov formulas.
\end{proof}

As a consequence, we obtain a transparent description of the si-logics whose degree of fmp is~1.

\begin{Corollary}\label{Cor:degre-one}
An si-logic $\mathsf{L}$ has the degree of fmp 1 iff it has the fmp and is axiomatizable by Jankov formulas.
\end{Corollary}

\begin{proof}
First suppose that $\textup{deg}(\mathsf{L}) = 1$. 
Since $\mathsf{L}, \mathsf{L}^-,\mathsf{L}^+ \in \textup{fmp}(\mathsf{L})$, this implies that $\mathsf{L} = \mathsf{L}^- = \mathsf{L}^+$. Because $\mathsf{L}^+$ has the fmp and $\mathsf{L}^-$ is axiomatizable by Jankov formulas, we conclude that $\mathsf{L}$ has the fmp and is axiomatizable by Jankov formulas.

To prove the converse, suppose that $\mathsf{L}$ has the fmp and is axiomatizable by Jankov formulas. 
By Theorem \ref{Thm:interval}, the only member of $\textup{fmp}(\mathsf{L})$ with the fmp is $\mathsf{L}^+$, and the only member of $\textup{fmp}(\mathsf{L})$ that is axiomatizable by Jankov formulas is $\mathsf{L}^-$. Since $\mathsf{L} \in \textup{fmp}(\mathsf{L})$, we obtain that $\mathsf{L} = \mathsf{L}^- = \mathsf{L}^+$. 
Therefore, with an application of Theorem \ref{Thm:interval}(\ref{item:interval:1}) we conclude that
\[
\textup{fmp}(\mathsf{L}) = [\mathsf{L}^-, \mathsf{L}^+] = [\mathsf{L}, \mathsf{L}] =  
\{ \mathsf{L}\},
\]
and hence $\textup{deg}(\mathsf{L}) = 1$.
\end{proof}

Examples of si-logics with the degree of fmp 1 include locally tabular logics. We recall that an si-logic $\mathsf{L}$ is {\em locally tabular} if for each $n<\aleph_0$ the Lindenbaum-Tarski algebra of $\mathsf{L}$ in $n$ variables is finite. Clearly each locally tabular logic has the fmp. Moreover, each locally tabular si-logic is axiomatizable by Jankov formulas (see, e.g., \cite[Thm.~3.4.24]{Bez-PhD}). Thus, we obtain:

\begin{Corollary}\label{Cor:locally-finite}
The degree of fmp of locally tabular si-logics is $1$. \color{black}
\end{Corollary}

Since there are continuum many locally tabular si-logics, the above corollary implies that there are also continuum many si-logics whose degree of fmp is 1. We point out that there are si-logics that are not locally tabular and yet have the degree of fmp 1. For example, $\sf IPC$ is such a logic. More examples will be given in Example~\ref{exa: KG and RN}.

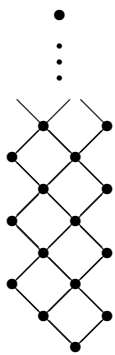
\begin{figure}
\begin{center}
\begin{picture}(62,144)
\put(38,3){\circle*{4}}
\put(38,3){\line(1,1){12}}
\put(38,3){\line(-1,1){12}}
\put(50,15){\circle*{4}}
\put(50,15){\line(-1,1){12}}
\put(26,15){\circle*{4}}
\put(26,15){\line(1,1){12}}
\put(38,27){\circle*{4}}
\put(12,5){}
\put(50,6){}
\put(42,25){}

\put(26,15){\line(-1,1){12}}
\put(38,27){\line(-1,1){12}}
\put(14,27){\circle*{4}}
\put(14,27){\line(1,1){12}}
\put(26,39){\circle*{4}}
\put(0,17){}
\put(30,37){}

\put(38,27){\line(1,1){12}}
\put(50,39){\circle*{4}}
\put(50,39){\line(-1,1){12}}
\put(26,39){\line(1,1){12}}
\put(38,51){\circle*{4}}
\put(50,30){}
\put(42,49){}

\put(26,39){\line(-1,1){12}}
\put(38,51){\line(-1,1){12}}
\put(14,51){\circle*{4}}
\put(14,51){\line(1,1){12}}
\put(26,63){\circle*{4}}
\put(0,41){}
\put(30,61){}

\put(38,51){\line(1,1){12}}
\put(50,63){\circle*{4}}
\put(50,63){\line(-1,1){12}}
\put(26,63){\line(1,1){12}}
\put(38,75){\circle*{4}}
\put(50,54){}
\put(42,73){}

\put(26,63){\line(-1,1){12}}
\put(38,75){\line(-1,1){12}}
\put(14,75){\circle*{4}}
\put(14,75){\line(1,1){12}}
\put(26,87){\circle*{4}}
\put(0,65){}
\put(30,85){}

\put(38,75){\line(1,1){12}}
\put(50,87){\circle*{4}}
\put(50,87){\line(-1,1){10}}
\put(26,87){\line(1,1){10}}
\put(50,78){}

\put(26,87){\line(-1,1){10}}

\put(32,105){\circle*{2}}
\put(32,111){\circle*{2}}
\put(32,117){\circle*{2}}

\put(32,129){\circle*{4}}
\put(36,127){}
\end{picture}
\end{center}
\caption{The Rieger-Nishimura lattice.}
\label{Fig:RN-algebra}
\end{figure}

\section{The Kuznetsov-Gerciu logic} \label{sec: KG}

In this section we briefly review the si-logic of Kuznetsov and Ger\v{c}iu \cite{GeKuz70,KuzGer70a}. 
We start by recalling (see \cite{Ri49,Ni60}) that the one-generated free Heyting algebra, known as the \textit{Rieger-Nishimura lattice} $\boldsymbol{RN}$, is the Heyting algebra depicted in Figure \ref{Fig:RN-algebra}.  

Let $\A$ and $\B$ be Heyting algebras. The \emph{sum} $\A + \B$ is the Heyting algebra obtained by pasting $\A$ \emph{below} $\B$ and gluing the top element of $\A$ to the bottom element of $\B$ \cite[Appendix A.9]{Esakia-book85}. As $+$ is clearly associative, there is no ambiguity in writing $\A_{1} + \dots + \A_{n}$ for finitely many Heyting algebras $\A_{1}, \dots, \A_{n}$, each glued to the next.

\begin{Definition}
The \emph{Kuznetsov-Ger\v{c}iu} logic $\mathsf{KG}$ is the si-logic 
of all Heyting algebras of the form $\A_{1} + \dots + \A_{n}$ where $\A_{1}, \dots, \A_{n}$ are one-generated. 
\end{Definition}

We will utilize that $\mathsf{KG}$ is a subframe logic. We recall that the theory of subframe modal logics was developed by Fine \cite{Fine85}, 
and that Zakharyaschev \cite{Zakha89} studied subframe si-logics. 
For the present purpose, we concentrate on subframe si-logics. 

With each finite rooted poset $X$ we can associate a formula $\beta(X)$ in the language of $\mathsf{IPC}$,  called the \emph{subframe formula} of $X$. 
Bearing in mind that \textit{frame} and \textit{poset} are synonyms in the context of si-logics, the next result 
motivates this terminology. 

\begin{Theorem}[Fine \& Zakharyaschev]\label{Thm:subframe-formulas}
Let $X$ be a finite rooted poset.  
\begin{enumerate}
\item \label{subframe-formulas 1} For every Esakia space $Y$ we have $Y \nvDash \beta(X)$ iff $X$ is a continuous p-morphic image of some clopen $Z\subseteq X$.
\item \label{subframe-formulas 2} For every poset $Y$ we have $Y \nvDash \beta(X)$ iff $X$ is a p-morphic image of some $Z\subseteq Y$.
\end{enumerate}
\end{Theorem}

\begin{proof}
See \cite[Thm.\ 9.40(ii)]{ChZa97}. Our formulation of the result differs from that of \cite[Thm.\ 9.40(ii)]{ChZa97} in that we require $Z$ to be  clopen (as opposed to a subframe). The fact that this is harmless follows from the proof of the result (see also \cite[Thm.~3.3.16]{Bez-PhD}).
\end{proof}

An si-logic is a {\em subframe logic} if it is axiomatizable  by subframe formulas.

\begin{Theorem}[Fine \& Zakharyaschev] \label{thm: fmp for subframe}
Each subframe si-logic has the fmp.
\end{Theorem}

\begin{proof}
See, e.g., \cite[Thm.~11.20]{ChZa97}.
\end{proof}

\begin{figure}
\[
\xymatrix@R=18pt @C=18pt @!0{
&&&
&&
& && *-{\bullet}\ar@{-}[d]
\\
&&&
*-{\bullet}\ar@{-}[d] && *-{\bullet}\ar@{-}[d]
& && *-{\bullet}\ar@{-}[d]
\\
*-{\bullet}\ar@{-}[dr]&*-{\bullet}\ar@{-}[d]&*-{\bullet}\ar@{-}[dl]&
*-{\bullet}\ar@{-}[dr] && *-{\bullet}\ar@{-}[dl]
&*-{\bullet}\ar@{-}[dr] && *-{\bullet}\ar@{-}[dl]
\\
&*-{\bullet}&&
&*-{\bullet}&
&&*-{\bullet}&
\\
&P_{1}&&
&P_{2}&
&&P_{3}&
}
\]
\caption{The posets $P_1, P_2$, and $P_3$.}
\label{Fig:KG-axiom}
\end{figure}
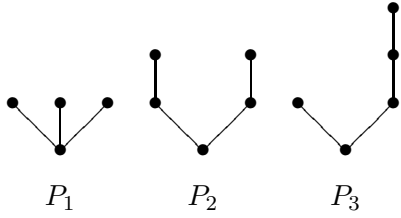

As we pointed out earlier in the section, $\sf KG$ is axiomatizable by subframe formulas (see, e.g., \cite{Kracht93MLQ} or \cite[Thm.~4.3.4]{Bez-PhD}):

\begin{Theorem}\label{Thm:axioms-subframe-KG}
$\mathsf{KG}$ is 
 axiomatized by the subframe formulas of the posets in Figure  \ref{Fig:KG-axiom}.
\end{Theorem}

For our purposes it is crucial that $\mathsf{KG}$ is also axiomatizable by Jankov formulas. For this we first recall the notion of width for posets. 

\begin{Definition} \label{def: width}
Let $1\leq n<\aleph_0$. 
The \textit{width} of a rooted poset $X$ is $n$ if
\[
n = \max \{ \kappa : \kappa \text{ is the cardinality of an antichain of }X \}.
\]
The \textit{width} of a poset $X$ is $n$ if all principal upsets of $X$ have width $\leq n$ and there is a principal upset of width $n$. 
 The empty poset will be assumed to have width zero. 
\end{Definition}

We next define the notion of width for Heyting algebras.

\begin{Definition}
Let $n<\aleph_0$. A Heyting algebra $\A$ has \textit{width} $n$ if $\A_\ast$ has width $n$. Let
\[
\mathsf{W}_n = \{ \A \in \mathsf{HA} : \A \text{ has width }\leq n \}.
\]
\end{Definition}

\begin{Definition}
For $n<\aleph_0$ let 
\[
{\sf bw}_n = \bigvee_{i=0}^n(p_i\to\bigvee_{j\ne i}p_j)
\] 
and define 
\[
{\sf BW}_n = {\sf IPC}+{\sf bw}_n.
\]
\end{Definition}

Sobolev \cite{So77} proved that a Heyting algebra $\A$ validates ${\sf bw}_n$ iff $\A\in\mathsf{W}_n$. Thus, the members of $\mathsf{W}_n$ are exactly the algebraic models of ${\sf BW}_n$. 

\begin{figure}
\[
\begin{tabular}{ccccccccc}
\begin{tikzpicture}
     \tikzstyle{point} = [shape=circle, thick, draw=black, fill=black , scale=0.35]
    \node[label=below:{$F_{n+1}$}]  (0) at (-1.5,-0.2)  {};    
    \node[label=above:{$1$}]  (v1) at (-2.5,1) [point] {};
    \node (0) at (-1.5,0) [point] {};
    \node[label=above:{$3$}]  (v3) at (-0.5,1) [point] {};
    \node[label=above:{$2$}] (v2) at (-1.5,1) [point] {};
   \node (down-dots) at (0.32,1) {$\cdots\cdots$};
    \node[label=above:{$n+1$}]  (v4) at (1.1,1) [point] {};

    \draw   (v1) -- (0) -- (v3) (v4) -- (0) -- (v2);
\end{tikzpicture}

\end{tabular}
\]\caption{The poset $F_{n+1}$.}
\label{Fig:Wn-subframe}
\end{figure}
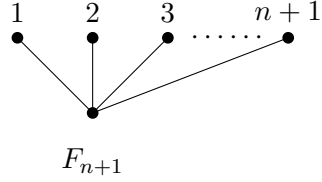

\begin{Theorem}[Fine \& Zakharyaschev]\label{thm: BWn} \label{Thm:subframe-axioms-for-Wn}
Each si-logic ${\sf BW}_n$ is axiomatized by the subframe formula of the poset depicted in Figure~\ref{Fig:Wn-subframe}.
Thus, each ${\sf BW}_n$ has the fmp.
\end{Theorem}

\begin{proof}
See \cite{Fine85,Zakha89}.
\end{proof}

We will use the following result of Kracht \cite[Prop.~23]{Kr93b}. 

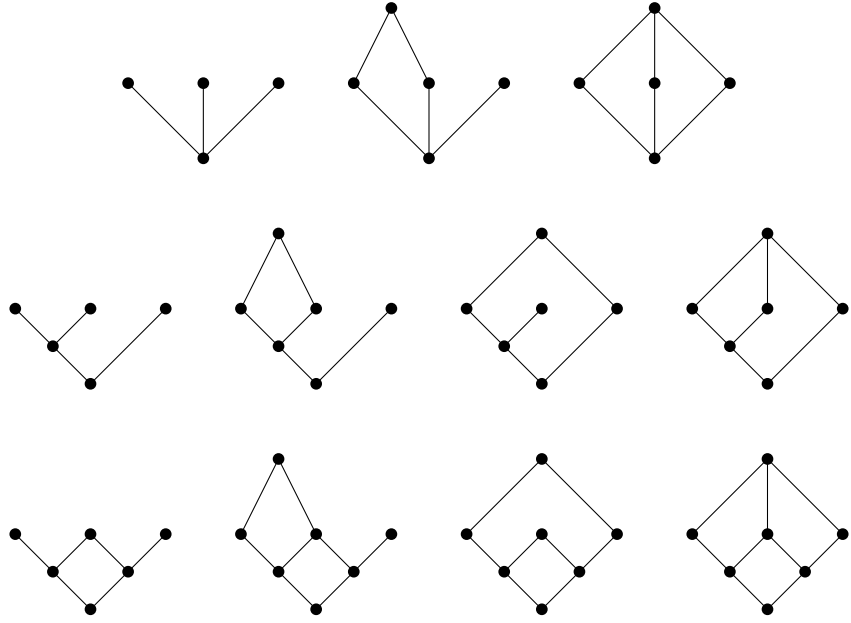
\begin{figure}[h]
\begin{tabular}{ccccccc}
\begin{tikzpicture}
    \tikzstyle{point} = [shape=circle, thick, draw=black, fill=black , scale=0.35]
   
    \node (a0) at (-3,0) [point] {};
    \node  (a1) at (-4,1) [point] {};
    \node (a2) at (-3,1) [point] {};
    \node (a3) at (-2,1) [point] {};
    \draw (a2) -- (a0) -- (a1)  ;
    \draw (a0) -- (a3) ;
    
    \node (b0) at (0,0) [point] {};
    \node  (b1) at (-1,1) [point] {};
    \node (b2) at (-0,1) [point] {};
    \node (b3) at (1,1) [point] {};
        \node (b4) at (-0.5,2) [point] {};

    \draw (b4) --  (b2) -- (b0) -- (b1) -- (b4) ;
    \draw (b0) -- (b3) ;
    
    \node (c0) at (3,0) [point] {};
    \node  (c1) at (4,1) [point] {};
    \node (c2) at (3,1) [point] {};
    \node (c3) at (2,1) [point] {};
        \node (c4) at (3,2) [point] {};

    \draw (c4) --  (c2) -- (c0) -- (c1) -- (c4) ;
    \draw (c0) -- (c3)-- (c4) ;
    
    \node (d0) at (-4.5,-3) [point] {};
    \node  (d1) at (-5.5,-2) [point] {};
    \node (d2) at (-4.5,-2) [point] {};
    \node (d3) at (-3.5,-2) [point] {};
        \node (d4) at (-5,-2.5) [point] {};
    \draw  (d4) -- (d2)  (d0) -- (d1)  ;
    \draw (d0) -- (d3)  ;
    
    \node (e0) at (-1.5,-3) [point] {};
    \node  (e1) at (-2.5,-2) [point] {};
    \node (e2) at (-1.5,-2) [point] {};
    \node (e3) at (-0.5,-2) [point] {};
        \node (e4) at (-2,-2.5) [point] {};
                \node (e5) at (-2,-1) [point] {};

    \draw  (e4) -- (e2) -- (e5) (e0) -- (e1) -- (e5)  ;
    \draw (e0) -- (e3)  ;
    
    \node (f0) at (1.5,-3) [point] {};
    \node  (f1) at (0.5,-2) [point] {};
    \node (f2) at (1.5,-2) [point] {};
    \node (f3) at (2.5,-2) [point] {};
        \node (f4) at (1,-2.5) [point] {};
                \node (f5) at (1.5,-1) [point] {};

    \draw  (f4) -- (f2)  (f0) -- (f1)  -- (f5);
    \draw (f0) -- (f3)-- (f5)   ;
    
    \node (g0) at (4.5,-3) [point] {};
    \node  (g1) at (3.5,-2) [point] {};
    \node (g2) at (4.5,-2) [point] {};
    \node (g3) at (5.5,-2) [point] {};
        \node (g4) at (4,-2.5) [point] {};
                \node (g5) at (4.5,-1) [point] {};

    \draw  (g4) -- (g2)  (g0) -- (g1)  -- (g5);
    \draw (g0) -- (g3)-- (g5) -- (g2)   ;
    
    
     \node (h0) at (-4.5,-6) [point] {};
    \node  (h1) at (-5.5,-5) [point] {};
    \node (h2) at (-4.5,-5) [point] {};
    \node (h3) at (-3.5,-5) [point] {};
        \node (h4) at (-5,-5.5) [point] {};
                \node (h9) at (-4,-5.5) [point] {};

    \draw  (h4) -- (h2) -- (h9) (h0) -- (h1)  ;
    \draw (h0) -- (h3)   ;
         
    \node (i0) at (-1.5,-6) [point] {};
    \node  (i1) at (-2.5,-5) [point] {};
    \node (i2) at (-1.5,-5) [point] {};
    \node (i3) at (-0.5,-5) [point] {};
        \node (i4) at (-2,-5.5) [point] {};
                \node (i5) at (-2,-4) [point] {};
                                \node (i9) at (-1,-5.5) [point] {};

    \draw  (i4) -- (i2) -- (i5) (i0) -- (i1) -- (i5)  ;
    \draw (i0) -- (i3)  ;
    
        \draw (i2) -- (i9)  ;
    
    \node (l0) at (1.5,-6) [point] {};
    \node  (l1) at (0.5,-5) [point] {};
    \node (l2) at (1.5,-5) [point] {};
    \node (l3) at (2.5,-5) [point] {};
        \node (l4) at (1,-5.5) [point] {};
                \node (l5) at (1.5,-4) [point] {};
                                \node (l9) at (2,-5.5) [point] {};

    \draw  (l4) -- (l2) -- (l9) (l0) -- (l1)  -- (l5);
    \draw (l0) -- (l3)-- (l5)   ;
        
    \node (m0) at (4.5,-6) [point] {};
    \node  (m1) at (3.5,-5) [point] {};
    \node (m2) at (4.5,-5) [point] {};
    \node (m3) at (5.5,-5) [point] {};
        \node (m4) at (4,-5.5) [point] {};
                \node (m5) at (4.5,-4) [point] {};
                                \node (m9) at (5,-5.5) [point] {};

    \draw  (m4) -- (m2) -- (m9) (m0) -- (m1)  -- (m5);
    \draw (m0) -- (m3)-- (m5) -- (m2)   ;

\end{tikzpicture}

\end{tabular}
\caption{The eleven posets whose Jankov formulas axiomatize $\mathsf{W}_2$.}
\label{fig:width2}
\end{figure}

\begin{Theorem}[Kracht]\label{Thm:width2-kracht}
The logic ${\sf BW}_2$ is axiomatized by the Jankov formulas of the posets in Figure~\ref{fig:width2}.
\end{Theorem}

\begin{proof}
Kracht proved this result in the setting of normal modal logics extending $\sf S4$. A natural adaptation of the proof yields the analogous result for si-logics.
\end{proof}

We use Theorem~\ref{Thm:width2-kracht} to prove that $\mathsf{KG}$ is also axiomatizable by Jankov formulas: 
 
\begin{Theorem}\label{Cor:KG-Jankov-axiom}
$\mathsf{KG}$ is axiomatizable by Jankov formulas.
\end{Theorem}

\begin{proof}
Since a similar result was sketched by Kracht in  \cite[Sec.~D]{Kr93b} (again in the setting of normal modal logics extending $\sf S4$) and because full proofs require lengthy combinatorial arguments, they are moved to the Appendix.
\end{proof}

\begin{exa} \label{exa: KG and RN}
The above theorem provides further examples of si-logics that are not locally tabular, but have the degree of fmp 1. Let $\sf RN$ be the logic of the Rieger-Nishimura lattice 
$\boldsymbol{RN}$. It is well known that both $\mathsf{KG}$ and $\mathsf{RN}$ have the fmp: for $\mathsf{KG}$ this follows from Theorems~\ref{thm: fmp for subframe} and \ref{Thm:axioms-subframe-KG}, while for $\mathsf{RN}$ see, e.g., \cite[Thm.\ 5.35]{BBdeJ08}. Moreover, 
$\sf RN$ is axiomatizable relative to $\mathsf{KG}$ by Jankov formulas \cite[Thm.~4.33]{BBdeJ08}. Therefore,
by Theorem \ref{Cor:KG-Jankov-axiom}, 
both $\sf KG$ and $\sf RN$ 
are axiomatizable by Jankov formulas. Thus, by Corollary \ref{Cor:degre-one}, 
both logics have the degree of fmp 1.
Clearly neither logic is locally tabular since
$\boldsymbol{RN}\vDash\sf RN,KG$. 
\end{exa}

\section{The countable case or ``anything goes''}

In this section we establish the countable case of the Antidichotomy Theorem. We do this by exhibiting for each cardinal $1\leq \kappa \leq \aleph_0$, an si-logic $\mathsf{L}$ such that $\textup{deg}(\mathsf{L}) = \kappa$. As we will see below, $\mathsf{L}$ can be chosen to be an extension of $\mathsf{KG}$. More precisely, we will prove the following:

\begin{Theorem}\label{Thm:countable-part}
For each cardinal $1\leq \kappa\leq \aleph_0$ there exists an si-logic $\mathsf{L}\supseteq\mathsf{KG}$  such that $\textup{deg}(\mathsf{L}) = \kappa$.
\end{Theorem}

\begin{figure}
\[
\begin{tabular}{ccccccc}
\begin{tikzpicture}
    \tikzstyle{point} = [shape=circle, thick, draw=black, fill=black , scale=0.35]
   
\node[label=left:{$w_0$}] (w0) at (11.2,2.1) [point] {};
\node[label=right:{$w_1$}] (w1) at (12.6,2.1) [point] {};
\node[label=left:{$w_2$}] (w2) at (11.2,1.4) [point] {};
\node[label=right:{$w_3$}] (w3) at (12.6,1.4) [point] {};
\node[label=left:{$w_4$}] (w4) at (11.2,0.7) [point] {};
\node[label=right:{$w_5$}] (w5) at (12.6,0.7) [point] {};
\node[label=left:{$w_6$}] (w6) at (11.2,0) [point] {};
\node[label=right:{$w_7$}] (w7) at (12.6,0) [point] {};
\node[label=left:{$w_8$}] (w8) at (11.2,-0.7) [point] {};
\node[label=right:{$w_9$}] (w9) at (12.6,-0.7) [point] {};
\node[label=left:{$w_{10}$}] (w10) at (11.2,-1.4) [point] {};
\node[label=right:{$w_{11}$}] (w11) at (12.6,-1.4) [point] {};

\node (w12) at (11.2,-2.1) [] {};
\node (w13) at (12.6,-2.1) [] {};

\node[label=below:{$\omega$}] (omega) at (11.9,-2.8) [point] {};

 \draw[dotted] (w12) -- (omega) -- (w13);

 \draw[dotted] (w10) --(w12)  (w11)--(w13);
  \draw (w0) --  (w2) --  (w4) --  (w6) --  (w8) -- (w10) ;

    \draw (w1) --  (w3) --  (w5) -- (w7) -- (w9) -- (w11)  ;
  
  \draw (w0) -- (w3) (w2) -- (w5) (w4) -- (w7) (w6) -- (w9) (w8) -- (w11) ;
  
    \draw (w1) -- (w4) (w3) -- (w6) (w5) -- (w8) (w7) -- (w10)  ;

\end{tikzpicture}
\end{tabular}
\]
\caption{The Rieger-Nishimura ladder $\mathfrak{L}$.}
\label{fig:Rieger-Nishimura-ladder}
\end{figure}
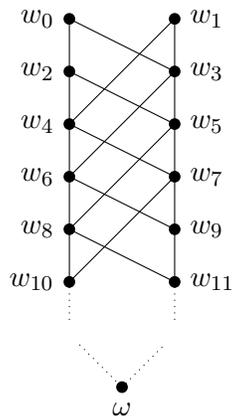

We will rely on several known facts about the Reieger-Nishimura lattice. 
We will use \cite{Bez-PhD} as our main reference, but these results can also be found in \cite{BBdeJ08}. 
The Esakia dual $\mathfrak{L}$ of the Rieger-Nishimura lattice $\boldsymbol{RN}$, often called the \textit{Rieger-Nishimura ladder}, is depicted in Figure \ref{fig:Rieger-Nishimura-ladder},  where the topology can be described as follows: a subset of $\mathfrak{L}$ is open iff either it misses $\omega$ or it is cofinite. In other words, each $w_n$ is an isolated point and $\omega$ is the only limit point. 

Using the labeling of Figure \ref{fig:Rieger-Nishimura-ladder}, for each $n\geq 0$ let $\mathfrak{L}_n$ be the subspace of $\mathfrak{L}$ whose underlying set is the upset ${\uparrow}w_n$. 
Let also ${\bf 1}$ be the one-point Esakia space and ${\bf 2}$ the Esakia space consisting of two incomparable elements.

For two Esakia spaces $X$ and $Y$, we denote by $X \oplus Y$ the Esakia space obtained by pasting $Y$ \emph{below} $X$. 
If $\A_1, \A_2$ are Heyting algebras with Esakia duals $X_1,X_2$, then $X_1\oplus X_2$ is the dual of the sum $\A_2 + \A_1$ (see e.g., \cite[Thm. 4.1.16]{Bez-PhD}). 
We will use this construction to produce models of~$\mathsf{KG}$. 

\begin{Definition}
For $m \geq 0$ and $n \geq 1$, define
\[
\mathfrak{C}_m = \underbrace{{\bf 1}\oplus\dots\oplus  {\bf 1}}_{m-times} \, \, \text{ and } \, \,\mathfrak{G}_n = {\bf 1}\oplus \mathfrak{L}\oplus \mathfrak{L}_4\oplus \mathfrak{C}_n.
\]
\end{Definition}

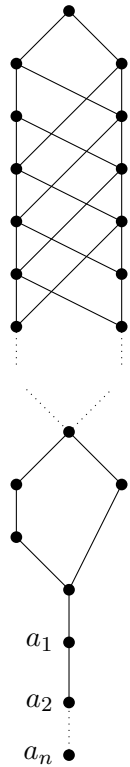
\begin{figure}
\[
\begin{tabular}{ccccccc}
\begin{tikzpicture}
    \tikzstyle{point} = [shape=circle, thick, draw=black, fill=black , scale=0.35]
   
   \node (top) at (11.9,2.8) [point] {};

\node (w0) at (11.2,2.1) [point] {};
\node (w1) at (12.6,2.1) [point] {};
\node (w2) at (11.2,1.4) [point] {};
\node (w3) at (12.6,1.4) [point] {};
\node (w4) at (11.2,0.7) [point] {};
\node (w5) at (12.6,0.7) [point] {};
\node (w6) at (11.2,0) [point] {};
\node (w7) at (12.6,0) [point] {};
\node (w8) at (11.2,-0.7) [point] {};
\node (w9) at (12.6,-0.7) [point] {};
\node (w10) at (11.2,-1.4) [point] {};
\node (w11) at (12.6,-1.4) [point] {};

\node (w12) at (11.2,-2.1) [] {};
\node (w13) at (12.6,-2.1) [] {};

\node (omega) at (11.9,-2.8) [point] {};

 \draw[dotted] (w12) -- (omega) -- (w13);

 \draw[dotted] (w10) --(w12)  (w11)--(w13);
  \draw (top) -- (w0) --  (w2) --  (w4) --  (w6) --  (w8) -- (w10) ;

    \draw (top) -- (w1) --  (w3) --  (w5) -- (w7) -- (w9) -- (w11)  ;

  \draw (w0) -- (w3) (w2) -- (w5) (w4) -- (w7) (w6) -- (w9) (w8) -- (w11) ;
  
    \draw (w1) -- (w4) (w3) -- (w6) (w5) -- (w8) (w7) -- (w10)  ;

\node (p0) at (11.2,-3.5) [point] {};
\node (p1) at (12.6, -3.5) [point] {};
\node (p2) at (11.2,-4.2) [point] {};
\node (p4) at (11.9,-4.9) [point] {};

\node[label=left:{$a_1$}] (q1) at (11.9,-5.6) [point] {};
\node[label=left:{$a_2$}] (q2) at (11.9,-6.4) [point] {};
\node[label=left:{$a_n$}] (qn) at (11.9,-7.1) [point] {};

  \draw  (omega)--(p0) --  (p2) --  (p4) -- (p1) --(omega) ;

  \draw[dotted] (qn) -- (q2);

    \draw (p4) -- (q1) -- (q2);
 
\end{tikzpicture}
\end{tabular}
\]
\caption{The poset underlying $\mathfrak{G}_n$.}
\label{fig:Gn-Nick-trick}
\end{figure}

The poset underlying $\mathfrak{G}_n$ is depicted in Figure \ref{fig:Gn-Nick-trick}. Notice that $\mathfrak{G}_n$ is the dual of the sum of Heyting algebras 
\[
\underbrace{\B_2 + \dots +  \B_2}_{n-times}+ \mathfrak{L}_4^\ast + \boldsymbol{RN}+\B_2,
\]
where $\B_2$ is the two-element Boolean algebra. Since each of the algebras $\B_2, \boldsymbol{RN}$, and $\mathfrak{L}_4^\ast$ \color{black} is one-generated, the Heyting algebra in the above display is a model of $\mathsf{KG}$, from which we deduce: 

\begin{Lemma}\label{Lem:model-of-KG}
For each  $n \geq 1$, we have 
 $\mathfrak{G}_n \vDash \mathsf{KG}$.
\end{Lemma}

We will rely on the following concept:

\begin{Definition}
For each $n\geq 1$, let $\mathcal{R}_n$ be the class of all finite rooted posets that (when endowed with the discrete topology) are continuous p-morphic images of closed upsets of $\mathfrak{G}_n$. 
\end{Definition}

From Lemma~\ref{lem: auxiliary}(\ref{item: auxiliary}) we deduce:

\begin{Lemma}\label{lem: finite-posets-m-n}
$\mathcal{R}_n = \mathsf{RFin}(\textup{Log}(\mathfrak{G}_n))$.
\end{Lemma}

We will make extensive use of the following class of Esakia spaces:

\begin{Definition}
Let $X$ be a finite Esakia space such that $X \vDash \mathsf{KG}$. 
\begin{enumerate}
\item $X$ is said to be \emph{simple} if it is a (possibly empty) finite sum of ${\bf 1}$ and~${\bf 2}$.
\item $X$ is said to be 
\emph{complex} if  $X$ is not isomorphic to $S\oplus \mathfrak{L}_k$ for any simple Esakia space $S$ and $k \geq 0$.
\end{enumerate}
\end{Definition}

The next result is a straightforward adaptation of \cite[Thm.~4.5.1]{Bez-PhD}:

\begin{Lemma}\label{lem : long structure lemma}
Let $X, Y$ be finite rooted Esakia spaces with $Y$ complex. The following are equivalent.
\begin{enumerate}
\item $X$ is a continuous p-morphic image of a closed upset of ${\bf 1}\oplus \mathfrak{L}\oplus Y$.
\item $X$ is isomorphic to ${\bf 1} \oplus S\oplus  \mathfrak{L}_k$ or  is a continuous p-morphic image of a rooted  upset of ${\bf 1}\oplus S\oplus {\bf 1}\oplus  Y$ for some simple Esakia space $S$ and $k \geq 0$.
\end{enumerate} 
\end{Lemma}

As a consequence, we obtain:

\begin{Theorem}\label{thm: struct}
For $n \geq 1$, if $X$ is a finite rooted continuous p-morphic image of a closed upset of $\mathfrak{G}_n$, then there are a simple Esakia space $S$, $k \geq 0$, and $m \leq n$ such that $X$ is isomorphic to
\[
{\bf 1}\oplus S\oplus  \mathfrak{L}_k\, \, \text{  or  }\, \, {\bf 1}\oplus S\oplus {\bf 1}\oplus \mathfrak{L}_4\oplus \mathfrak{C}_m\, \, 
\text{  or  }\, \, {\bf 1} \oplus \mathfrak{L}_4\oplus \mathfrak{C}_m\, \, \text{  or  } \, \,  {\bf 1}. 
\]
\end{Theorem}
\begin{proof}
Let $X$ be a continuos p-morphic image of a closed upset of $\mathfrak{G}_n = {\bf 1} \oplus \mathfrak{L} \oplus  \mathfrak{L}_4 \oplus \mathfrak{C}_n$. Letting $Y = \mathfrak{L}_4 \oplus \mathfrak{C}_n$, we have  $\mathfrak{G}_n = {\bf 1} \oplus \mathfrak{L} \oplus Y$. Since $Y$ is complex, Lemma \ref{lem : long structure lemma} implies that $X$  
is either  
isomorphic to ${\bf 1}\oplus S\oplus  \mathfrak{L}_k$ or  is a continuous p-morphic image of a rooted upset $U$ of ${\bf 1}\oplus S\oplus {\bf 1}\oplus Y$ for some simple Esakia space $S$ and $k \geq 0$. 

In the former case, we are done. In the latter case, as $U$ is a rooted upset  of ${\bf 1}\oplus S\oplus {\bf 1}\oplus Y$, it is of the form
\[
{\bf 1}\oplus S \oplus {\bf 1} \oplus \mathfrak{L}_4\oplus \mathfrak{C}_t \, \, \text{ or } \, \, {\bf 1}\oplus S' \oplus {\bf 1} \, \, \text{ or } \, \, {\bf 1}
\]
for some $t \leq n$ and simple Esakia space $S'$. The continuous p-morphic images of the Esakia spaces in the above display are of the form
\[
 {\bf 1}\oplus S'' \oplus {\bf 1} \oplus \mathfrak{L}_4\oplus \mathfrak{C}_m\, \, 
\text{  or  }\, \, {\bf 1} \oplus \mathfrak{L}_4\oplus \mathfrak{C}_m \, \, \text{ or }\, \, {\bf 1}\oplus S''\oplus {\bf 1} \, \,  \, \, \text{ or }\, \, \bf{1}
\]
for some simple Esakia space $S''$ and $m \leq t \leq n$. Since ${\bf 1}\oplus S''\oplus {\bf 1}$ is of the form ${\bf 1}\oplus S\oplus  \mathfrak{L}_k$ for $k = 0, 1$, the result follows.
\end{proof}

As a consequence, we obtain the following  characterization of the posets in $\mathcal{R}_n$. 

\begin{Corollary}\label{cor: struct}
A finite rooted poset $X$ belongs to $\mathcal{R}_n$ iff $X$ is isomorphic to 
\[
{\bf 1}\oplus S\oplus  \mathfrak{L}_k \, \, \text{  or  }\, \, {\bf 1}\oplus S\oplus {\bf 1} \oplus \mathfrak{L}_4\oplus \mathfrak{C}_m \, \, \text{  or  }\, \,  {\bf 1} \oplus \mathfrak{L}_4\oplus \mathfrak{C}_m\, \,  \text{  or  }\, \,  {\bf 1}
\]
\color{black} for some simple Esakia space $S$, $k \geq 0$, and $m\leq n$.
\end{Corollary}
\begin{proof}
To prove the implication from left to right, suppose that $X \in \mathcal{R}_n$. By the Dual Jankov Lemma, $X \not\vDash \mathcal{J}(X)$. Therefore,  $\mathfrak{G}_n\not\vDash \mathcal{J}(X)$ by Lemma \ref{lem: finite-posets-m-n}.. 
Thus, 
$X$ is a continuous p-morphic image of a closed upset of $\mathfrak{G}_n$  by the Dual Jankov Lemma. Consequently, we can apply Theorem~\ref{thm: struct} to  obtain that $X$ is isomorphic to one of the posets in the above display. 

Next we turn to proving the implication from right to left. Suppose that $X$ is isomorphic to one of the posets in the above display. Let  $Y = \mathfrak{L}_4\oplus \mathfrak{C}_n$. Then $Y$ is a finite rooted complex poset such that $\mathfrak{G}_n = {\bf 1}\oplus \mathfrak{L}\oplus Y$. From Lemma \ref{lem : long structure lemma} it follows that $X$ is a continuous p-morphic image of a closed upset of $\mathfrak{G}_n$. By the definition of $\mathcal{R}_n$ this amounts to $X \in \mathcal{R}_n$.
\end{proof}

\begin{Definition}
For $n\geq 1$ define: 
\begin{align*}
\mathsf{L}_0^n & = \textup{Log} (\mathcal{R}_n);\\
\mathsf{L}_1^n & = \textup{Log} (\mathcal{R}_n\cup \{\mathfrak{G}_1\});\\
\mathsf{L}_2^n & = \textup{Log} (\mathcal{R}_n\cup \{\mathfrak{G}_2\});\\
&\vdots\\
 \mathsf{L}_n^n & = \textup{Log} (\mathcal{R}_n\cup \{\mathfrak{G}_n\}) =  \textup{Log}(\mathfrak{G}_n),
\end{align*}
where the equality $\textup{Log} (\mathcal{R}_n\cup \{\mathfrak{G}_n\})=\textup{Log}(\mathfrak{G}_n)$ holds by Lemma \ref{lem: finite-posets-m-n}. 
\end{Definition}

When the integer $n$ is clear from the context, we will drop the superscript and write $\mathsf{L}_0, \dots, \mathsf{L}_n$ instead of $\mathsf{L}_0^n, \dots, \mathsf{L}_n^n$. In view of Lemma \ref{Lem:model-of-KG} 
we have:

\begin{Lemma}
$\mathsf{L}_0, \dots, \mathsf{L}_n$ are extensions of $\mathsf{KG}$.
\end{Lemma}

Our aim is to prove that the fmp span of $\textup{Log}(\mathfrak{G}_n)$ is precisely the set $\{ \mathsf{L}_0, \dots, \mathsf{L}_n \}$, and hence that the degree of fmp of $\textup{Log}(\mathfrak{G}_n)$ is $n+1$. Since $n\geq 1$ was arbitrary  and we already proved that $\mathsf{KG}$ has the degree of fmp 1, this will show that there are extensions of $\mathsf{KG}$ with an arbitrary \textit{finite} degree of fmp. 

We begin by the following simple observation. 

\begin{Lemma}\label{lem: p-m}
The following conditions hold.
\begin{enumerate}
\item \label{item:pm:1} If $m \leq n$, then $\mathfrak{G}_m$ is a continuous p-morphic image of $\mathfrak{G}_n$.
\item \label{item:pm:2} For each $n \geq 1$ the Esakia space ${\bf 1}\oplus \mathfrak{L}_4\oplus \mathfrak{C}_n$ is a continuous p-morphic image of $\mathfrak{G}_n$.
\end{enumerate}
\end{Lemma}
\begin{proof}
(\ref{item:pm:1}) Define $\alpha \colon \mathfrak{G}_n \to \mathfrak{G}_m$ by sending the least $n-m$ points of $\mathfrak{G}_n$ to the root of $\mathfrak{G}_m$ and any other point in $\mathfrak{G}_n$ to its copy in  $\mathfrak{G}_m$. It is straightforward to check that $\alpha$ is an onto continuous p-morphism. 

(\ref{item:pm:2}) Let $X = {\bf 1}\oplus \mathfrak{L}_4\oplus \mathfrak{C}_n$ and let $\alpha \colon \mathfrak{G}_n\to X$ be the map 
that sends ${\bf 1}\oplus \mathfrak{L}$ to the top element of $X$ and is the identity on the rest of 
$\mathfrak{G}_n$. It is straightforward to check that $\alpha$ is an onto continuous p-morphism. 
\end{proof}

We will show that the  $\mathsf{L}_1, \dots, \mathsf{L}_n$   form a descending chain of logics with the same finite models.

\begin{Lemma}\label{Lem:VI-Rn-description} The following conditions hold.
\begin{enumerate}
\item \label{item:VI-Rn:1} $\mathsf{L}_n\subsetneq \dots \subsetneq \mathsf{L}_0$. 
\item \label{item:VI-Rn:2} For all $i \leq n$ we have $\mathcal{R}_n = \mathsf{RFin}(\mathsf{L}_i)$.
\end{enumerate} 
\end{Lemma}

\begin{proof}
(\ref{item:VI-Rn:1}) By Lemma~\ref{lem: p-m}(\ref{item:pm:1}) we have the inclusions $\mathsf{L}_n\subseteq \dots \subseteq \mathsf{L}_0$. To show that these inclusions are proper, consider the subframe formula $\beta(\mathfrak{L}_6)$.\footnote{We note that the same proof would work if $\mathfrak{L}_6$ is replaced by any $\mathfrak{L}_k$ with $k \geq 6$.}
  For each $i\leq n-1$ let 
\[
\varphi_i \coloneqq \beta(\mathfrak{L}_6)\vee \mathcal{J}({\bf 1}\oplus \mathfrak{L}_4\oplus \mathfrak{C}_{i+1}).\]
It is enough to show that $\varphi_i\in \mathsf{L}_i \smallsetminus \mathsf{L}_{i+1}$ for each $i \leq n-1$. To this end, let  $i \leq n-1$.

To see  that $\varphi_i\notin  \mathsf{L}_{i+1}$, it is enough to show that $\mathfrak{G}_{i+1}\nvDash \varphi_i$. Since $\mathfrak{L}_6$ is a clopen subset of $\mathfrak{G}_{i+1}$, Theorem~\ref{Thm:subframe-formulas}(\ref{subframe-formulas 1}) implies that 
$\mathfrak{G}_{i+1}\nvDash \beta(\mathfrak{L}_6)$. Moreover,  
 ${\bf 1}\oplus \mathfrak{L}_4\oplus \mathfrak{C}_{i+1}$ is a continuous p-morphic image of $\mathfrak{G}_{i+1}$ by Lemma~\ref{lem: p-m}(\ref{item:pm:2}).
  Thus, $\mathfrak{G}_{i+1}\not \vDash \mathcal{J}({\bf 1}\oplus \mathfrak{L}_4\oplus \mathfrak{C}_{i+1})$ by the Dual Jankov Lemma. Since a disjunction holds in a rooted Esakia space iff one of the disjuncts does, we conclude that $\mathfrak{G}_{i+1}\not \vDash \varphi_i$, and hence $\varphi_i\notin  \mathsf{L}_{i+1}$.
 
To prove that $\varphi_i\in  \mathsf{L}_i$, it is sufficient to show that $\mathcal{R}_n\cup \{\mathfrak{G}_i\}\vDash \varphi_i$.\ From Theorem~\ref{thm: struct} it follows that ${\bf 1}\oplus \mathfrak{L}_4\oplus \mathfrak{C}_{i+1}$ is not a contionuous p-morphic image of a closed upset of $\mathfrak{G}_i$. Therefore, 
$\mathfrak{G}_i\vDash \mathcal{J}({\bf 1}\oplus \mathfrak{L}_4\oplus \mathfrak{C}_{i+1})$ by the Dual Jankov Lemma, and hence $\mathfrak{G}_i\vDash \varphi_i$.
Next, let $X\in \mathcal{R}_n$. If $X \vDash \mathcal{J}({\bf 1}\oplus \mathfrak{L}_4\oplus \mathfrak{C}_{i+1})$, then $X \vDash \varphi_i$ as desired. Therefore, we may assume 
that $X \nvDash \mathcal{J}({\bf 1}\oplus \mathfrak{L}_4\oplus \mathfrak{C}_{i+1})$. By the Fine Lemma, ${\bf 1}\oplus \mathfrak{L}_4\oplus \mathfrak{C}_{i+1}$ is a p-morphic image of an upset of $X$. Together with Corollary~\ref{cor: struct}, this implies that
\[
X \cong {\bf 1}\oplus S\oplus {\bf 1} \oplus \mathfrak{L}_4\oplus \mathfrak{C}_{m}  \mbox{ or } 
X \cong  {\bf 1} \oplus \mathfrak{L}_4\oplus \mathfrak{C}_{m}
\]
 for some simple Esakia space $S$ and $i+1\leq m\leq n$ (because the other configurations in  Corollary~\ref{cor: struct} cannot have   
${\bf 1}\oplus \mathfrak{L}_4\oplus \mathfrak{C}_{i+1}$ as a p-morphic image of one of their upsets \cite[Thm.~4.4.12(1)]{Bez-PhD}). 
The above display guarantees that 
$\mathfrak{L}_6$ is not a p-morphic image of  a subposet of $X$. 
Therefore, $X\vDash  \beta(\mathfrak{L}_6)$ by Theorem~\ref{Thm:subframe-formulas}(\ref{subframe-formulas 2}), and hence
$X\vDash  \varphi_i$.  
Thus, $\mathcal{R}_n\cup \{\mathfrak{G}_i\}\vDash \varphi_i$, yielding that 
the inclusions are proper. 

 (\ref{item:VI-Rn:2}) By the definition of $ \mathsf{L}_i$ we have that $\mathcal{R}_n\subseteq \mathsf{RFin}(\mathsf{L}_i)$.  To prove the other inclusion, let  $X \in \mathsf{RFin}(\mathsf{L}_i)$.  
Since $X\nvDash \mathcal{J}(X)$ by the Dual Jankov Lemma, we have  $\mathcal{J}(X)\notin \mathsf{L}_i$. Because $\mathsf{L}_i = \textup{Log} (\mathcal{R}_n\cup \{\mathfrak{G}_i\})$, either there exists some  $Y\in \mathcal{R}_n$ such that $Y\nvDash \mathcal{J}(X)$ or  
$\mathfrak{G}_i\nvDash \mathcal{J}(X)$. By the Dual Jankov Lemma, $X$ is a continuous p-morphic image of a closed upset of either some $Y\in \mathcal{R}_n$ or $\mathfrak{G}_i$. In the former case, it is clear that $X \in \mathcal{R}_n$. In the latter case, apply Lemma~\ref{lem: p-m}(\ref{item:pm:1}) to obtain that $X$ is a continuous p-morphic image of a closed upset of a p-morphic image of $\mathfrak{G}_n$. This easily implies that $X$ is also a continuous p-morphic image of a closed upset of $\mathfrak{G}_n$, and hence is a member of $\mathcal{R}_n$.
\end{proof}

Together with Lemma \ref{lem: finite-posets-m-n} this yields that $\mathsf{RFin}(\textup{Log}(\mathfrak{G}_n)) = \mathsf{RFin}(\mathsf{L}_i)$ for every $i \leq n$. As a consequence, we obtain the following.

\begin{Lemma}\label{Cor:n+1-logics-in-span}
$\mathsf{L}_0, \dots, \mathsf{L}_n$ are $n+1$ distinct elements of the fmp span of $\textup{Log}(\mathfrak{G}_n)$.
\end{Lemma}

To show that there are no other logics in the fmp span of $\textup{Log}(\mathfrak{G}_n)$, we rely on the following observations. 

\begin{Lemma}\label{lem: extr}
Let $X, Y, Z$ be Esakia spaces. The following conditions hold.
\begin{enumerate}
\item \label{item:extr:1} $X \oplus \mathfrak{L}\oplus Z$ is a continuous p-morphic image of $X \oplus \mathfrak{L}\oplus Y\oplus Z$.
\item \label{item:extr:2} If $X, Z\vDash \mathsf{KG}$, then $\textup{Log}(X\oplus \mathfrak{L}\oplus \mathfrak{L}\oplus Z) = \textup{Log}(X\oplus \mathfrak{L}\oplus Z)$. 
\end{enumerate}
\end{Lemma}

\begin{proof}
(\ref{item:extr:1}) The map $\alpha \colon X \oplus \mathfrak{L}\oplus Y\oplus Z \to X \oplus \mathfrak{L}\oplus Z$ that sends the points of $Y$ to the bottom element of $\mathfrak{L}$ and is the identity on the rest of the points is an onto continuous p-morphism.

(\ref{item:extr:2}) By (\ref{item:extr:1}) we have that $X\oplus \mathfrak{L}\oplus Z $  is a continuous p-morphic image of 
$X\oplus \mathfrak{L}\oplus \mathfrak{L}\oplus Z$. Thus, $\textup{Log}(X\oplus \mathfrak{L}\oplus \mathfrak{L}\oplus Z) \subseteq \textup{Log}(X\oplus \mathfrak{L}\oplus Z) $. The reverse inclusion 
follows from \cite[Lem~4.4.9(4)]{Bez-PhD}. 
\end{proof}

Consider the si-logic
\[
\mathsf{RN.KC} = \mathsf{RN} + (\neg p \vee\neg\neg p).
\]
It is well known that a rooted Esakia space $X$ validates $\mathsf{RN.KC}$ iff it validates $\mathsf{RN}$ 
and it has 
a maximum. 

Now, for each si-logic $\mathsf{L}$, let $\mathsf{FGR}(\mathsf{L})$ be the class of Esakia spaces $X$ such that $X^\ast$ is finitely generated, SI, and $X\vDash \mathsf{L}$. Clearly $\mathsf{FGR}(\mathsf{IPC})$ is the class of all Esakia spaces $X$ with an isolated root  such that $X^\ast$ is finitely generated. 

\begin{Theorem}\label{thm: KG-struct}
The following conditions hold for each $X\in \mathsf{FGR}(\mathsf{IPC})$. 
\begin{enumerate}
\item\label{item:KG-struct:1} $X\vDash \mathsf{KG}$ iff $X$ is isomorphic to $X_1\oplus \dots \oplus X_n\oplus {\bf 1}$, where each $X_i$ is isomorphic to $\mathfrak{L}$ or to a finite upset of $\mathfrak{L}$.
\item\label{item:KG-struct:2} $X\vDash \mathsf{RN}$ iff $X$ is isomorphic to $X_1\oplus \dots \oplus X_n\oplus \mathfrak{L}_k$, where $k\geq 0$ and each $X_i$ is isomorphic to ${\bf 1}$, ${\bf 2}$, or $\mathfrak{L}$.
\item\label{item:KG-struct:3} $X\vDash \mathsf{RN.KC}$ iff $X$ is isomorphic to ${\bf 1} \oplus X_1\oplus \dots \oplus X_n\oplus  \mathfrak{L}_k$, where $k\geq 0$ and each $X_i$ is isomorphic to ${\bf 1}$, ${\bf 2}$, or $\mathfrak{L}$
\item\label{item:KG-struct:4} If $X$ is infinite and  $X \vDash \mathsf{RN.KC}$, then $\textup{Log}(X) = \textup{Log}(\mathfrak{{\bf 1}\oplus L})$ and $\textup{Log}(X)$ has the fmp.
\end{enumerate}
\end{Theorem}

\begin{proof}
(\ref{item:KG-struct:1}) By \cite[Thm.~4.3.9]{Bez-PhD}, $X\vDash \mathsf{KG}$ iff $X$ is isomorphic to $X_1\oplus \dots \oplus X_n\oplus \mathfrak{L}_k$, where $k \geq 0$ and each $X_i$ 
is isomorphic to $\mathfrak{L}$ or to a finite upset of $\mathfrak{L}$. Next observe that if $k \in \{ 0,1 \}$, then $\mathfrak{L}_k$ is isomorphic to ${\bf 1}$; if $k=2$, then $\mathfrak{L}_k$ is isomorphic to ${\bf 1}\oplus {\bf 1}$; and if $k>2$, then $\mathfrak{L}_k$ is isomorphic to $Y \oplus {\bf 1}$ where $Y$ is the upset of $\mathfrak{G}$ generated by $\{w_{k-2}, w_{k-3}\}$. Thus, $\mathfrak{L}_k$ is always isomorphic to $Z\oplus {\bf 1}$ for a (possibly empty) finite upset $Z$ of $\mathfrak{L}$.

 (\ref{item:KG-struct:2}) This is \cite[Thm.~4.4.12(1)]{Bez-PhD}.
 
 (\ref{item:KG-struct:3}) This follows from (\ref{item:KG-struct:2})  by observing that a rooted Esakia space validates $\neg p\vee \neg\neg p$ iff it has the maximum.

(\ref{item:KG-struct:4}) We recall that $ \textup{Log}(\mathfrak{{\bf 1}\oplus L}) = \mathsf{RN.KC}$ (see \cite[Thm.~4.6.4]{Bez-PhD}). Therefore, from $X \vDash \mathsf{RN.KC}$ it follows that $\textup{Log}(\mathfrak{{\bf 1}\oplus L})   \subseteq \textup{Log}(X)$. Because  $X$ is infinite, it follows from  (\ref{item:KG-struct:3})  that $X$ is isomorphic to ${\bf 1} \oplus Y\oplus \mathfrak{L} \oplus Z\oplus \mathfrak{L}_k$, where $k\geq 0$ and $Y$ and $Z$ are possibly empty Esakia spaces. By Lemma~\ref{lem: extr}(\ref{item:extr:1}), ${\bf 1} \oplus Y\oplus \mathfrak{L}$ is a continuous p-morphic image of ${\bf 1} \oplus Y\oplus \mathfrak{L} \oplus Z\oplus \mathfrak{L}_k$. Next identify the points in $Y$ with the maximum to obtain that ${\bf 1} \oplus \mathfrak{L}$ is a continuous p-morphic image of ${\bf 1} \oplus Y\oplus \mathfrak{L}$. Thus, ${\bf 1} \oplus \mathfrak{L}$ is also a continuous p-morphic image of $X$, and hence $\textup{Log}(X) \subseteq \textup{Log}(\mathfrak{{\bf 1}\oplus L})$. This  shows that $\textup{Log}(X) = \textup{Log}(\mathfrak{{\bf 1}\oplus L})$. Finally, since every extension of $\mathsf{RN}$ has the fmp \cite[Thm.~4.4.13]{Bez-PhD}, we conclude that $\textup{Log}(X)$ has the fmp. 
\end{proof}

We will make use of the following observation.

\begin{Lemma}\label{lem: nnnnn}
\cite[Cor.~4.2.7]{Bez-PhD} If $S$ is a simple Esakia space, then $S\oplus \mathfrak{L}$ and $S\oplus {\bf 1}$ are continuous p-morphic images of $\mathfrak{L}$.
\end{Lemma}

We are ready for the key result of this section.

\begin{Theorem}\label{Thm:main-trick-n+1}
Let $\mathsf{L}$ be an extension of $\mathsf{KG}$. If $\mathsf{L}\in\textup{fmp}(\textup{Log}(\mathfrak{G}_n))$, then $\mathsf{L} = \mathsf{L}_i$ for some $i\leq n$. 
\end{Theorem}

\noindent \textit{Proof.}
Let $\mathsf{L}$ be an extension of $\mathsf{KG}$ such that $\mathsf{L}\in\textup{fmp}(\textup{Log}(\mathfrak{G}_n))$. Then $\mathsf{RFin}(\mathsf{L}) = \mathcal{R}_n$ by Lemma~\ref{lem: finite-posets-m-n}.
Since $\mathsf{L}' = \textup{Log}(\mathsf{FGR}(\mathsf{L}'))$ for every si-logic $\mathsf{L}'$, this implies that
  $$\mathsf{L} = \textup{Log}(\mathsf{FGR}(\mathsf{L})) = \bigcap  
\{\textup{Log}(\mathcal{R}_n \cup \{X\}): X\in \mathsf{FGR}(\mathsf{L})\}.$$
 Thus, in order to prove that $\mathsf{L} = \mathsf{L}_i$ for some $i \leq n$, it is sufficient to show 
that for each  $X\in \mathsf{FGR}(\mathsf{L})$, there is $j \leq n$ with $\textup{Log}(\mathcal{R}_n\cup \{ X\})= \mathsf{L}_j$. For 
in this case, we can take $i$ to be the maximum of the $j$ by Lemma~\ref{Lem:VI-Rn-description}(\ref{item:VI-Rn:1}). 

Let $X\in \mathsf{FGR}(\mathsf{L})$. 
First suppose that $X$ is finite. Then $X \in \mathcal{R}_n$, which implies that 
\[
\textup{Log}(\mathcal{R}_n\cup\{ X\}) =  \textup{Log}(\mathcal{R}_n) =\mathsf{L}_0.
\]
Next suppose that $X$ is infinite. 
By Theorem~\ref{thm: KG-struct}(\ref{item:KG-struct:1}) we may assume that $X = X_1\oplus \dots \oplus X_m\oplus {\bf 1}$,  where $m\geq 1$ and  each $X_i$  is either $\mathfrak{L}$ or a finite upset of $\mathfrak{L}$. \color{black}

\begin{Claim}
$X_1 = {\bf 1}$.
\end{Claim}
 \begin{proof}[Proof of the Claim]
We first show that $X$ has a maximum.  If not, then a finite rooted upset of $\mathfrak{L}$ containing two maximal points is a rooted upset of $X$, so it  is in $\mathcal{R}_n$. But 
none of these belongs to $\mathcal{R}_n$ by  Corollary~\ref{cor: struct}, a contradiction. Thus, we may assume that  $X_1$ is either ${\bf 1}$  or the two-element chain $\mathfrak{L}_2$. But $\mathfrak{L}_2$ is isomorphic to ${\bf 1}\oplus {\bf 1}$.
So without the loss of generality, we may assume that $X_1 = {\bf 1}$ (otherwise we renumber the summands: the new $X_2$ becomes ${\bf 1}$, the new $X_3$ becomes the old $X_2$, etc.).
\end{proof}
  
  Since $X$ is infinite, one of the $X_k$ must be $\mathfrak{L}$. Let $k$ be the least such. Then
  \[
  X  = {\bf 1}\oplus S_1\oplus \mathfrak{L}\oplus X_{{k}+1}\oplus \dots \oplus X_m\oplus {\bf 1},
  \]
   where 
$S_1 = X_2\oplus \dots\oplus X_{{k}-1}$.

\begin{Claim}
$S_1$ is a simple Esakia space. 
\end{Claim}

 \begin{proof}[Proof of the Claim]
 Suppose the contrary, with a view to contradiction. Then some poset among $X_2, \dots, X_{{k}-1}$ is not simple. Recall that these posets are all different from $\mathfrak{L}$ by assumption and, therefore, each of them is a  finite upset of $\mathfrak{L}$. Consequently, some poset among $X_2, \dots, X_{{k}-1}$ is a finite non-simple upset of $\mathfrak{L}$. Now, define $Y = {\bf 1}\oplus S_1 \oplus \mathfrak{L}_4$. Since $S_1 = X_2\oplus \dots\oplus X_{{k}-1}$ and both $\mathfrak{L}_4$ and some poset among $X_2, \dots, X_{{k}-1}$ are a non-simple finite upsets of $\mathfrak{L}$, we obtain that $Y$ can be written as a finite sum of posets two of which are finite non-simple upsets of $\mathfrak{L}$. Together with Corollary~\ref{cor: struct}, this yields that $Y \notin \mathcal{R}_n$. On the other hand, from the definition of $Y$ and the equality
\[
X  = {\bf 1}\oplus S_1\oplus \mathfrak{L}\oplus X_{{k}+1}\oplus \dots \oplus X_m\oplus {\bf 1}
\]
it follows that $Y$ is a finite principal upset of $X$, hence belongs to $\mathsf{RFin}(\mathsf{L})$.  Since $\mathsf{RFin}(\mathsf{L}) = \mathcal{R}_n$ by assumption, we conclude that $Y \in \mathcal{R}_n$, a contradiction.
\end{proof}

Iterating the argument described above, we obtain the following:

\begin{Claim}\label{Claim:X:Nick}
One of the following conditions holds.
\begin{enumerate}
\item \label{item:Claim:X:Nick:1} There exist simple Esakia spaces $S_1, \dots, S_p$ such that
\[
X  \cong {\bf 1}\oplus S_1\oplus \mathfrak{L}\oplus S_2\oplus \dots \oplus \mathfrak{L} \oplus S_p \oplus \mathbf{1};
\]
\item \label{item:Claim:X:Nick:2}There exist simple Esakia spaces $S_1, \dots, S_p$ and finite upsets $X_{i}, \dots, X_m$ of $\mathfrak{L}$ with $X_i$ non-simple such that
\[
X  \cong {\bf 1}\oplus S_1\oplus \mathfrak{L}\oplus S_2\oplus \dots \oplus \mathfrak{L} \oplus S_p \oplus X_{i}\oplus \dots \oplus X_m\oplus {\bf 1}.
\]
\end{enumerate}
\end{Claim}

Suppose first that Condition (\ref{item:Claim:X:Nick:1}) of Claim \ref{Claim:X:Nick} holds. Then $X \vDash \mathsf{RN.KC}$ by Theorem \ref{thm: KG-struct}(\ref{item:KG-struct:3}). By  Condition (\ref{item:KG-struct:4}) of the same theorem we obtain that $\textup{Log}(X)$ has the fmp. Now, from the assumption that $X \in \mathsf{FGR}(\mathsf{L})$ and $\mathsf{RFin}(\mathsf{L}) = \mathcal{R}_n$ it follows that $\mathsf{RFin}(\textup{Log}(X)) \subseteq \mathsf{RFin}(\mathsf{L}) = \mathcal{R}_n$. As $\textup{Log}(X)$ has the fmp, this implies that $\textup{Log}(\mathcal{R}_n \cup \{ X \}) = \textup{Log}(\mathcal{R}_n) = \mathsf{L}_0$.

It only remains to consider the case where Condition (\ref{item:Claim:X:Nick:2}) of Claim \ref{Claim:X:Nick} holds.  If there exists $t \geq 0$ such that 
\[
X  \cong {\bf 1}\oplus S_1\oplus \mathfrak{L}\oplus S_2\oplus \dots \oplus \mathfrak{L} \oplus S_p \oplus \mathfrak{L}_t,
\]
then we can repeat the argument detailed in the previous paragraph. Therefore, we may assume that there is no $t \geq 0$ for which the above display holds. Under this assumption, we will prove the following:

\begin{Claim}\label{Claim:X:b:Nick}
There exists $k \leq n$ such that
\[
X  \cong {\bf 1}\oplus S_1\oplus \mathfrak{L}\oplus S_2\oplus \dots \oplus \mathfrak{L} \oplus S_p \oplus \mathfrak{L}_4\oplus \mathfrak{C}_k.
\]
\end{Claim}

\begin{proof}[Proof of the Claim]
Recall from Condition (\ref{item:Claim:X:Nick:2}) of Claim \ref{Claim:X:Nick} that
\[
X  \cong {\bf 1}\oplus S_1\oplus \mathfrak{L}\oplus S_2\oplus \dots \oplus \mathfrak{L} \oplus S_p \oplus X_{i}\oplus \dots \oplus X_m\oplus {\bf 1},
\]
where each $S_j$ is simple, each $X_j$ is a finite upset of $\mathfrak{L}$, and $X_i$ is not simple. We may also assume that each $X_j$ is nonempty.

We have two cases depending on whether $i = m$ or $i < m$. First suppose that $i = m$. Then
\[
X  \cong {\bf 1}\oplus S_1\oplus \mathfrak{L}\oplus S_2\oplus \dots \oplus \mathfrak{L} \oplus S_p \oplus X_{i}\oplus {\bf 1}.
\]
Since we assumed that $X$ is not of the form ${\bf 1}\oplus S_1\oplus \mathfrak{L}\oplus S_2\oplus \dots \oplus \mathfrak{L} \oplus S_p \oplus \mathfrak{L}_t$ for any $t \geq 0$, the non-simple finite upset $X_i$ of $\mathfrak{L}$ must be rooted (otherwise $X_i \oplus \mathbf{1} = \mathfrak{L}_t$ for some $t \geq 0$). Therefore, $X_i = \mathfrak{L}_s$ for some $s \geq 0$. Because $X_i$ is not simple, we must have $s \geq 4$. If $s = 4$, then 
\[
X = {\bf 1}\oplus S_1\oplus \mathfrak{L}\oplus S_2\oplus \dots \oplus \mathfrak{L} \oplus S_p \oplus \mathfrak{L}_{4}\oplus \mathfrak{C}_1
\]
and
 we are done. 
We will show that the case where $s > 4$ cannot happen. Suppose the contrary. Then 
\[
X  \cong {\bf 1}\oplus S_1\oplus \mathfrak{L}\oplus S_2\oplus \dots \oplus \mathfrak{L} \oplus S_p \oplus \mathfrak{L}_{s}\oplus \mathbf{1}.
\]
Therefore, ${\bf 1}\oplus \mathfrak{L}_{s}\oplus \mathbf{1}$ is a continuous p-morphic image of $X$ obtained by collapsing the top part of $X$. Since $X \vDash \mathsf{L}$, we have ${\bf 1}\oplus \mathfrak{L}_{s}\oplus \mathbf{1} \in \mathsf{RFin}(\mathsf{L}) = \mathcal{R}_n$. Together with $s \geq 5$, this contradicts Corollary \ref{cor: struct}.

Now suppose $i < m$. We begin by proving that
\begin{equation}\label{Eq:Xi+1equalXmequals1}
X_{i+1} \oplus \dots \oplus X_m \oplus \mathbf{1}  \cong \mathfrak{C}_{q}
\end{equation}
for some $q \geq 2$. Suppose the contrary. Then there exists some $j > i$ which contains two incomparable elements. In this case, one of the following is a continuous p-morphic image of $X$, depending on whether $i+1 < j$ or $i +1 = j$:
\[
\mathbf{1} \oplus X_i \oplus \mathbf{1}\oplus X_j \oplus \mathbf{1} \, \, \text{ or } \, \, \mathbf{1} \oplus X_i \oplus X_j \oplus \mathbf{1}.
\]
Now, since $X_j$ is a finite upset of $\mathfrak{L}$ containing two incomparable points, $X_j$ must contain the two maximal elements of $\mathfrak{L}$. It is therefore easy to show that if $X_j$ consists of two disjoint points or of the disjoint union of the two element chain and a point, then $\mathbf{2}$ is a continuous  p-morphic image of $X_j$, and in all other cases $\mathbf{2} \oplus \mathbf{1}$ is a continuous p-morphic image of $X_j$.
By collapsing the summand $X_j$ in this manner in the Esakia spaces in the above display, we obtain that one of the following is a continuous p-morphic image of $X$:
\[
\mathbf{1} \oplus X_i \oplus \mathbf{1} \oplus \mathbf{2} \oplus \mathbf{1} \, \, \text{ or } \, \, \mathbf{1} \oplus X_i \oplus \mathbf{1} \oplus \mathbf{2}\oplus \mathbf{1} \oplus \mathbf{1} \, \, \text{ or } \, \, \mathbf{1} \oplus X_i \oplus \mathbf{2} \oplus \mathbf{1}\, \, \text{ or } \, \, \mathbf{1} \oplus X_i \oplus \mathbf{2} \oplus \mathbf{1} \oplus \mathbf{1}.
\]
Since the finite continuous p-morphic images of $X$ belong to $\mathsf{RFin}(\mathsf{L}) = \mathcal{R}_n$, one of the Esakia spaces above should belong to $\mathcal{R}_n$.  Observe that none of them is of the form ${\bf 1}$  or $ {\bf 1} \oplus \mathfrak{L}_4\oplus \mathfrak{C}_s$ for some $s \leq n$. We also claim that none of them is of the form 
 ${\bf 1}\oplus S\oplus {\bf 1} \oplus \mathfrak{L}_4\oplus \mathfrak{C}_s$ for some $s \leq n$ and simple $S$. 
This is clear for the first two posets displayed above. 
For the other two this is a consequence of the fact that $\mathfrak{L}_4$ cannot be isomorphic to a poset of the form $Y\oplus {\bf 2} \oplus \mathfrak{C}_s$ for any $Y$.   Since none of the posets in the above display is of the form ${\bf 1}\oplus S\oplus {\bf 1} \oplus \mathfrak{L}_4\oplus \mathfrak{C}_s$, we can apply Corollary \ref{cor: struct} to obtain that one of these posets must be of the form ${\bf 1}\oplus S\oplus  \mathfrak{L}_s$ for some $s \geq 0$ and simple $S$. By inspecting these posets and using the structure of the posets $\mathfrak{L}_s$, we see that $\mathfrak{L}_s$ should be either $\mathbf{1}$ or $\mathbf{2}\oplus \mathbf{1}$ or $\mathbf{1}\oplus \mathbf{1}$. In all three cases, $S$ should contain $X_i$ as a subposet, which is impossible because $S$ is simple and $X_i$ is not by assumption. This establishes  (\ref{Eq:Xi+1equalXmequals1}).

Consequently, we obtain that
\begin{equation}\label{Eq:Xi+1equalXmequals2}
X  \cong {\bf 1}\oplus S_1\oplus \mathfrak{L}\oplus S_2\oplus \dots \oplus \mathfrak{L} \oplus S_p \oplus X_{i}\oplus \mathfrak{C}_q
\end{equation}
for some $q \geq 2$. From the above display it follows that $\mathbf{1} \oplus X_i \oplus \mathfrak{C}_{q}$ is a continuous p-morphic image of $X$, and so belongs to $\mathcal{R}_n$. Since $\mathbf{1} \oplus X_i \oplus \mathfrak{C}_{q}$ is not isomorphic to ${\bf 1}$  (because $q \geq 2$) and belongs to $\mathcal{R}_n$, we can apply Corollary \ref{cor: struct}, obtaining that there exist a simple Esakia space $S$, $s \geq 0$, and $k \leq n$ such that one of the following conditions holds:
\benroman
\item\label{item : tom : addition : 1} $\mathbf{1} \oplus X_i \oplus \mathfrak{C}_{q} \cong {\bf 1}\oplus S\oplus  \mathfrak{L}_s$;
\item\label{item : tom : addition : 2} $\mathbf{1} \oplus X_i \oplus \mathfrak{C}_{q} \cong {\bf 1}\oplus S\oplus {\bf 1} \oplus \mathfrak{L}_4\oplus \mathfrak{C}_k$;
\item\label{item : tom : addition : 3} $\mathbf{1} \oplus X_i \oplus \mathfrak{C}_{q} \cong {\bf 1} \oplus \mathfrak{L}_4\oplus \mathfrak{C}_k$.
\eroman

We will show that Conditions (\ref{item : tom : addition : 1}) and (\ref{item : tom : addition : 2}) lead to a contradiction and, therefore, do not hold. To this end, suppose first that Condition (\ref{item : tom : addition : 1}) holds. Since $q \geq 2$, we obtain that $\mathfrak{L}_s$ is either $\mathbf{1}$ or $\mathbf{1} \oplus \mathbf{1}$. Together with $\mathbf{1} \oplus X_i \oplus \mathfrak{C}_{q} \cong {\bf 1}\oplus S\oplus  \mathfrak{L}_s$ and $q \geq 2$, this yields that $X_i$ is isomorphic to a subposet of $S$, which is false because $S$ is simple and $X_i$ is not. Next we turn to proving that Condition (\ref{item : tom : addition : 2}) leads to a contradiction too. In this case, we have  $X_i \oplus \mathfrak{C}_{q}  \cong S\oplus {\bf 1} \oplus \mathfrak{L}_4\oplus \mathfrak{C}_k$. As a consequence, either $X_i$ is an upset of $S$ or $S$ is an upset of $X_i$. The former is impossible because $S$ is simple, while $X_i$ is not. Therefore, $S$ is an upset of $X_i$. Since $X_i$ is an upset of $\mathfrak{L}$ by assumption, we obtain that $S$ is a simple upset of $\mathfrak{L}$.  By looking at the structure of $\mathfrak{L}$, this implies that $S$ is either empty or $\mathbf{1}$ or $\mathbf{1} \oplus \mathbf{1}$ or $\mathbf{2} \oplus \mathbf{1}$. First, suppose that $S \in \{ \varnothing, \mathbf{1}, \mathbf{1} \oplus \mathbf{1}\}$. Together with $X_i \oplus \mathfrak{C}_{q}  \cong S\oplus {\bf 1} \oplus \mathfrak{L}_4\oplus \mathfrak{C}_k$ and $X_i \ne \varnothing$ (the latter because $X_i$ is not simple), this implies that $X_i$ has only one maximal element. Since $X_i$ is an upset of $\mathfrak{L}$, the inspection of the structure of $\mathfrak{L}$ yields that $X_i$ is simple, a contradiction. It only remains to consider the case where $S = \mathbf{2} \oplus \mathbf{1}$. In this case,  $X_i \oplus \mathfrak{C}_{q}  \cong S\oplus {\bf 1} \oplus \mathfrak{L}_4\oplus \mathfrak{C}_k  \cong \mathbf{2} \oplus \mathbf{1}\oplus {\bf 1} \oplus \mathfrak{L}_4\oplus \mathfrak{C}_k$. Since $S$ is an upset of $X_i$, this yields that $X_i  \cong \mathbf{2} \oplus \mathbf{1} \oplus Y$ for some poset $Y$. Together with the fact that $X_i$ is an upset of $\mathfrak{L}$, this yields that $Y  = \varnothing$ and, therefore, $X_i  \cong \mathbf{2} \oplus \mathbf{1}$, contradicting 
 the assumption that $X_i$ is not simple. 
 
This establishes that Conditions (\ref{item : tom : addition : 1}) and (\ref{item : tom : addition : 2}) do not hold. As a consequence, we obtain that Condition (\ref{item : tom : addition : 3}) holds, that is $\mathbf{1} \oplus X_i \oplus \mathfrak{C}_{q} \cong {\bf 1} \oplus \mathfrak{L}_4\oplus \mathfrak{C}_k$.
Therefore, 
$X_i \oplus \mathfrak{C}_{q}  \cong \mathfrak{L}_4 \oplus \mathfrak{C}_k$. Together with (\ref{Eq:Xi+1equalXmequals2}) this implies that $X  \cong {\bf 1}\oplus S_1\oplus \mathfrak{L}\oplus S_2\oplus \dots \oplus \mathfrak{L} \oplus S_p \oplus \mathfrak{L}_4\oplus \mathfrak{C}_k$.
\end{proof}

We thus established that  
there exists   $k \leq n$ such that
\[
X  \cong {\bf 1}\oplus S_1\oplus \mathfrak{L}\oplus S_2\oplus \dots \oplus \mathfrak{L} \oplus S_p \oplus \mathfrak{L}_4\oplus \mathfrak{C}_k.
\]
Then $\mathbf{1} \oplus \mathfrak{L} \oplus \mathfrak{L}_4 \oplus \mathfrak{C}_k$ is a continuous p-morphic image of $X$, which is  obtained 
by identifying 
${\bf 1}\oplus S_1\oplus \mathfrak{L}\oplus S_2\oplus \mathfrak{L} \oplus \dots \oplus S_{p-1}$ with the maximum of $X$  
and $S_p$ with the minimum of the remaining copy of $\mathfrak{L}$. 
Since $\mathfrak{G}_k = \mathbf{1} \oplus \mathfrak{L} \oplus \mathfrak{L}_4 \oplus \mathfrak{C}_k$, this implies that $\textup{Log}(X) \subseteq \textup{Log}(\mathfrak{G}_k)$. To see that the other inclusion also holds,
recall from Lemma~\ref{lem: nnnnn} that each $S_j\oplus \mathfrak{L}$ is a  continuous p-morphic image of $\mathfrak{L}$. 
In addition, if $S_p$ rooted and nonempty, then  $S_p$ is also a continuous p-morphic image of $\mathfrak{L}$. 

If $S_p$ is nonempty, but not rooted, then we note that since  
$\mathbf{1} \oplus S_p \oplus \mathfrak{L}_4 \oplus \mathfrak{C}_k$ is a continuous p-morphic image of $X$, we have  
$\mathbf{1} \oplus S_p \oplus \mathfrak{L}_4 \oplus \mathfrak{C}_k\in \mathcal{R}_n$, which is a contradiction by Corollary~\ref{cor: struct}. 
Therefore, if $S_p$ is nonempty, $X$ is a continuous p-morphic image of 
\[
Y \coloneqq {\bf 1} \oplus \underbrace{\mathfrak{L}\oplus \dots\oplus \mathfrak{L}}_{p-times}\oplus \mathfrak{L}_4\oplus \mathfrak{C}_k;
\]
and if $S_p$ is empty,  $X$ is a continuous p-morphic image of 
\[
Z \coloneqq {\bf 1} \oplus \underbrace{\mathfrak{L}\oplus \dots\oplus \mathfrak{L}}_{(p-1)-times}\oplus \mathfrak{L}_4\oplus \mathfrak{C}_k.
\]
By Lemma~\ref{lem: extr}(\ref{item:extr:2}), $\textup{Log}(Y) =  \textup{Log}(Z)$. Thus, we may  concentrate on the 
case where $X$ is a continuous p-morphic image of $Y$.

Because $X$ is a continuous p-morphic image of $Y$, we have  $\textup{Log}(Y) \subseteq \textup{Log}(X)$. Since at least one copy of $\mathfrak{L}$ appears as a summand of $X$ (because $X$ is infinite), and hence the same holds in the above decomposition of $Y$, we can apply 
Lemma~\ref{lem: extr}(\ref{item:extr:2}) to obtain that 
\[
\textup{Log}(Y) = \textup{Log}({\bf 1} \oplus  \mathfrak{L}\oplus \mathfrak{L}_4\oplus \mathfrak{C}_k) = \textup{Log}(\mathfrak{G}_k).
\]
Thus, $\textup{Log}(\mathfrak{G}_k) = \textup{Log}(Y)\subseteq \textup{Log}(X)$, and hence $\textup{Log}(X) = \textup{Log}(\mathfrak{G}_k)$. Since $k \leq n$, we conclude that 
\[
\pushQED{\qed}\textup{Log}(\mathcal{R}_n \cup \{ X \}) = \textup{Log}(\mathcal{R}_n \cup \{ \mathfrak{G}_k)) = \mathsf{L}_k.\qedhere \popQED
\]

As a consequence, we obtain our desired result: 

\begin{Theorem}\label{Thm: fin fmp span} 
The fmp span of $\textup{Log}(\mathfrak{G}_n)$ is $\{ \mathsf{L}_0, \dots, \mathsf{L}_n \}$ and its degree of fmp   is $n+1$.
\end{Theorem}

\begin{proof}
In view of Lemma \ref{Cor:n+1-logics-in-span}, in order to prove that
\[
\textup{fmp}(\textup{Log}(\mathfrak{G}_n)) = \{ \mathsf{L}_0, \dots, \mathsf{L}_n \} \text{ and } \textup{deg}(\textup{Log}(\mathfrak{G}_n)) = n+1,
\]
it suffices to show that if $\mathsf{L} \in \textup{fmp}(\textup{Log}(\mathfrak{G}_n))$, then $\mathsf{L} = \mathsf{L}_i$ for some $i \leq n$.

Let $\mathsf{L} \in \textup{fmp}(\textup{Log}(\mathfrak{G}_n))$. We first show that $\mathsf{L}$ is an extension of $\mathsf{KG}$. In view of Theorem \ref{Cor:KG-Jankov-axiom}, it suffices to prove that $\mathsf{L}$ contains all the Jankov formulas in $\mathsf{KG}$. Let $X$ be a finite rooted poset such that $\mathcal{J}(X) \in \mathsf{KG}$. Since $\mathsf{Fin}(\mathsf{L}) = \mathsf{Fin}(\textup{Log}(\mathfrak{G}_n))$, from Lemma~\ref{Lem:model-of-KG} it follows that 
$\mathsf{Fin}(\mathsf{L}) \vDash \mathsf{KG}$. Because $\mathcal{J}(X) \in \mathsf{KG}$, we obtain $\mathsf{Fin}(\mathsf{L})  \vDash \mathcal{J}(X)$. By the Fine Lemma, 
$X \notin \mathsf{Fin}(\mathsf{L})$. The application of the Dual Jankov Lemma now yields $\mathcal{J}(X) \in \mathsf{L}$, and hence $\mathsf{KG} \subseteq \mathsf{L}$. Therefore, we can invoke Theorem \ref{Thm:main-trick-n+1} to conclude that $\mathsf{L}= \mathsf{L}_i$ for some $i \leq n$. 
\end{proof}

As we mentioned earlier, this establishes the following. 

\begin{Theorem}\label{Prop: fin}
For each $1 \leq n < \aleph_0$, there exists an si-logic $\mathsf{L}$
with $\textup{deg}(\mathsf{L}) = n$. 
\end{Theorem}

To complete the proof of Theorem \ref{Thm:countable-part}, we require the following.

\begin{Theorem}\label{Prop: aleph}
There exists an si-logic $\mathsf{L}$ such that 
$\textup{deg}(\mathsf{L}) = \aleph_0$. 
\end{Theorem}

\begin{proof}
Define
\[
\mathcal{R} = \bigcup_{1\leq n }\mathcal{R}_n \, \, \text{ and } \, \, \mathsf{L}_0^\ast =  \textup{Log} (\mathcal{R}).
\] 
Clearly $\mathsf{L}_0^\ast$ is an si-logic. We show that its degree of fmp is $\aleph_0$. For this, consider the following extensions of $\mathsf{KG}$, where $n\geq 1$:
\begin{align*}
\mathsf{L}_n^\ast & =  \textup{Log} (\mathcal{R}\cup \{\mathfrak{G}_n\});\\
\mathsf{L}_\infty^\ast & =  \textup{Log} (\mathcal{R}\cup \{\mathfrak{G}_n: 1\leq n\}).
\end{align*}

A proof similar to that of
Lemma~\ref{Lem:VI-Rn-description} shows that
\begin{equation}\label{Eq:the:chain:of:logics:infty:case}
\mathsf{L}_\infty^\ast\subsetneq \dots\subsetneq \mathsf{L}_n^\ast \subsetneq \dots   \subsetneq  \mathsf{L}_1^\ast \subsetneq \mathsf{L}_0^\ast,
\end{equation}
and that all the logics in the above display belong to $\textup{fmp}(\mathsf{L}_0^\ast)$. Therefore, in order to prove that $\textup{deg}(\mathsf{L}_0^\ast) = \aleph_0$, it suffices to show that 
\begin{equation}\label{Eq:the:chain:of:logics:infty:case:2}
\textup{fmp}(\mathsf{L}_0^\ast) \subseteq \{ \mathsf{L}_n^\ast : n \geq 0 \} \cup  \{\mathsf{L}_\infty^\ast\}.
\end{equation}
 
An argument similar
to the one in the proof of Theorem~\ref{Thm:main-trick-n+1} shows that 
for each  extension $\mathsf{L}$ of $\mathsf{KG}$ such that $\mathsf{RFin}(\mathsf{L}) = \mathcal{R}$ and for each  $X\in \mathsf{FGR}(\mathsf{L})$,
we have $\textup{Log}(\mathcal{R}\cup \{ X \})= \mathsf{L}^\ast_j$ for some $j \geq 0$, where $j$ may possibly be $\infty$. This, by (\ref{Eq:the:chain:of:logics:infty:case}), implies that $\mathsf{L} = \mathsf{L}_p^\ast$ for
\[
p \coloneqq \max (\{ q \in \{ 0, 1, 2, \dots, \infty \} : \textup{Log}(\mathcal{R} \cup \{X\}) = {\sf L}_q^\ast \text{ for some }X \in \mathsf{FGR}(\mathsf{L}) \}).
\]
Therefore, every extension of $\mathsf{KG}$ in $\textup{fmp}(\mathsf{L}_0^\ast)$ belongs to $\{ \mathsf{L}_n^\ast : n \geq 0 \} \cup  \{\mathsf{L}_\infty^\ast\}$. As every member of $\textup{fmp}(\mathsf{L}_0^\ast)$ is an extension of $\mathsf{KG}$ (which can be shown as in the proof of Theorem~\ref{Thm: fin fmp span}), we conclude that  (\ref{Eq:the:chain:of:logics:infty:case:2}) holds. 
\end{proof}

\section{The continuum case} \label{sec: BWn}

In order to complete the proof of the antidichotomy theorem, it suffices to exhibit an si-logic  whose degree of fmp is $2^{\aleph_0}$. 
We will do this by proving the following:

\begin{Theorem}\label{Thm:LICS-continuum-case-3}
If $2 < n < \aleph_0$, then $\textup{deg}({\sf BW}_n) =  2^{\aleph_0}$.
\end{Theorem}

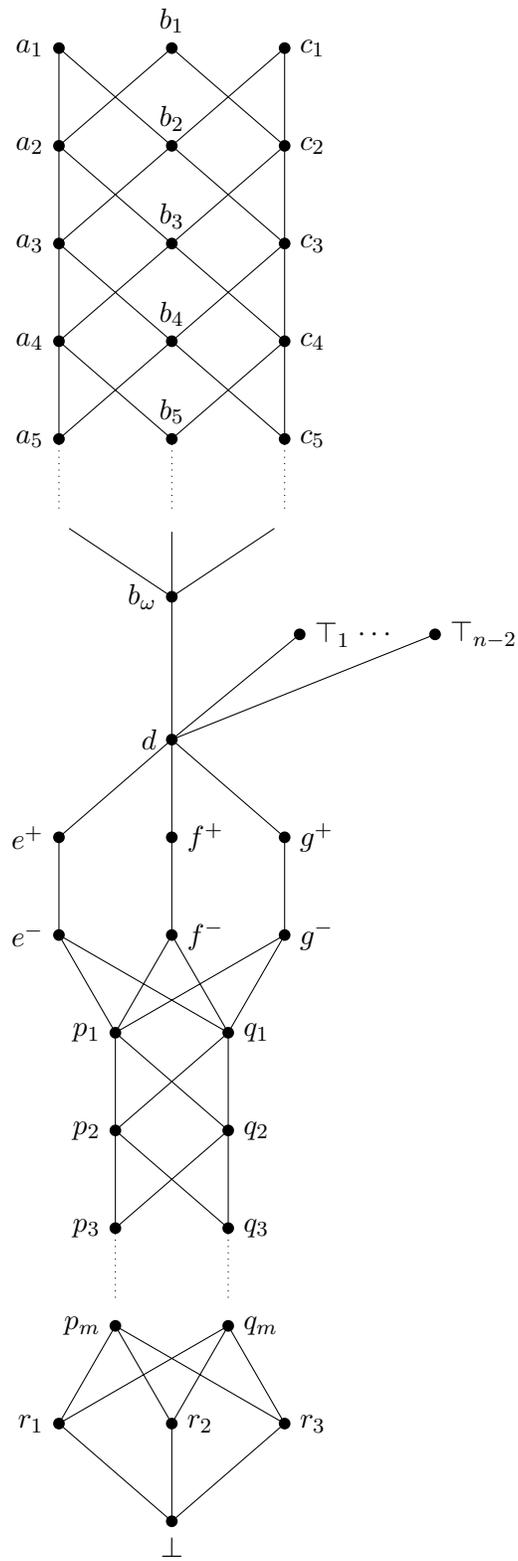
\begin{figure}
\[
\begin{tabular}{ccccccc}
\begin{tikzpicture}
    \tikzstyle{point} = [shape=circle, thick, draw=black, fill=black , scale=0.35]
    
\node[label=left:{$a_1$}] (a1) at (0.5,8) [point] {};
\node[label=above:{$b_1$}] (b1) at (2,8) [point] {};
\node[label=right:{$c_1$}] (c1) at (3.5,8) [point] {};
\node[label=left:{$a_2$}] (a2) at (0.5,6.7) [point] {};
\node[label=above:{$b_2$}] (b2) at (2,6.7) [point] {};
\node[label=right:{$c_2$}] (c2) at (3.5,6.7 ) [point] {};
\node[label=left:{$a_3$}] (a3) at (0.5,5.4) [point] {};
\node[label=above:{$b_3$}] (b3) at (2,5.4) [point] {};
\node[label=right:{$c_3$}] (c3) at (3.5,5.4) [point] {};
\node[label=left:{$a_4$}] (a4) at (0.5,4.1) [point] {};
\node[label=above:{$b_4$}] (b4) at (2,4.1) [point] {};
\node[label=right:{$c_4$}] (c4) at (3.5,4.1) [point] {};
\node[label=left:{$a_5$}] (a5) at (0.5,2.8) [point] {};
\node[label=above:{$b_5$}] (b5) at (2,2.8) [point] {};
\node[label=right:{$c_5$}] (c5) at (3.5,2.8) [point] {};

\node (a6) at (0.5,1.7) {};
\node (b6) at (2,1.7) {};
\node (c6) at (3.5,1.7) {};

\node[label=left:{$b_\omega$}] (bomega) at (2,0.7) [point] {};

\node[label=right:{$\top_1$}] (d1) at (3.7,0.2) [point] {};
\node[label=right:{$\cdots$}] (pp) at (4.2,0.2)  {};
\node[label=right:{$\top_{n-2}$}] (dn) at (5.5,0.2) [point] {};

\node[label=left:{$d$}] (bot) at (2,-1.2) [point] {};

\node[label=left:{$e^+$}] (e+) at (0.5,-2.5) [point] {};
\node[label=right:{$f^+$}] (f+) at (2,-2.5) [point] {};
\node[label=right:{$g^+$}] (g+) at (3.5,-2.5) [point] {};

\node[label=left:{$e^-$}] (e-) at (0.5,-3.8) [point] {};
\node[label=right:{$f^-$}] (f-) at (2,-3.8) [point] {};
\node[label=right:{$g^-$}] (g-) at (3.5,-3.8) [point] {};

\node[label=left:{$p_1$}] (p1) at (1.25,-5.1) [point] {};
\node[label=right:{$q_1$}] (q1) at (2.75,-5.1) [point] {};

\node[label=left:{$p_2$}] (p2) at (1.25,-6.4) [point] {};
\node[label=right:{$q_2$}] (q2) at (2.75,-6.4) [point] {};

\node[label=left:{$p_3$}] (p3) at (1.25,-7.7) [point] {};
\node[label=right:{$q_3$}] (q3) at (2.75,-7.7) [point] {};

\node (p4) at (1.25,-8.8) {};
\node (q4) at (2.75,-8.8) {};

\node[label=left:{$p_m$}] (pm) at (1.25,-9) [point] {};
\node[label=right:{$q_m$}] (qm) at (2.75,-9) [point] {};

\node[label=left:{$r_1$}] (r1) at (0.5,-10.3) [point] {};
\node[label=right:{$r_2$}] (r2) at (2,-10.3) [point] {};
\node[label=right:{$r_3$}] (r3) at (3.5,-10.3) [point] {};
\node[label=below:{$\bot$}] (bott) at (2,-11.6) [point] {};

    \draw[dotted] (p3) -- (p4) (q3) -- (q4) ;
    \draw (r1) -- (bott) -- (r2) (bott) -- (r3);
\draw (qm) -- (r1) -- (pm) -- (r2) -- (qm) -- (r3) -- (pm) ;
\draw (p3) -- (p2) -- (q3) -- (q2) -- (p3) ;
\draw (p1) -- (p2) -- (q1) -- (q2) -- (p1);
	\draw   (g-)  -- (q1) -- (f-) -- (p1) -- (g-) -- (g+) -- (bot) -- (e+) -- (e-)-- (p1) (q1) -- (e-) (bot) -- (f+) -- (f-);
    \draw  (b1) -- (c2) -- (c1) -- (b2) -- (a1) -- (a2) -- (b1);
    \draw  (b2) -- (c3) -- (c2) -- (b3) -- (a2) -- (a3) -- (b2);
    \draw  (b3) -- (c4) -- (c3) -- (b4) -- (a3) -- (a4) -- (b3);
    \draw  (b4) -- (c5) -- (c4) -- (b5) -- (a4) -- (a5) -- (b4);

    \draw[dotted] (a5) -- (a6) (b5) -- (b6) (c5) -- (c6) ;
    \draw (b6) -- (bomega)-- (a6) (bot) -- (bomega)-- (c6) ;
	\draw (d1) -- (bot) -- (dn);
    
\end{tikzpicture}
\end{tabular}
\]
\caption{The poset $X_m$.}
\label{Fig:bottom-part-W3}
\end{figure}

To establish the above result, let $2 < n < \aleph_0$ and let $\mathbb{Z}^+$ be the set of positive integers. For each $m \in \mathbb{Z}^+$, let $X_m$ be the poset in Figure \ref{Fig:bottom-part-W3} (we point out that  ${\uparrow}b_\omega$ is infinite and ${\downarrow}b_\omega$ is finite).\footnote{Posets similar to the upper part of $X_m$ have been considered in the literature (see, e.g., \cite[p.~319]{ChZa97}).}

We define a topology 
on $X_m$ by letting 
a subset $U$ of $X_m$ be open provided $b_{\omega} \in U$ implies $U$ is cofinite. Therefore, $b_\omega$ is the only limit point of $X_m$ and all the other points are isolated. It is routine to verify that this turns $X_m$ into an Esakia space, which we also denote by $X_m$. 

For each subset $M$ of $\mathbb{Z}^+$ let $\mathsf{L}_M$ be the si-logic of 
the class of Heyting algebras
\[
\{ X_{m}^{\ast} : m \in M \} \cup \{ \A \in \mathsf{W}_n : \A \text{ is finite} \}.
\] 
In order to prove Theorem \ref{Thm:LICS-continuum-case-3}, it suffices to establish the following.

\begin{Proposition}\label{Prop:LICS-intermediate-step-continuum}
The set $\{ \mathsf{L}_M : M \subseteq \mathbb{Z}^+ \}$ has the cardinality $2^{\aleph_0}$ 
and is a subset of $\textup{fmp}({\sf BW}_n)$.
\end{Proposition}

We split the proof of the above proposition in two parts, first showing that the cardinality of $\{ \mathsf{L}_M : M \subseteq \mathbb{Z}^+ \}$ is $2^{\aleph_0}$ and then that $\{ \mathsf{L}_M : M \subseteq \mathbb{Z}^+ \}$ is a subset of $\textup{fmp}({\sf BW}_n)$. 

In order to prove that the cardinality of $\{ \mathsf{L}_M : M \subseteq \mathbb{Z}^+ \}$ is $2^{\aleph_0}$, recall that $F_{n+1}$ is the poset depicted in Figure \ref{Fig:Wn-subframe}. We let $Y_m = {\downarrow}d$. 
It is enough to establish the following result.

\begin{Lemma}\label{Lem:continuum-LICS-1}
Let $M, N \subseteq \mathbb{Z}^+$ and $m \in M \smallsetminus N$. Then $\beta(F_{n+1}) \lor \mathcal{J}(Y_m) \in \mathsf{L}_N \smallsetminus \mathsf{L}_M$.
\end{Lemma}

\begin{proof}
We first show that $\beta(F_{n+1}) \lor \mathcal{J}(Y_m)$ is refuted on $X_m$. Since $X_m$ is rooted, it suffices to show that neither $\beta(F_{n+1})$ nor $\mathcal{J}(Y_m)$ is valid in $X_m$. First, as $X_m$ has width $n+1$ and $\beta(F_n)$ axiomatizes the si-logic $\mathsf{BW}_n$ of Esakia spaces of width $\leq n$, we obtain that $X_m \nvDash \beta(F_{n+1})$. On the other hand, when endowed with the discrete topology, $Y_m$ is a continuous p-morphic image of $X_m$ obtained by collapsing all the elements of ${\uparrow}d$ to $d$. By the Dual Jankov Lemma, this yields that $X_m \nvDash \mathcal{J}(Y_m)$. Thus, $\beta(F_{n+1}) \lor \mathcal{J}(Y_m)$ does not belong to~$\mathsf{L}_M$. 

It remains to prove that $\beta(F_{n+1}) \lor \mathcal{J}(Y_m)$ belongs to $\mathsf{L}_N$. By the definition of $\mathsf{L}_N$, it suffices to show that $\beta(F_{n+1}) \lor \mathcal{J}(Y_m)$ is valid in the class of algebras 
\[
\{ X_{k}^{\ast} : k \in N \} \cup \{ \A \in \mathsf{W}_n : \A \text{ is finite} \}.
\] 
Since $\beta(F_{n+1})$ axiomatizes ${\sf BW}_n$, the disjunction $\beta(F_{n+1}) \lor \mathcal{J}(Y_m)$ holds in the finite members of $\mathsf{W}_n$. Therefore, it only remains to show that  $\beta(F_{n+1}) \lor \mathcal{J}(Y_m)$ is valid in $\{ X_{k} : k \in N \}$. 

Suppose the contrary, with a view to contradiction. Then there exists $k \in N$ such that $X_k \nvDash \beta(F_{n+1}) \lor \mathcal{J}(Y_m)$. As a consequence,  $X_k \nvDash \mathcal{J}(Y_m)$. By the Dual Jankov Lemma there exist a closed upset $U$ of $X_k$ and a surjective continuous p-morphism $\alpha \colon U \to Y_m$ where the poset $Y_m$ is endowed with the discrete topology.

We will show that
\begin{equation}\label{Eq:why-Y-cannot-be-ontain-Wn-continuum}
\alpha^{-1}(\{ e^{+}, e^-, f^+, f^-, g^+, g^-\}) = \{ e^{+}, e^-, f^+, f^-, g^+, g^-\}.
\end{equation}
To prove the inclusion from left to right, consider $x \in\alpha^{-1}(\{ e^{+}, e^-, f^+, f^-, g^+, g^-\})$. By symmetry, we may assume that $\alpha(x) \in \{ e^+, e^- \}$. Since $\alpha \colon U \to Y_m$ is surjective, there are $x_f^+, x_f^-, x_g^+, x_g^- \in U$ such that
\[
 \alpha(x_f^+) = f^+ \quad \alpha(x_f^-) = f^- \quad \alpha(x_g^+) = g^+ \quad \alpha(x_g^-) = g^-.
\] 
Furthermore, as $\alpha$ is a p-morphism and $e^+, e^-, f^+, f^-, g^+$ and $g^-$ are not maximal,  the elements $x, x_f^+, x_f^-, x_g^+, x_g^-$ are also not maximal. Now, since $\alpha$ is order preserving, $\alpha(x) \in \{ e^+, e^- \}$, and $e^+, e^-$ are incomparable with $f^+, f^-, g^+$ and $g^-$, the element $x$ must be incomparable with $x_f^+, x_f^-, x_g^+$ and $x_g^-$. By the same token, $x_f^+, x_f^-$ are incomparable with $x_g^+, x_g^-$. In brief, $x$ is a nonmaximal element that is incomparable with four distinct nonmaximal elements $x_f^+, x_f^-, x_g^+$ and $x_g^-$ such that $x_f^+, x_f^-$ are incomparable with $x_g^+, x_g^-$. Examining the pictorial definition of $X_m$, it is easy to see that 
\[
x \in \{ e^{+}, e^-, f^+, f^-, g^+, g^- \}.
\]
This establishes the inclusion from left to right in (\ref{Eq:why-Y-cannot-be-ontain-Wn-continuum}). 

We prove the reverse inclusion by contradiction.
By symmetry, we may assume that there exists $x \in \{ e^{+}, e^-, f^+, f^-, g^+, g^- \}$ such that either $x \notin U$ or ($x \in U$ and $\alpha(x) = d$ or $\alpha(x) \leq p_1$). By symmetry, we may assume that $x \in \{ e^-, e^+ \}$. First suppose that $x \in U$ and $\alpha(x) = d$. The monotonicity of $\alpha$ implies that $\alpha(y) = d$ for each $y \in U$ such that $y \geq x$. Since $\alpha \colon U \to Y_m$ is surjective, this implies that the restriction $\alpha \colon (U\smallsetminus {\uparrow} x) \to (Y_m \smallsetminus \{ d \})$ is also surjective. But the assumption that $x \in \{ e^{+}, e^- \}$ implies that
\[
\vert U \smallsetminus {\uparrow} x\vert \leq \vert U \smallsetminus {\uparrow}e^+ \vert \leq  \vert Y_m \vert - 2 < \vert Y_m \smallsetminus \{ d \} \vert,
\]
a contradiction to the surjectivity of $\alpha \colon (U\smallsetminus {\uparrow} x) \to (Y_m \smallsetminus \{ d \})$. Therefore, $\alpha(x) \leq p_1 \leq e^-, f^-$. Since $\alpha$ is a p-morphism and $e^-, f^-$ are non-maximal, there must be two incomparable non-maximal elements $y, z \geq x$ such that $\alpha(y) = e^-$ and $\alpha(z) = f^-$. From $x \in \{ e^+, e^- \}$ it follows that $a_p \leq y, z$ for some $p \in \mathbb{Z}^+$. Moreover, since $\{f^+, f^-\}$ and $\{ g^+, g^- \}$ are two-element chains in $Y_m$ that are incomparable with each other and with $\alpha(y) = e^-$, the surjectivity of $\alpha$ implies the existence of two two-element chains in $U$ that are incomparable with each other as well as with $y$. But the fact that $y \geq a_p$ makes this impossible. Therefore, it only remains to consider the case when $x \notin U$. Since $U$ is an upset of $X_k$, from $x \in \{ e^+, e^- \}$ it follows that $U \subseteq  {\uparrow}\{ e^+, f^-, g^- \}$. Together with the assumption that $\alpha \colon U \to Y_m$ is a surjective p-morphism, this guarantees the existence of some $p \in \mathbb{Z}^+$ and $y \geq a_p$ such that $\alpha(y) = e^-$. This allows us to repeat the argument detailed above to obtain the desired contradiction. Thus, the inclusion from right to left in (\ref{Eq:why-Y-cannot-be-ontain-Wn-continuum}) holds.

Now, given a subset $V$ of $X_k$ and $W \subseteq V$, let
\[
{\uparrow^V} W = \{ x \in V : x \geq y \text{ for some }y \in W \} \, \, \text{ and }\, \, {\downarrow^V} W = \{ x \in V : x \leq y \text{ for some }y \in W \}.
\]
Let also $A = \{e^{+}, e^-, f^+, f^-, g^+, g^- \}$ and notice that $A \subseteq U$ by 
(\ref{Eq:why-Y-cannot-be-ontain-Wn-continuum}). Moreover, $U$ can be partitioned into the disjoint sets ${\downarrow}^{U}A$ and ${\uparrow}^{U}d$. Since $\alpha \colon U \to Y_m$ is order preserving, from 
(\ref{Eq:why-Y-cannot-be-ontain-Wn-continuum}) it follows that $\alpha({\downarrow}^{U}A) \subseteq {\downarrow}^{Y_m}A$. Because $d \in Y_m \smallsetminus {\downarrow}^{Y_m}A$ and $\alpha \colon U \to Y_m$ is surjective, we obtain that $d \in \alpha({\uparrow}^{U}d)$. We show that $\alpha({\uparrow}^{U}d) \subseteq \{ d \}$. Suppose the contrary. As $d$ is the maximum of $Y_m$ and $\alpha$ is order preserving, this implies $\alpha(d) < d$. By the definition of $Y_m$ we obtain that $\alpha(d) \in {\downarrow}^{Y_m}A$. By symmetry, we may assume that $\alpha(d) \leq e^+$. Moreover, recall from 
(\ref{Eq:why-Y-cannot-be-ontain-Wn-continuum}) that $\alpha(A) = A$. Let $x \in A$ be such that $\alpha(x) = f^+$. Since $x \in A$, we have $x \leq d$. Because $\alpha$ is order preserving, this implies $f^+ = \alpha(x) \leq \alpha(d) \leq e^+$, a contradiction. Therefore, we conclude that $\alpha({\uparrow}^{U}d) \subseteq \{ d \}$; that is,  $\alpha({\uparrow}^Ud) = \{ d \}$. 

In brief, $\alpha$ sends all elements of ${\downarrow}^U A$ to elements of $Y_m$ that are strictly less than $d$ and all elements of ${\uparrow}^U d$ to $d$. Since $d$ is an upper bound of $A$ in $U$, this implies that the restriction $\alpha \colon (\{ d \} \cup {\downarrow}^U A) \to Y_m$ is a surjective p-morphism. Because $\{ d \} \cup {\downarrow}^U A = Y_k \cap U$, we obtain that the map $\alpha \colon (Y_k \cap U) \to Y_m$ is also a surjective p-morphism. But, inspecting Figure \ref{Fig:bottom-part-W3}, it is easy to see that $m \ne k$ makes this impossible. Hence, we conclude that the disjunction $\beta(F_{n+1}) \lor \mathcal{J}(Y_m)$ is valid in $\{ X_k : k \in N \}$.
\end{proof}

The second part of Proposition \ref{Prop:LICS-intermediate-step-continuum} requires to prove that $\{ \mathsf{L}_M : M \subseteq \mathbb{Z}^+ \}$ is a subset of the fmp span of ${\sf BW}_n$, which amounts to the following.

\begin{Lemma}\label{Lem:LICS-continuum-technical-lem-2}
For each $M \subseteq \mathbb{Z}^+$ and finite poset  $X$,
\[
X \vDash {\sf BW}_n \, \, \mbox{ iff } \, \, X \vDash \mathsf{L}_M.
\]
\end{Lemma}

\begin{proof}
Let $M \subseteq \mathbb{Z}^+$. For the left to right implication, if $X \vDash {\sf BW}_n$, then $X$ has width $\leq n$, and hence $\mathsf{Up}(X)$ is a finite member of $\mathsf{W}_n$. Together with the definition of $\mathsf{L}_M$, this implies that $X \vDash \mathsf{L}_M$. 

For the right to left implication, it suffices to prove that if $X$ 
is a finite poset such that $X \vDash \mathsf{L}_M$, then 
$X \vDash \mathsf{BW}_n$. Suppose the contrary, with a view to contradiction. Then there exists a finite rooted poset $X$ of width $> n$ such that $X \vDash \mathsf{L}_M$. Since $X$ is finite and rooted, we can consider the Jankov formula $\mathcal{J}(X)$. Now, from $X \vDash \mathsf{L}_M$ it follows that $\mathsf{L}_M \nvDash\mathcal{J}(X)$. The definition of $\mathsf{L}_M$ implies that $\mathcal{J}(X)$ fails either in some $X_{m}$ with $m \in M$ or in some finite member of $\mathsf{W}_n$. Since $\mathsf{W}_n \vDash \mathcal{J}(X)$ because $\mathsf{Up}(X) \notin \mathsf{W}_n$, we conclude that there exists $m \in M$ such that $X_{m} \nvDash \mathcal{J}(X)$. Therefore, the Dual Jankov Lemma implies that there exist a closed upset $U$ of $X_m$ and an E-partition $R$ of $U$ such that $U / R$ is isomorphic to $X$. As $X$ is rooted, we may assume that $U$ is also rooted. 

Furthermore, as $X$ is not of width $\leq n$ and $U / R \cong X$, the set $U$ contains an $(n+1)$-element antichain. An inspection of the pictorial definition of $X_m$ shows that $U$ must contain an antichain of the form $\{\top_1, \dots, \top_{n-2}, a_k, b_k, c_k\}$ for some $k \in \mathbb{Z}^+$. Bearing in mind that $U$ is rooted, this implies that $U$ contains $d$ and, therefore, ${\uparrow}d \subseteq U$ because $U$ is an upset. In brief, $U$ is a rooted upset of $X_m$ such that ${\uparrow}d \subseteq U$ and $R$ is an E-partition of $U$ such that $U / R$ is finite and has an $(n+1)$-element antichain.

Examining again the pictorial definition of $X_m$, it is easy to see that there must be some $k \in \mathbb{Z}^+$ such that $\{ a_k, b_k, c_k \} \subseteq U$ and 
\begin{equation}\label{Eq:the-antichain-of-the-first-part-of-the-proof}
\{ [a_k], [b_k], [c_k], [\top_1], \dots, [\top_{n-2}] \}
\end{equation}
is an $(n+1)$-element antichain of $U /R$ (notice that $\top_1, \dots, \top_{n-2}\in U$ because ${\uparrow}d \subseteq U$). Consequently, $\{ [a_k], [b_k], [c_k]\}$ is a three-element antichain.

\begin{Claim}
There exists the  largest $j \in \mathbb{Z}^+$ such that $\{ [a_j], [b_j], [c_j] \}$ is a three-element antichain. 
\end{Claim}

\begin{proof}[Proof of the Claim.]
Suppose the contrary and recall that ${\uparrow}b_\omega \subseteq U$. We show that the equivalence class $[b_\omega]$ does not contain any 
$a_i, b_i,$ or $c_i$ for $i \in \mathbb{Z}^+$. 
If $[b_\omega]$ contains 
$x \geq b_\omega$, then it also contains the interval $[b_\omega, x]$. In particular, if $[b_\omega]$ contains 
$a_i, b_i$ or $c_i$, then $[b_{\omega}, a_{i+2}] \subseteq [b_\omega]$. This means that for each $t \geq i+2$ we have $[a_t] = [b_t] = [c_t] = [b_\omega]$. Hence,
\[
p \coloneqq \max \{ t \in \mathbb{Z}^+ : t < i+2 \text{ and }\{ [a_t], [b_t], [c_t] \}\text{ is a three-element antichain} \}
\]
exists (because $k \leq i+2$ and $\{ [a_k], [b_k], [c_k]\}$ is a three-element antichain) and is the largest positive integer $t$ such that $\{ [a_t], [b_t], [c_t] \}$ is a three-element antichain. The obtained contradiction proves 
that $[b_\omega]$ does not contain any of
$a_i, b_i$ or $c_i$. 

As a consequence,
\[
[b_\omega] \subseteq  \{b_\omega, \top_1, \dots, \top_{n-2}\} \cup {\downarrow}d,
\]
 whence $[b_\omega]$ is finite. Now, recall that $U / R \cong X$ and that the topology of $X$ is discrete because $X$ is finite. Therefore, $[b_\omega]$ is an isolated point of $U / R$. Since the map $x \mapsto [x]$ is a continuous p-morphism from $U$ to $U / R$, it follows that $[b_\omega]$ is a clopen subset of $U$. But since $U$ is a closed upset of $X_m$ containing ${\uparrow}b_\omega$,  the definition of the topology of $X_m$ guarantees that $[b_\omega]$ must contain infinitely many elements of ${\uparrow}b_\omega$. Therefore, $[b_\omega]$ is infinite, a contradiction. 
\end{proof}

Let $j$ be the largest positive integer such that $\{ [a_j], [b_j], [c_j] \}$ is a three element-antichain, which exists by the Claim. Then $\{ [a_{j+1}], [b_{j+1}], [c_{j+1}] \}$ is not a three element antichain. By symmetry, we may assume that
\[
[a_{j+1}] \leq [b_{j+1}].
\]
From $b_{j+1} \leq c_j$ it follows that $[b_{j+1}] \leq [c_j]$, and so $[a_{j+1}] \leq [c_j]$. Therefore, there exist $x \in [a_{j+1}]$ and $y \in [c_j]$ such that $x \leq y$. Since $[a_{j+1}] = [x]$ and $x \leq y$, the definition of an E-partition guarantees the existence of some $z \in U$ such that $a_{j+1} \leq z$ and $[z] = [y]$. Since $[y] = [c_{j}]$, we obtain 
\[
a_{j+1} \leq z\text{ and } [z] = [c_j].
\]
Notice that every element of ${\uparrow}a_{j+1}$ is comparable with $a_j$ or $b_j$.\ In particular, $z$ must be comparable with $a_j$ or $b_j$. Together with $[z] = [c_j]$, this implies that $[c_j]$ is comparable with $[a_j]$ or $[b_j]$. But this contradicts the assumption that $\{ [a_j], [b_j], [c_j] \}$ is a three element antichain.
\end{proof}

From Lemmas \ref{Lem:continuum-LICS-1} and \ref{Lem:LICS-continuum-technical-lem-2} it follows that Proposition \ref{Prop:LICS-intermediate-step-continuum} holds. Therefore, $\textup{deg}(\mathsf{BW}_n) = 2^{\aleph_0}$. Since we proved this equality for an arbitrary $2 < n < \aleph_0$, this establishes Theorem \ref{Thm:LICS-continuum-case-3}. Together with Theorem \ref{Thm:countable-part}, this concludes the proof of the Antidichotomy Theorem. 

We close this section with an observation about the logics \nolinebreak ${\sf BW}_n$. 

\begin{Proposition}\label{Prop:jankov-axiom-width-charact}
For each $n < \aleph_0$, the logic ${\sf BW}_n$ can be axiomatized by Jankov formulas iff $n \leq 2$.
\end{Proposition}

\begin{proof}
Since $\mathsf{BW}_0$ is the trivial variety, it is axiomatizable by the Jankov formula of the two-element Boolean algebra. 
Also, since $\mathsf{BW}_1$ is the G\"odel-Dummett logic \cite{Dm59}, 
it is well known that 
$\mathsf{BW}_1$ is axiomatizable by the Jankov formulas of the posets in Figure~\ref{fig:goedel-Jankov} (see, e.g., \cite[Thm.~4.23(4)]{BB22}).

\begin{figure}
\[
\begin{tabular}{ccccccc}
\begin{tikzpicture}
    \tikzstyle{point} = [shape=circle, thick, draw=black, fill=black , scale=0.35]
       \node (0) at (-2.5,-7) [point] {};    
       \node (a4) at (-3.5,-6) [point] {};
       \node (a9) at (-1.5,-6) [point] {};

    \draw  (a9) -- (0) -- (a4)    ;
    
      \tikzstyle{point} = [shape=circle, thick, draw=black, fill=black , scale=0.35]
       \node (00) at (1.5,-7) [point] {};
    \node (b2) at (1.5,-5) [point] {};
        \node (b4) at (0.5,-6) [point] {};
        \node (b9) at (2.5,-6) [point] {};

    \draw  (b9) -- (00) -- (b4) -- (b2) -- (b9)   ;
  
\end{tikzpicture}
\end{tabular}
\]
\caption{The two posets whose Jankov formulas axiomatize ${\sf BW}_1$.}
\label{fig:goedel-Jankov}
\end{figure}
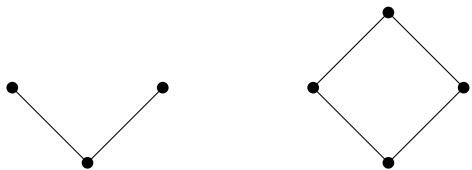

If $n=2$, then it follows from Theorem \ref{Thm:width2-kracht}
that ${\sf BW}_2$ is also axiomatizable by Jankov formulas. Finally, let $n > 2$. By Theorem \ref{Thm:subframe-axioms-for-Wn}, ${\sf BW}_n$ has the fmp. Since   ${\sf BW}_n$ has the degree of fmp $2^{\aleph_0}$ by Theorem \ref{Thm:LICS-continuum-case-3}, we can use Corollary \ref{Cor:degre-one} to deduce that $\mathsf{BW}_n$ cannot be axiomatized by Jankov formulas.
\end{proof}

In view of the Antidichotomy Theorem, every nonzero cardinal $\kappa$ such that $\kappa \leq \aleph_0$ or  $\kappa = 2^{\aleph_0}$ may occur as the degree of fmp of some si-logic. 
An opposite scenario appears if we restrict our attention to the logics ${\sf BW}_n$. 

\begin{Theorem}[\textbf{Width Dichotomy Theorem}]
For each $n<\aleph_0$, we have
\[
\textup{deg}({\sf BW}_n) = \left\{ \begin{array}{ll}
 1 & \text{if $n \leq 2$;}\\
 2^{\aleph_{0}} & \text{if $n > 2$.}
  \end{array} \right. 
 \]
\end{Theorem}

\begin{proof}
By Theorem~\ref{Thm:subframe-axioms-for-Wn}, each ${\sf BW}_n$ has the fmp. Therefore, from Proposition \ref{Prop:jankov-axiom-width-charact} and Corollary \ref{Cor:degre-one} it follows that if $n \leq 2$, then $\textup{deg}({\sf BW}_n) = 1$. The case where $n > 2$ is a consequence of Theorem \ref{Thm:LICS-continuum-case-3}.
\end{proof}

\section{Degrees of fmp for modal logics}

In this section we investigate the degree of fmp for normal extensions of some prominent modal logics. To this end, we denote the class of Kripke frames (resp.~finite Kripke frames) validating a normal modal logic $\mathsf{L}$ by $\mathsf{Fr}(\mathsf{L})$ (resp.~$\mathsf{Fin}(\mathsf{L})$). The \emph{degree of incompleteness} (resp.~\emph{the degree of fmp}) of $\mathsf{L}$ is the number of normal modal logics $\mathsf{L}'$ such that $\mathsf{Fr}(\mathsf{L}) = \mathsf{Fr}(\mathsf{L}')$ (resp.~$\mathsf{Fin}(\mathsf{L}) = \mathsf{Fin}(\mathsf{L}')$).  

As a consequence of the Blok Dichotomy Theorem, we obtain a dichotomy theorem for the degree of fmp of normal extensions of the basic modal logic \nolinebreak $\mathsf{K}$. 

\begin{Theorem}[\textbf{FMP Dichotomy Theorem}]\label{Thm:strict:dichotomy:modal}
The degree of the fmp of a normal modal logic $\mathsf{L}$ is either $1$ or $2^{\aleph_0}$.
\end{Theorem}

\begin{proof}
Let $\mathsf{L}$ be a normal modal logic. By Blok Dichotomy Theorem, its degree of incompleteness is either 1 or $2^{\aleph_0}$.\ First suppose that the degree of incompleteness of $\mathsf{L}$ is $2^{\aleph_0}$. Then there are $2^{\aleph_0}$ normal modal logics $\mathsf{L}'$ such that $\mathsf{Fr}(\mathsf{L}')=\mathsf{Fr}(\mathsf{L})$.
Since $\mathsf{Fin}(\mathsf{L})\subseteq \mathsf{Fr}(\mathsf{L})$, it follows that $\mathsf{Fin}(\mathsf{L}')=\mathsf{Fin}(\mathsf{L})$.
Thus, 
the degree of fmp of $\mathsf{L}$ is also $2^{\aleph_0}$.\color{black}

Next suppose that the degree of incompleteness of $\mathsf L$ is 1. Then 
$\mathsf{L}$ is a join-splitting logic (see, e.g.,  \cite[Thm.~10.59]{ChZa97}), and hence $\mathsf{L}$ has the fmp (see, e.g.,  \cite[Thm.~10.54]{ChZa97}). 

 Recall that a Kripke frame $X$ is said to be \textit{cycle free} if there is no path of length $n >0$ from a point of $X$ to itself.  For each finite rooted  cycle-free   Kripke frame $X$ we denote by $\mathcal{J}(X)$ an analogue of the Jankov formula in the language of modal logic \cite[p.~362]{ChZa97}. The join-splitting normal modal logics are precisely those axiomatized by formulas of the form $\mathcal{J}(X)$ where $X$ is a finite rooted cycle free Kripke frame (see, e.g.,  \cite[Thm.\ 10.53]{ChZa97}).
In particular, since $\mathsf{L}$ is a join-splitting logic, $\mathsf{L}=\mathsf{K}+\{\mathcal{J}(X_i) : i \in I\}$ for some set $\{ X_i : i\in I \}$ of finite rooted 
cycle free Kripke frames.  

Let $\mathsf{L}'$ be a normal modal logic such that $\mathsf{Fin}(\mathsf{L}) = \mathsf{Fin}(\mathsf{L}')$. We show that $\mathsf{L} = \mathsf{L}'$. Since $\mathsf{L}$ has the fmp, from  $\mathsf{Fin}(\mathsf{L}) = \mathsf{Fin}(\mathsf{L}')$ it follows that  $\mathsf{L}' \subseteq \mathsf{L}$. On the other hand, the modal analogue of Lemma~\ref{lem: auxiliary}(\ref{item: splittings}) 
yields that $\mathcal{J}(X_i)\in\mathsf{L}'$ for each $i\in I$. Therefore, 
$\mathsf{L} \subseteq \mathsf{L}'$, and hence $\mathsf{L} = \mathsf{L}'$. Thus, the dichotomy theorem holds for degrees of fmp. 
\end{proof}

However, the situation changes dramatically if we relativize the notion of the degree of fmp to stronger normal modal logics.
Following the terminology of \cite{ChZa97}, given a normal modal logic $\mathsf{L}$, let $\textup{Next}\,\mathsf{L}$ be the lattice of normal extensions of $\mathsf{L}$. 

\begin{Definition}
For a normal extension $\mathsf{L}$ of a normal modal logic $\mathsf{M}$, let
\begin{align*}
\textup{fmp}_{\mathsf{M}}(\mathsf{L}) & = \{ \mathsf{L}' \in \textup{Next}\,\mathsf{M} : \mathsf{Fin}(\mathsf{L}') = \mathsf{Fin}(\mathsf{L}) \};\\
\textup{deg}_{\mathsf{M}}(\mathsf{L}) & = \vert \textup{fmp}_{\mathsf{M}}(\mathsf{L}) \vert.
\end{align*}
\end{Definition}

Recall that the \textit{Grzegorczyk logic} $\mathsf{Grz}$ is the normal extension of $\mathsf{S4}$ 
by the formula
\[
\Box (\Box (p \to \Box p)  \to p) \to \Box p)
\]
(see, e.g., \cite[pp.\ 74 and 93]{ChZa97}).
\begin{Theorem}[\textbf{Modal Antidichotomy Theorem}]\label{Thm:modal:dichotomy:main}
Let $\mathsf{M} \subseteq \mathsf{Grz}$ be a normal modal logic with the fmp \color{black} such that $\mathsf{Grz}$ is a join-splitting in $\textup{Next}\,\mathsf{M}$. For each nonzero cardinal $\kappa$ such that $\kappa \leq \aleph_0$ or $\kappa = 2^{\aleph_0}$ there is a normal extension $\mathsf{L}$ of $\mathsf{M}$ with $\textup{deg}_{\mathsf{M}}({\sf L}) = \kappa$. 
\end{Theorem}

Before proving the Modal Antidichotomy Theorem, we point out that it holds for $\mathsf{S4}$ and $\mathsf{K4}$. For recall that $\mathsf{Grz}$ is a join-splitting in $\textup{Next}\,\mathsf{S4}$ \cite[Exmp.~1.11]{ChZaWoCh01} and that $\mathsf{S4}$ is a join-splitting in $\textup{Next}\,\mathsf{K4}$ \cite[Exmp.~10.48]{ChZa97}.  Consequently, $\mathsf{Grz}$ is also a join-splitting in $\textup{Next}\,\mathsf{K4}$. Since both $\mathsf{S4}$ and $\mathsf{K4}$ have the fmp, we obtain  that the modal antidichotomy theorem holds in both $\textup{Next}\,\mathsf{S4}$ and $\textup{Next}\,\mathsf{K4}$:

\begin{Corollary}\label{Cor:antidic:S4:K4}
For each nonzero cardinal $\kappa$ such that $\kappa \leq \aleph_0$ or $\kappa = 2^{\aleph_0}$ there is $\mathsf{L} \in \textup{Next}\,\mathsf{S4}$ with
\[
\textup{deg}_{\mathsf{K4}}({\sf L}) = \textup{deg}_{\mathsf{S4}}({\sf L}) = \kappa.
\]
 \end{Corollary}

\begin{Remark}
In particular, a normal extension of $\mathsf{K4}$ has the degree of fmp 1 iff it has the fmp and is axiomatizable by Jankov formulas. 
The proof of this result is analogous to Corollary~\ref{Thm:interval} since the machinery of Jankov formulas is available for 
$\mathsf{K4}$ (see, e.g., \cite[Ch.~9]{ChZa97}).  
Consequently, since each locally tabular normal extension of $\mathsf{K4}$ is axiomatizable by Jankov formulas, we obtain an analogue of Corollary~\ref{Cor:locally-finite}: the degree of fmp of locally tabular normal extensions of $\mathsf{K4}$ is $1$.
\end{Remark}

In order to prove the Modal Antidichotomy Theorem, we recall that the \textit{G\"odel translation}, associating   with \color{black} each intuitionistic formula $\varphi$ the modal formula $\varphi^t$, is defined recursively as follows:
\begin{align*}
p^t & = \Box p \text{ for each propositional variable }p\\
\bot^t &= \bot\\
(\chi \land \psi)^t & = \chi^t \land \psi^t\\
(\chi \lor \psi)^t & = \chi^t \lor \psi^t\\
(\chi \to \psi)^t & = \Box( \chi^t \to \psi^t).
\end{align*}
By \cite{McKT48}, for each intuitionistic formula $\varphi$, we have
\[
\varphi \in \mathsf{IPC} \, \, \text{ iff } \, \, \varphi^t \in \mathsf{S4}.
\]
Let $\mathsf{L}$ be an si-logic and $\mathsf{M}$ a normal extension of $\mathsf{S4}$. Following the standard terminology (see, e.g., \cite[Sec.~9.6]{ChZa97}), we say that $\mathsf{M}$ is a \textit{modal companion} of $\mathsf{L}$ provided for each intuitionistic formula $\varphi$, we have
\[
\varphi \in \mathsf{L}\, \,  \text{ iff }\, \, \varphi^t \in \mathsf{M}.
\]
Notably, each si-logic $\mathsf{L}$ has the least and greatest modal companions, denoted by $\tau(\mathsf{L})$ and $\sigma(\mathsf{L})$. 
For example, $\tau(\mathsf{IPC}) = \mathsf{S4}$ and $\sigma(\mathsf{IPC})=\mathsf{Grz}$. More generally, $\tau(\mathsf{L}) = \mathsf{S4} + \{\varphi^t: \varphi\in \mathsf{L}\}$ and $\sigma(\mathsf{L}) = \mathsf{Grz} + \{\varphi^t: \varphi\in \mathsf{L}\}$ (see, e.g., \cite[Sec.~9.6]{ChZa97}).
The latter is a consequence of an important result in modal logic, known as the Blok-Esakia theorem.

\begin{Theorem}[\textbf{Blok-Esakia Theorem}]
The map $\sigma \colon \textup{Ext}\,\mathsf{IPC} \to \textup{Next}\,\mathsf{Grz}$ is an 
isomorphism. 
\end{Theorem}

\begin{proof}
See Blok \cite{Blok76} and Esakia \cite{Esakia76,MR579150}.
\end{proof}

When dealing with the degree of fmp, the following observation will also be useful. 

\begin{Proposition}\label{Prop:same-finite-frames-Grz}
For every si-logic $\mathsf{L}$,
\[
\mathsf{Fin}(\mathsf{L}) = \mathsf{Fin}(\sigma(\mathsf{L})).
\]
\end{Proposition}

\begin{proof}
We recall from \cite[Cor.~3.5.10]{Esakia-book85} that $\mathsf{Fin}(\mathsf{Grz})$ is the class of all finite posets. This yields the 
result together with the fact that a finite poset validates an intuitionistic formula $\varphi$ iff the same poset, when viewed as a Kripke frame, validates the modal formula $\varphi^t$.
\end{proof}

\begin{Proposition}\label{Prop:anti-dichotomy-Grz}
For each nonzero cardinal $\kappa$ such that $\kappa \leq \aleph_0$ or $\kappa = 2^{\aleph_0}$ there is a normal extension $\mathsf{L}$ of $\mathsf{Grz}$ with $\textup{deg}_{\mathsf{Grz}}({\sf L}) =  \kappa$.
\end{Proposition}
\color{black}
\begin{proof}
Consider a nonzero cardinal $\kappa$ such that $\kappa \leq \aleph_0$ or $\kappa = 2^{\aleph_0}$. By the Antidichotomy Theorem, there is an si-logic $\mathsf{L}$ such that $\textup{deg}(\mathsf{L}) = \kappa$. To complete the proof, it suffices to show that $\textup{deg}_{\mathsf{Grz}}(\sigma(\mathsf{L})) = \kappa$.
Since $\sigma \colon \textup{Ext}\,\mathsf{IPC} \to \textup{Next}\,\mathsf{Grz}$ is an isomorphism, it is enough to prove that
\[
\textup{fmp}_{\mathsf{Grz}}(\sigma(\mathsf{L})) = \{ \sigma(\mathsf{L}') : \mathsf{L}' \in \textup{fmp}(\mathsf{L}) \}. 
\]

The inclusion from right to left is an immediate consequence of Proposition \ref{Prop:same-finite-frames-Grz}. Indeed, if $\mathsf{L}$ and $\mathsf{L}'$ share the class of finite posets, then $\sigma(\mathsf{L})$ and 
$\sigma(\mathsf{L}')$ also share the same class of finite Kripke frames by Proposition \ref{Prop:same-finite-frames-Grz}.  

To prove the other inclusion, let $\mathsf{S} \in \textup{fmp}_{\mathsf{Grz}}(\sigma(\mathsf{L}))$. By the Blok-Esakia Theorem, there is a unique $\mathsf{L}'\in\textup{Ext}\,\mathsf{IPC}$ such that $\mathsf{S} = \sigma(\mathsf{L}')$. Moreover, $\mathsf{Fin}(\mathsf{S}) = \mathsf{Fin}(\mathsf{L}')$ by Proposition \ref{Prop:same-finite-frames-Grz}. Therefore,
\[
\mathsf{Fin}(\mathsf{L}) = \mathsf{Fin}(\sigma(\mathsf{L})) = \mathsf{Fin}(\mathsf{S}) = \mathsf{Fin}(\mathsf{L}').
\]
Thus, $\mathsf{L}' \in \textup{fmp}(\mathsf{L})$, which together with $\sigma(\mathsf{L}') = \mathsf{S}$ yields that $\mathsf{S} \in \{ \sigma(\mathsf{L}') : \mathsf{L}' \in \textup{fmp}(\mathsf{L}) \}$. 
\end{proof}

We are now ready to prove the Modal Antidichotomy Theorem.

\begin{proof}
Consider a nonzero cardinal $\kappa$ such that $\kappa \leq \aleph_0$ or $\kappa = 2^{\aleph_0}$. By Proposition \ref{Prop:anti-dichotomy-Grz}, there exists a normal extension $\mathsf{L}$ of $\mathsf{Grz}$ such that $\textup{deg}_{\mathsf{Grz}}(\mathsf{L}) = \kappa$. Therefore, to conclude the proof it suffices to show that
\[
\textup{fmp}_{\mathsf{M}}(\mathsf{L}) = \textup{fmp}_{\mathsf{Grz}}(\mathsf{L}).
\]

The inclusion from right to left is obvious because $\mathsf{M} \subseteq \mathsf{Grz}$ by assumption. To prove the other inclusion, consider $\mathsf{L}' \in \textup{fmp}_{\mathsf{M}}(\mathsf{L})$. Since $\mathsf{Fin}(\mathsf{L}') = \mathsf{Fin}(\mathsf{L})$, in order to prove that $\mathsf{L}' \in \textup{fmp}_{\mathsf{Grz}}(\mathsf{L})$ it remains to show that $\mathsf{L}'$ extends $\mathsf{Grz}$.

Since $\mathsf{M}$ has the fmp and $\mathsf{Grz}$ is a join-splitting in $\textup{Next}\,\mathsf{M}$, it follows from a general result of McKenzie \cite[Thm.~4.3]{MR313141} 
that there exists a set $K$ of finite rooted Kripke frames validating $\mathsf{M}$ such that $\mathsf{Grz}$ is the least normal extension $\mathsf{N}$ of $\mathsf{M}$ with $\mathsf{Fin}(\mathsf{N}) \cap K = \varnothing$. In particular, since $\mathsf{L}$ is an extension of $\mathsf{Grz}$, we have $\mathsf{Fin}(\mathsf{L}) \cap K = \varnothing$. Together with the assumption that $\mathsf{Fin}(\mathsf{L}) = \mathsf{Fin}(\mathsf{L}')$, this yields that $\mathsf{Fin}(\mathsf{L}') \cap K = \varnothing$. As $\mathsf{L}'$ is a normal extension of $\mathsf{M}$, we conclude that $\mathsf{Grz} \subseteq \mathsf{L}'$ as desired.
\end{proof}

\section{Conclusions}

In this paper, we introduced the notion of the degree of fmp for superintuionistic and modal logics in analogy with the classic notion of the degree of incompleteness for these logics. 
We proved the Antidichitomy Theorem for the degree of fmp for superintuionistic and transitive modal logics. Namely, for every nonzero cardinal $\kappa$ such that $\kappa\leq \aleph_0$ or $\kappa = 2^{\aleph_0}$ there is a  
superintuitionistic or transitive modal logic 
$\mathsf L$ such that the degree of fmp of $\mathsf L$ is $\kappa$. We conclude by discussing possible future research directions that could originate from this work. 

\begin{enumerate}
\item 
By assuming the Continuum Hypothesis (CH), our results show that the degree of fmp of si-logics can be any positive integer, $\aleph_0$, or $2^{\aleph_0}$. However, for proving this result, the power of CH may not be necessary. We leave it as an open problem whether the assumption of CH can be avoided. 

\item In this paper we determined what cardinalities can be realized as degrees of fmp for superintuitionistic and modal logics. However, it still remains an open problem to characterize the degree of fmp of a given si-logic or an extension of $\mathsf{S4}$ or $\mathsf{K4}$. Note that for an extension $\mathsf{L}$  of $\mathsf{K}$, the degree of fmp, as well as the degree of incompleteness of $\mathsf{L}$, is determined by the Blok Dichotomy Theorem: If $\mathsf{L}$ is join-splitting, then it is $1$; otherwise it is $2^{\aleph_0}$. In analogy with this, we showed that if a logic $\mathsf{L}$ has the fmp and is join-splitting, then its degree of fmp is $1$. But if these conditions are not met, the exact degree of fmp of $\mathsf{L}$ is still unclear. 

\item The first step to answer (2) would be to determine the degree of fmp of a given extension of $\mathsf{KG}$. In particular, it is still unclear whether the continuum degree of fmp can occur above $\mathsf{KG}$. 

\item We also find it interesting to study the degree of fmp for other prominent deductive systems \color{black} such as bi-intuionistic logic, tense and temporal  logics, and fixpoint logics such as $\mathsf{PDL}$ and the modal $\mu$-calculus. In fact, one can define and investigate the degree of fmp for any logic (or a variety of algebras thereof)
that has finite models. 
\end{enumerate}

More generally, one can apply this perspective to other logically interesting properties. 
For a given logic $\mathsf{L}$, let $\mathcal{S}$  be a semantics of $\mathsf{L}$ (relational, topological, algebraic, etc.). 
For a property $P$, 
 the {\em $P$-degree of the $\mathcal{S}$-semantics} is the cardinality of the set of 
logics $\mathsf{L}'$ such that $\mathsf{L}$ and $\mathsf{L}'$ share the same class of $\mathcal{S}$-models satisfying property $P$. The degree of fmp is then the $P$-degree of the $\mathcal{S}$-semantics when the $\mathcal{S}$-semantics is Kripke semantics and $P$ is the property of being finite. Note that being finite can be replaced by other properties; for example, by being countable, etc.

Since every si-logic or modal logic $\mathsf{L}$ is complete with respect to its algebraic semantics, the $P$-degree of the $\mathcal{S}$-semantics of each $\mathsf{L}$ is  always 1 when the $\mathcal{S}$-semantics is  the algebraic semantics and $P$ is any property true in each algebraic model. Indeed, in this case two logics have the same $P$-degree if they have the same class of algebraic models. Every such logic is complete with respect to its algebraic models. Hence, every logic  has the $P$-degree $1$. 
However, if $\mathcal{S}$-semantics is the  topological semantics, then the situation changes drastically since it is well known that there exist topologically incomplete modal logics (see, e.g., \cite{MR1688513,She05}) 
 and it remains an outstanding open problem whether there exist topologically incomplete si-logics. In a recent paper \cite{MR4520555}, it was shown that there exist (continuum many) extensions of the bi-intuitionistic logic that are topologically incomplete. 

In topological semantics of modal logic, it is customary to interpret $\Diamond$ as topological closure. Under such interpretation, $\mathsf{S4}$ is the least topologically complete modal logic, and the degree of topological fmp coincides with the degree of fmp in $\textup{NExt}\,\mathsf{S4}$ since finite topological spaces are in one-to-one correspondence with finite $\mathsf{S4}$-frames. On the other hand, if we interpret $\Diamond$ as topological derivative \cite[Appendix~I]{McKT44}
(the so-called {\em d-semantics}; see \cite{MR2195500}), 
then 
it makes sense to investigate  
the degree of topological fmp (which modal logics have the same class of finite topological models).
 
In our opinion, the study of $P$-degrees of $\mathcal{S}$-semantics for non-classical logics is a promising direction for future research.

\vspace{3mm}

\paragraph{\bfseries Acknowledgements.}

We are very grateful to the referee for  careful reading and useful comments. 
The authors acknowledge the support of the MSCA-RISE-Marie Skłodowska-Curie Research and Innovation Staff Exchange (RISE) project MOSAIC 101007627 funded by Horizon 2020 of the European Union. Part of this work was conducted during the visit of the third author to the Institute for Logic, Language and Computation of the University of Amsterdam, supported by the Dutch NWO visitor's grant. The third author was supported by the
proyecto PID2022-141529NB-C21 de investigaci\'on financiado por MICIU/AEI/ 10.13039/501100
011033 y por FEDER, UE. He was also supported by the Research Group in Mathematical Logic,
2021SGR00348 funded by the Agency for Management of University and Research Grants of
the Government of Catalonia.

\bibliographystyle{plain}

\begin{thebibliography}{10}

\bibitem{Ab05b}
S.~Abramsky.
\newblock A {C}ook's tour of the finitary non-well-founded sets.
\newblock Invited Lecture at BCTCS. Available at arXiv:1111.7148., 1988.

\bibitem{BB22}
G.~Bezhanishvili and N.~Bezhanishvili.
\newblock Jankov formulas and axiomatization techniques for intermediate
  logics.
\newblock In {\em V. {A}. {Y}ankov on non-classical logics, history and
  philosophy of mathematics}, volume~24 of {\em Outst. Contrib. Log.}, pages
  71--124. Springer, Cham, [2022] \copyright 2022.

\bibitem{BBdeJ08}
G.~Bezhanishvili, N.~Bezhanishvili, and D.~de~Jongh.
\newblock The {K}uznetsov-{G}erciu and {R}ieger-{N}ishimura logics: the
  boundaries of the finite model property.
\newblock {\em Logic and Logical Philosophy}, 17:73--110, 2008.

\bibitem{MR2195500}
G.~Bezhanishvili, L.~Esakia, and D.~Gabelaia.
\newblock Some results on modal axiomatization and definability for topological
  spaces.
\newblock {\em Studia Logica}, 81(3):325--355, 2005.

\bibitem{MR4520555}
G.~Bezhanishvili, D.~Gabelaia, and M.~Jibladze.
\newblock A negative solution of {K}uznetsov's problem for varieties of
  bi-{H}eyting algebras.
\newblock {\em J. Math. Log.}, 22(3):Paper No. 2250013, 21, 2022.

\bibitem{Bez-PhD}
N.~Bezhanishvili.
\newblock {\em Lattices of Intermediate and Cylindric Modal Logics}.
\newblock PhD thesis, University of Amsterdam, 2006.
\newblock
  \url{https://eprints.illc.uva.nl/id/eprint/2049/1/DS-2006-02.text.pdf}.

\bibitem{BezMor19}
N.~Bezhanishvili and T.~Moraschini.
\newblock {C}itkin's description of hereditarily structurally complete
  intermediate logics via {E}sakia duality.
\newblock {\em {S}tudia {L}ogica}, 111:174--186, 2023.

\bibitem{Blok76}
W.~J. Blok.
\newblock {\em Varieties of Interior Algebras}.
\newblock PhD thesis, University of Amsterdam, 1976.

\bibitem{B78}
W.~J. Blok.
\newblock On the degree of incompleteness of modal logics.
\newblock {\em Bulletin of the Section of Logic}, 7:167--175, 1978.

\bibitem{Blok-dichotomy78}
W.~J. Blok.
\newblock On the degree of incompleteness of modal logics and the covering
  relation in the lattice of modal logics.
\newblock Technical Report 78--07, Department of Mathematics, University of
  Amsterdam, 1978.

\bibitem{ChZa97}
A.~Chagrov and M.~Zakharyaschev.
\newblock {\em Modal Logic}, volume~35 of {\em Oxford Logic Guides}.
\newblock Oxford University Press, 1997.

\bibitem{MR1688513}
L.~Chagrova.
\newblock On the degree of neighborhood incompleteness of normal modal logics.
\newblock In {\em Advances in modal logic, {V}ol.\ 1 ({B}erlin, 1996)}, CSLI Lecture Notes, pages
  63--72. 1998.

\bibitem{MR0197372}
D.~H.~J. de~Jongh and A.~S. Troelstra.
\newblock On the connection of partially ordered sets with some
  pseudo-{B}oolean algebras.
\newblock {\em Indag. Math.}, pages
  317--329, 1966.

\bibitem{Dm59}
M.~Dummett.
\newblock A propositional calculus with denumerable matrix.
\newblock {\em J. Symbolic Logic}, 24:97--106, 1959.

\bibitem{Es74}
L.~Esakia.
\newblock Topological {K}ripke models.
\newblock {\em Soviet Math. Dokl.}, 15:147--151, 1974.

\bibitem{Esakia76}
L.~Esakia.
\newblock On modal ``companions" of superintuitionistic logics.
\newblock In {\em VII Soviet Symposium on Logic (Russian) (Kiev, 1976)}, pages
  135--136. 1976.

\bibitem{MR579150}
L.~Esakia.
\newblock On the variety of {G}rzegorczyk algebras.
\newblock In {\em Studies in nonclassical logics and set theory ({R}ussian)},
  pages 257--287. ``Nauka'', Moscow, 1979.

\bibitem{Esakia-book85}
L.~Esakia.
\newblock {\em {H}eyting {A}lgebras. {D}uality {T}heory}.
\newblock Springer, English translation of the original 1985 book. 2019.

\bibitem{Fin74}
K.~Fine.
\newblock An incomplete logic containing {${\rm S}4$}.
\newblock {\em Theoria}, 40:23--29, 1974.

\bibitem{Fin74b}
K.~Fine.
\newblock Logics containing {$K4$}. {P}art {I}.
\newblock {\em J. Symbolic Logic}, 34:31--42, 1974.

\bibitem{Fine85}
K.~Fine.
\newblock Logics containing {${\rm K}4$}. {P}art {II}.
\newblock {\em J. Symbolic Logic}, 50(3):619--651, 1985.

\bibitem{GeKuz70}
V.~Ja. Ger\v{c}iu and A.~V. Kuznetsov.
\newblock The finitely axiomatizable superintuitionistic logics.
\newblock {\em Soviet Math. Dokl.}, 11:1654--1658, 1970.

\bibitem{Jankov63for}
V.~A. Jankov.
\newblock On the relation between deducibility in intuitionistic propositional
  calculus and finite implicative structures.
\newblock {\em Doklady Akademii Nauk SSSR}, 151:1293--1294, 1963.
\newblock (In Russian).

\bibitem{Jankov68}
V.~A. Jankov.
\newblock The construction of a sequence of strongly independent
  superintuitionistic propositional calculi.
\newblock {\em Soviet Math. Dokl.}, 9:806--807, 1968.

\bibitem{Jankov69}
V.~A. Jankov.
\newblock Conjunctively irresolvable formulae in propositional calculi.
\newblock {\em Izvestiya Akademii Nauk SSSR. Seriya Matematicheskaya},
  33:18--38, 1969.
\newblock (In Russian).

\bibitem{Kracht93MLQ}
M.~Kracht.
\newblock Prefinitely axiomatizable modal and intermediate logics.
\newblock {\em {M}ath. {L}ogic {Q}uaterly}, 39:301--322, 1993.

\bibitem{Kr93b}
M.~Kracht.
\newblock Splittings and the finite model property.
\newblock {\em J. Symbolic Logic}, 58(1):139--157, 1993.

\bibitem{KupKuVe04}
C.~Kupke, A.~Kurz, and Y.~Venema.
\newblock Stone coalgebras.
\newblock {\em Theoretical Computer Science}, 327(1-2):109--134, 2004.

\bibitem{KuzGer70a}
A.~V. Kuznetsov and V.~Ja. Ger\v{c}iu.
\newblock Superintuitionistic logics and finite approximability.
\newblock {\em Soviet Math. Dokl.}, 11:1614--1619, 1970.

\bibitem{Lit08}
T.~Litak.
\newblock Stability of the {B}lok theorem.
\newblock {\em Algebra Universalis}, 58(4):385--411, 2008.

\bibitem{MR313141}
R.~McKenzie.
\newblock Equational bases and nonmodular lattice varieties.
\newblock {\em Trans. Amer. Math. Soc.}, 174:1--43, 1972.

\bibitem{McKT44}
J.~C.~C. McKinsey and A.~Tarski.
\newblock The algebra of topology.
\newblock {\em Ann. of Math.}, 45:141--191, 1944.

\bibitem{McKT48}
J.~C.~C. McKinsey and A.~Tarski.
\newblock Some theorems about the sentential calculi of {L}ewis and {H}eyting.
\newblock {\em J. Symbolic Logic}, 13:1--15, 1948.

\bibitem{Ni60}
I.~{}Nishimura.
\newblock On formulas of one variable in intuitionistic propositional calculus.
\newblock {\em J. Symbolic Logic}, 25:327--331 (1962), 1960.

\bibitem{RZW06}
W.~Rautenberg, M.~Zakharyaschev, and F.~Wolter.
\newblock Willem {B}lok and modal logic.
\newblock {\em Studia Logica}, 83(1-3):15--30, 2006.

\bibitem{Ri49}
L.~Rieger.
\newblock On the lattice theory of {B}rouwerian propositional logic.
\newblock {\em Acta Fac. Nat. Univ. Carol., Prague}, 1949(189, 40 pages), 1949.

\bibitem{She05}
V.~Shehtman.
\newblock On neighbourhood semantics 30 years later.
\newblock In S.~N. Artemov, H.~Barringer, A.~S. d'Avila Garcez, L.~C.`Lamb, and
  J.~Woods, editors, {\em We Will Show Them! Essays in Honour of Dov Gabbay,
  Volume Two}, pages 663--692. College Publications, 2005.

\bibitem{So77}
S.~K. Sobolev.
\newblock On finite-dimensional superintuitionistic logics.
\newblock {\em Izv. Akad. Nauk SSSR Ser. Mat.}, 41(5):963--986, 1977.
\newblock (In Russian).

\bibitem{Zakha89}
M.~Zakharyaschev.
\newblock Syntax and semantics of superintuitionistic logics.
\newblock {\em {A}lgebra and {L}ogic}, 28(4):262--282, 1989.

\bibitem{ChZaWoCh01}
M.~Zakharyaschev, F.~Wolter, and A.~Chagrov.
\newblock {\em Advanced modal logic}, volume~3 of {\em Handbook of
  Philosophical Logic}.
\newblock Springer, 2001.

\end{thebibliography}

\section*{Appendix} 

The aim of the Appendix is to prove Theorem~\ref{Cor:KG-Jankov-axiom} that $\mathsf{KG}$ is axiomatizable by Jankov formulas. For this we utilize the following classic result of Fine \cite{Fin74b}. For $n \geq 1$ we say that an si-logic $\mathsf{L}$ is of {\em width} $\leq n$ if each Esakia space $X$ validating $\mathsf{L}$ is of width $\leq n$ (see Definition~\ref{def: width} for the definition of the width of a poset). We call $\mathsf{L}$ of {\em finite width} if there is $n$ such that $\mathsf{L}$ is of width $\leq n$. Clearly $\mathsf{L}$ is of finite width provided ${\sf BW}_n\subseteq\mathsf{L}$ for some $n$. 

We recall that a poset is \textit{Noetherian} if it has no infinite strictly ascending chains. We then have (see, e.g., \cite[Thm.~10.45]{ChZa97}):

\begin{Theorem}[\textbf{Fine Completeness Theorem}]
If $\mathsf{L}$ is an si-logic of width $\leq n$, then
there is a class $K$ of rooted Noetherian posets of width $\leq n$ such that $\mathsf{L} = \textup{Log}(K)$.
\end{Theorem}

The proof of Theorem~\ref{Cor:KG-Jankov-axiom} is based on the following two combinatorial observations. 

\begin{figure}
\[ 
\begin{tabular}{ccccccc}
\begin{tikzpicture}
    \tikzstyle{point} = [shape=circle, thick, draw=black, fill=black , scale=0.35] 
  
  \node[label=above:{$K_1$}]  at (1.4,-1)  {};
\node (a1) at (1.4,0) [point] {};
\node (a2) at (0.7,0.7) [point] {};
\node (a3) at (2.1,0.7) [point] {};
\node (a4) at (0.7,1.4) [point] {};
\node (a5) at (2.1,1.4) [point] {};

  \node[label=above:{$K_2$}]  at (3.5,-1)  {};
\node (b1) at (3.5,0) [point] {};
\node (b2) at (2.8,0.7) [point] {};
\node (b3) at (4.2,0.7) [point] {};
\node (b4) at (2.8,1.4) [point] {};
\node (b5) at (4.2,1.4) [point] {};
\node (b6) at (4.2,2.1) [point] {};

  \node[label=above:{$K_3$}]  at (5.6,-1)  {};

\node (c1) at (5.6,0) [point] {};
\node (c2) at (4.9,0.7) [point] {};
\node (c3) at (6.3,0.7) [point] {};
\node (c4) at (4.9,1.4) [point] {};
\node (c5) at (6.3,1.4) [point] {};
\node (c6) at (5.6,2.1) [point] {};

  \node[label=above:{$K_4$}]  at (7.7,-1)  {};

\node (d1) at (7.7,0) [point] {};
\node (d2) at (7,0.7) [point] {};
\node (d3) at (8.4,0.7) [point] {};
\node (d4) at (7,1.4) [point] {};
\node (d5) at (8.4,1.4) [point] {};
\node (d6) at (8.4,2.1) [point] {};
\node (d7) at (7.7,2.8) [point] {};

  \node[label=above:{$K_5$}]  at (9.8,-1)  {};

\node (e1) at (9.8,0) [point] {}; 
\node (e2) at (9.1,0.7) [point] {};
\node (e3) at (10.5,0.7) [point] {};
\node (e4) at (9.1,1.4) [point] {};
\node (e5) at (10.5,1.4) [point] {};
\node (e6) at (10.5,2.1) [point] {};

  \node[label=above:{$K_6$}]  at (11.9,-1)  {};

\node (f1) at (11.9,0) [point] {}; 
\node (f2) at (11.2,0.7) [point] {};
\node (f3) at (12.6,0.7) [point] {};
\node (f4) at (11.2,1.4) [point] {};
\node (f5) at (12.6,1.4) [point] {};
\node (f6) at (12.6,2.1) [point] {};
\node (f7) at (11.2,2.1) [point] {};

  \node[label=above:{$K_7$}]  at (14,-1)  {};

\node (g1) at (14,0) [point] {};
\node (g2) at (13.3,0.7) [point] {};
\node (g3) at (14.7,0.7) [point] {};
\node (g4) at (13.3,1.4) [point] {};
\node (g5) at (14.7,1.4) [point] {};
\node (g6) at (14.7,2.1) [point] {};
\node (g7) at (14.7,2.8) [point] {};

\draw (a4) -- (a2) -- (a1) -- (a3) -- (a5) ;
  \draw (b4) -- (b2) -- (b1) -- (b3) -- (b5) -- (b6) (b3) -- (b4) ;
  \draw (c4) -- (c2) -- (c1) -- (c3) -- (c5)-- (c6) -- (c4) ;
  \draw (d6) -- (d7) -- (d4) -- (d2) -- (d1) -- (d3) -- (d5) -- (d6) (d3) -- (d4) ;
    \draw (e4) -- (e2) -- (e1) -- (e3) -- (e5) (e2) -- (e6) -- (e5) ;
  \draw (f7) -- (f4) -- (f2) -- (f1) -- (f3) -- (f5)-- (f2) (f5) -- (f6) -- (f4) ;
  \draw (g4) -- (g2) -- (g1) -- (g3) -- (g5) -- (g6)-- (g7) -- (g2) (g3) -- (g4) ;
\end{tikzpicture}
\end{tabular}
\]
\caption{The posets $K_1, \dots, K_7$.}
\label{Fig:posets-K}
\end{figure}
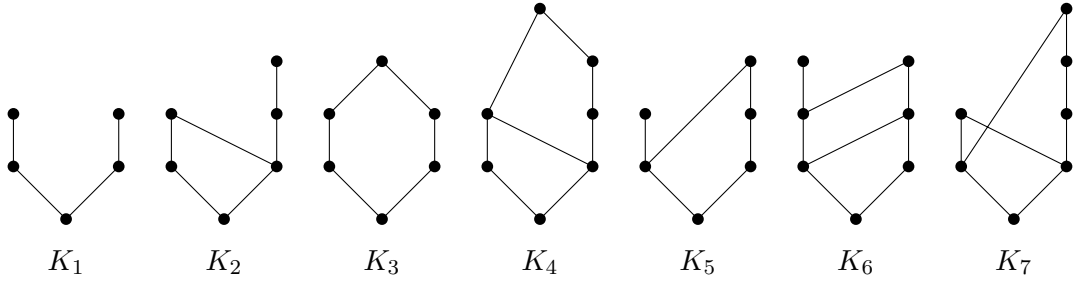

\begin{Lemma}\label{Lem:jankov-axioms-KG}
A rooted Noetherian poset $X$ of width $\leq 2$ validates $\beta(P_2)$ iff it validates the Jankov formulas of the posets in Figure \ref{Fig:posets-K}.
\end{Lemma}

\begin{figure}
  
\[
\begin{tabular}{ccccccc}
\begin{tikzpicture}
    \tikzstyle{point} = [shape=circle, thick, draw=black, fill=black , scale=0.35]
          
\node[label=above:{$G_1$}]  at (1.4,-1)  {};
\node (k1) at (1.4,0) [point] {};
\node (k2) at (0.7,0.7) [point] {};
\node (k3) at (2.1,0.7) [point] {};
\node (k5) at (2.1,1.4) [point] {};
\node (k6) at (2.1,2.1) [point] {};
    \draw  (k2) -- (k1) -- (k3) -- (k5) -- (k6); 
   
\node[label=above:{$G_2$}]  at (3.5,-1)  {};   
\node (e1) at (3.5,0) [point] {};
\node (e2) at (2.8,0.7) [point] {};
\node (e3) at (4.2,0.7) [point] {};
\node (e4) at (2.8,1.4) [point] {};
\node (e5) at (4.2,1.4) [point] {};
\node (e6) at (4.2,2.1) [point] {};

\node[label=above:{$G_3$}]  at (5.6,-1)  {};
\node (g1) at (5.6,0) [point] {};
\node (g2) at (4.9,0.7) [point] {};
\node (g3) at (6.3,0.7) [point] {};
\node (g4) at (4.9,1.4) [point] {};
\node (g5) at (6.3,1.4) [point] {};
\node (g6) at (6.3,2.1) [point] {};
\node (g7) at (5.6,2.8) [point] {};

\node[label=above:{$G_4$}]  at (7.7,-1)  {};
\node (h1) at (7.7,0) [point] {};
\node (h2) at (7,0.7) [point] {};
\node (h3) at (8.4,0.7) [point] {};
\node (h5) at (8.4,1.4) [point] {};
\node (h6) at (8.4,2.1) [point] {};
\node (h7) at (7.7,2.8) [point] {};

\node[label=above:{$G_5$}]  at (9.8,-1)  {};
\node (p1) at (9.8,0) [point] {};
\node (p2) at (9.1,0.7) [point] {};
\node (p3) at (10.5,0.7) [point] {};
\node (p4) at (9.1,2.1) [point] {};
\node (p5) at (10.5,1.4) [point] {};
\node (p6) at (10.5,2.1) [point] {};

\node[label=above:{$G_6$}]  at (11.9,-1)  {};
\node (q1) at (11.9,0) [point] {};
\node (q2) at (11.3,0.7) [point] {};
\node (q3) at (12.6,0.7) [point] {};
\node (q4) at (11.3,2.1) [point] {};
\node (q5) at (12.6,1.4) [point] {};
\node (q6) at (12.6,2.1) [point] {};
\node (q7) at (11.3,1.4) [point] {};

   \draw (q2) -- (q5) -- (q4) --(q2) -- (q1) -- (q3) -- (q5) (q2)  (q5) -- (q6) ;

  \draw (p5) -- (p4) --(p2) -- (p1) -- (p3) -- (p5) (p2)  (p5) -- (p6) ;
  \draw (e4) -- (e2) -- (e1) -- (e3) -- (e5) (e2) -- (e5) -- (e6) ;
  \draw (g4) -- (g2) -- (g1) -- (g3) -- (g5) -- (g6)-- (g7)  (g2)-- (g5)  (g4) -- (g7) ;
  \draw (h7) --(h2) -- (h1) -- (h3) -- (h5) -- (h6) -- (h7); 
\end{tikzpicture}
\end{tabular}
\]
\caption{The posets $G_1, \dots, G_6$.}
\label{Fig:posets-G}
\end{figure}
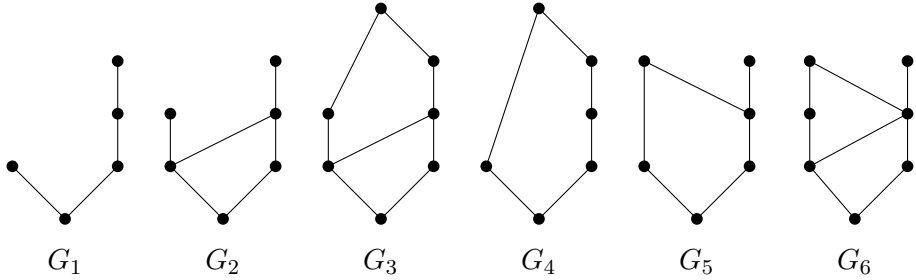

\begin{Lemma}\label{Lem:jankov-axioms-KG-2}
Let $X$ be a rooted Noetherian poset of width $\leq 2$ validating $\beta(P_2)$. Then $X$ validates $\beta(P_3)$ iff it validates the Jankov formulas of the posets in Figure \ref{Fig:posets-G}.
\end{Lemma}

We point out that the posets $K_3$ and $K_4$ are obtained by adding a new top to $K_1$ and $K_2$, respectively. Moreover, $G_3$ and $G_4$ are obtained in a similar manner from $G_2$ and $G_1$. 

In order to shorten the proofs of Lemmas~\ref{Lem:jankov-axioms-KG} and~\ref{Lem:jankov-axioms-KG-2}, we use the following equivalent formulation of Condition (\ref{item:E-partition a}) of Definition~\ref{def: E partition} (see \cite[Rem.~3.1]{BezMor19}): 
\begin{enumerate} 
\item[(1a')] Suppose that $\langle x, y \rangle \in R$ and $x, y \in X$ are distinct. If there is $z \geq x$ such that $y \nleq z$ and $\langle x, z \rangle \notin R$, then there is $u \in X$ such that $y \leq u$ and $\langle z, u \rangle \in R$.
\end{enumerate}

We will also use repeatedly 
that for every poset $X$ and upset $U$, 
identifying $U$ into a point is an E-partition on $X$.

\begin{proof}[Proof of Lemma~\ref{Lem:jankov-axioms-KG}. ]
Let $X$ be a rooted Noetherian poset of width $\leq 2$. Suppose first that $X\vDash \beta(P_2)$. Consider a nonnegative integer $i \leq 7$. By the Fine Lemma, 
to show that  $X \vDash \mathcal{J}(K_i)$, it suffices to prove that $K_i$ is not a p-morphic image of any upset of $X$. Suppose the contrary. Then there exists an upset $U$ of $X$ and a surjective p-morphism $\alpha \colon U \to K_i$. As a consequence, $K_i$ validates all the formulas valid in $X$ and, in particular, $\beta(P_2)$. But in view of Theorem~\ref{Thm:subframe-formulas}(\ref{subframe-formulas 2}) this is false because $P_2$ is isomorphic to a subposet of $K_i$, as it can be checked by inspecting the posets in Figure \ref{Fig:posets-K}.

To prove the converse, assume that $X$ validates the Jankov formulas of the posets $K_1, \dots, K_7$ in Figure~\ref{Fig:posets-K} and suppose, with a view to contradiction, that $X \nvDash \beta(P_2)$. By Theorem \ref{Thm:subframe-formulas}(\ref{subframe-formulas 2}) this implies that $P_2$ is isomorphic to a p-morphic image of a subposet of $X$. The definition of a p-morphism and the structure of $P_2$ imply that actually  $P_2$ is isomorphic to a subposet of $X$. We name the elements of this subposet as in Figure~\ref{Fig:who-is-P2-in-X}.

\begin{figure}
\begin{tabular}{ccccccc}
\begin{tikzpicture}
    \tikzstyle{point} = [shape=circle, thick, draw=black, fill=black , scale=0.35]

\node[label=below:{$\bot$}] (a1) at (2,0) [point] {};
\node[label=left:{$c$}] (a2) at (1,1) [point] {};
\node[label=right:{$d$}] (a3) at (3,1) [point] {};
\node[label=left:{$a$}] (a4) at (1,2) [point] {};
\node[label=right:{$b$}] (a5) at (3,2) [point] {};

  \draw  (a4) -- (a2) -- (a1) -- (a3) -- (a5);

\end{tikzpicture}
\end{tabular}
\caption{The poset $P_2$ viewed as a subposet of $X$. }
\label{Fig:who-is-P2-in-X}
\end{figure}
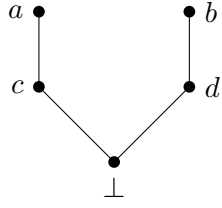

As $X$ is Noetherian and every element in the interval $[c, a]$ is incomprable with $b$ and $d$, we may assume that $c$ is an immediate predecessor of $a$ (otherwise we replace $c$ by a maximal element in $[c, a)$). Similarly, we may assume that $d$ is an immediate predecessor of $b$. By the same token, we may assume that $\bot$ is maximal in ${\downarrow}c \cap {\downarrow}d$, whence we obtain that for every $x \in X$,
\begin{equation}\label{Eq:KG-bot-maximal}
(\text{if }\bot < x \leq c, \text{ then }x \nleq d)\, \, \text{ and }\, \, (\text{if }\bot < x \leq d, \text{ then }x \nleq c).
\end{equation}
Lastly, as the upset of $X$ generated by $\bot$ validates all the formulas valid in $X$, we may also assume that $\bot$ is the minimum of $X$ (otherwise we replace $X$ by ${\uparrow}\bot$). 

Since $({\downarrow}\{a, b\})^c$ is an upset, the following relation is an E-partition of $X$:
\[
R \coloneqq \{ \langle x, y \rangle \in X \times X : x = y \text{ or } x, y \notin {\downarrow}\{ a, b\} \}.
\] 
Notice that the subposet of $X / R$ with the universe $\{ [a], [b], [c], [d], [\bot] \}$ is isomorphic to the subposet of $X$ with the universe $\{ a, b, c, d, \bot \}$. Therefore, since $X / R$ validates all the formulas valid in $X$, we may assume that $R$ is the identity relation (otherwise we replace $X$ by $X / R$). Consequently,
\[
\text{either }({\downarrow}\{ a, b\})^c = \varnothing \, \, \text{ or }\, \, ({\downarrow}\{ a, b\})^c = \{ \top \}
\]
for some element $\top \in X$.

Notice that if $({\downarrow}\{ a, b\})^c = \{ \top \}$, then obviously $\top \nleq a$ and $\top \nleq b$. As $X$ is a  rooted poset of width $\leq 2$ and $a$ and $b$ are incomparable, by symmetry we may assume that $b$ and $\top$ are comparable which, together with $\top \nleq b$, yields $b < \top$. Therefore, one of the following Conditions holds.
\begin{enumerate}
\item $({\downarrow}\{ a, b\})^c = \varnothing$;
\item $a, b < \top$;
\item $a$ and $\top$ are incomparable and $b < \top$.
\end{enumerate}
If $({\downarrow}\{ a, b\})^c = \{ \top \}$, the subposet of $X$ with the universe $\{ \bot, a, b, c, d, \top \}$ is one of the three $X_i$ depicted in Figure~\ref{Fig:who-is-P2}. Thus, 
\begin{equation}\label{Eq:what:I:need-now}
X = {\downarrow}\{ a, b\} \, \,  \text{ or } \, \, (X = \{ \top \} \cup {\downarrow}\{ a, b \} \text{ and the poset }\{ \bot, a, b, c, d, \top \}\text{ is one of the }X_i).
\end{equation}
As $a$ and $b$ are incomparable, this implies that for every $x \in X$,
\begin{equation}\label{Eq:where-top-lies-in-K}
\text{if }a < x \text{ or } b < x\text{, then }x = \top > b.
\end{equation}

\begin{figure}
\begin{tabular}{ccccccc}
\begin{tikzpicture}
    \tikzstyle{point} = [shape=circle, thick, draw=black, fill=black , scale=0.35]

\node[label=above:{$X_1$}]  at (2,-1.5)  {};
\node[label=below:{$\bot$}] (a1) at (2,0) [point] {};
\node[label=left:{$c$}] (a2) at (1,1) [point] {};
\node[label=right:{$d$}] (a3) at (3,1) [point] {};
\node[label=left:{$a$}] (a4) at (1,2) [point] {};
\node[label=right:{$b$}] (a5) at (3,2) [point] {};
\node[label=above:{$\top$}] (a6) at (2,3) [point] {};

\node[label=above:{$X_2$}]  at (5.5,-1.5)  {};
\node[label=below:{$\bot$}] (b1) at (5.5,0) [point] {};
\node[label=left:{$c$}] (b2) at (4.5,1) [point] {};
\node[label=right:{$d$}] (b3) at (6.5,1) [point] {};
\node[label=left:{$a$}] (b4) at (4.5,2) [point] {};
\node[label=right:{$b$}] (b5) at (6.5,2) [point] {};
\node[label=right:{$\top$}] (b6) at (6.5,3) [point] {};

\node[label=above:{$X_3$}]  at (9,-1.5)  {};
\node[label=below:{$\bot$}] (c1) at (9,0) [point] {};
\node[label=left:{$c$}] (c2) at (8,1) [point] {};
\node[label=right:{$d$}] (c3) at (10,1) [point] {};
\node[label=left:{$a$}] (c4) at (8,2) [point] {};
\node[label=right:{$b$}] (c5) at (10,2) [point] {};
\node[label=right:{$\top$}] (c6) at (10,3) [point] {};

  \draw  (a6) -- (a4) -- (a2) -- (a1) -- (a3) -- (a5) -- (a6);
    \draw   (b4) -- (b2) -- (b1) -- (b3) -- (b5) -- (b6);
    \draw   (c4) -- (c2) -- (c1) -- (c3) -- (c5) -- (c6)-- (c2);

\end{tikzpicture}
\end{tabular}
\caption{The posets $X_1, X_2$ and $X_3$. }
\label{Fig:who-is-P2}
\end{figure}
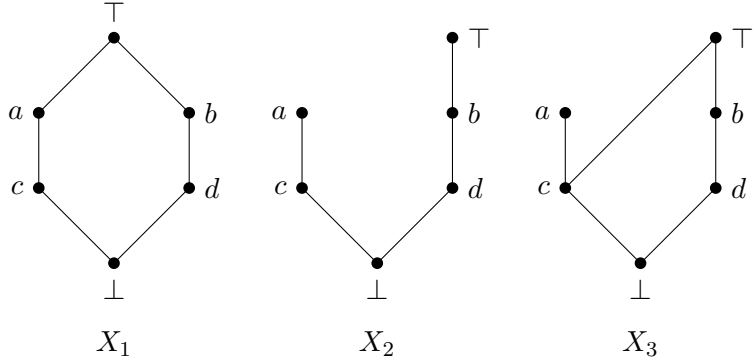

Given a pair $y_1, y_2$ of elements of a poset $Y$, we denote by $(y_1, y_2)$ the open interval $\{ z \in Y : y_1 < z < y_2 \}$. We will prove that
\begin{align}\label{Eq:order-4}
    \begin{split}
X = \left\{ \begin{array}{ll}
 \{ \bot, a, b, c, d\} \cup (\bot, c) \cup (\bot, d) & \text{if $\top$ does not exist}\\
 \{ \bot, a, b, c, d, \top\} \cup (\bot, c) \cup (\bot, d) & \text{if $\top$ exists.}\\
  \end{array} \right.
\end{split}
\end{align} 
The inclusion from right to left is obvious. To prove the other inclusion, consider some $x \in X$ other than $\top$. In view of Condition (\ref{Eq:what:I:need-now}), we have $x \leq a$  or $x \leq b$. If $x \in \{ a, b \}$, we are done. Therefore, we may assume that
\[
\text{either }x < a \, \, \text{ or }\, \, x < b.
\]
Now, if $x \in {\downarrow}c \cup {\downarrow}d$, we are done because ${\downarrow}c \cup {\downarrow}d \subseteq \{ \bot, c, d \} \cup (\bot, c) \cup (\bot, d)$. Consequently, we may assume that $x \notin {\downarrow} \{ c, d \}$. Since $c$ and $d$ are incomparable and $X$ is a rooted poset of width $\leq 2$, this yields 
\[
c < x\, \, \text{ or }\, \, d < x.
\]
As $a$ and $c$ are incomparable with $b$ and $d$, the two displays above imply that
\[
\text{either }c < x < a \, \, \text{ or } \, \, d < x < b.
\]
But this contradicts the assumption that $c$ (resp.\ $d$) is an immediate predecessor of $a$ (resp.\ $b$). Hence,  Condition (\ref{Eq:order-4}) holds as desired.

Now, we consider the sets
\begin{align*}
Y_1 &= (\bot, c] \smallsetminus ({\downarrow}\top \cup {\downarrow} b);\\
Y_2 &= ((\bot, c] \cap {\downarrow}\top) \smallsetminus {\downarrow}b;\\
Y_3 &= (\bot, c] \cap {\downarrow} b; \\
Y_4 &= (\bot, d] \smallsetminus {\downarrow}a;\\
Y_5 &= (\bot, d] \cap {\downarrow}a.
\end{align*}
If $\top$ does not exist, the expression ${\downarrow}\top$ in the above definition should be interpreted as denoting the empty set.

\begin{Claim}\label{Claim : collapsing Yi : correct partition}
The following relation is an E-partition of $X$:
\[
S = \{ \langle x, y \rangle \in X \times X : x = y \text{ or } x, y \in Y_i \text{ for some }i \leq 5 \}.
\]
\end{Claim}

\begin{proof}[Proof of the Claim.]
Since the various $Y_i$ are pairwise disjoint, $S$ is an equivalence relation on $X$. To prove that it is also an E-partition, it suffices to show that there are no distinct $x, y \in X$ such that $\langle x, y \rangle \in S$ and there exists $z \in X$ such that $x \leq z$ and $y \nleq z$ and $\langle x, z \rangle \notin S$ (see Condition (1a') if necessary). Suppose the contrary, with a view to contradiction. Since $x$ and $y$ are distinct and related by $S$, we have
\[
\bot < x, y \leq c \, \, \text{ or } \, \, \bot < x, y \leq d.
\]

Suppose first that $\bot < x, y \leq c$. As $d \nleq c$, we have $d \nleq x, y$. Furthermore, from Condition~(\ref{Eq:KG-bot-maximal})  it follows that $x, y \nleq d$. Thus, $x$ and $y$ are incomparable with $d$. As $X$ is a rooted poset of width $\leq 2$, this implies that $x$ and $y$ are comparable. Since by assumption $x \leq z$ and $y \nleq z$, we conclude that $x < y$. Now, by applying the assumption that $X$ has width $\leq 2$ to the fact that $y$ and $d$ are incomparable, we obtain that $z$ is comparable with either $y$ or $d$. We will prove that $z$ is incomparable with $y$. On the one hand, by assumption $y \nleq z$. On the other hand, if $z \leq y$, then we would have $x \leq z \leq y$, because $x \leq z$ by assumption. Since $\langle x, y \rangle \in S$, the equivalence class $[x]$ contains the interval $(x, y)$. In particular, $\langle x, z \rangle \in S$, a contradiction. Thus, we conclude that $z \nleq y$. Consequently, $y$ and $z$ are incomparable, which in turn means that $d$ and $z$ are comparable. Since $x$ and $d$ and incomparable and $x \leq z$, this means that $d < z$. A similar argument shows that if $\bot < x, y \leq d$, then $c < z$. Thus, we obtain that 
\[
(\bot < x, y \leq c \text{ and } d < z)\, \, \text{ or } \, \, (\bot < x, y \leq d \text{ and }c < z).
\]

We need to prove that both cases lead to a contradiction. First suppose that $\bot < x, y \leq c$ and $d < z$. Since $d < z$, by Condition~(\ref{Eq:order-4}) we obtain that $z \in \{ b, \top \}$. Recall that $x, y \in Y_i$ for some $i \leq 5$ because $x$ and $y$ are different and related by $S$. Furthermore, as $\bot < x \leq c$, Condition~(\ref{Eq:KG-bot-maximal})  implies that $Y_i$ is $Y_1$, $Y_2$, or $Y_3$. We have two cases: either $z = b$ or $z = \top$. If $z = \top$, then $Y_i \ne Y_1$ because $x \in Y_i$ and $x \leq z = \top$ and $Y_1 \subseteq ({\downarrow}\top)^c$. Therefore, $Y_i$ is $Y_2$ or $Y_3$. Since $Y_2 \cup Y_3 \subseteq {\downarrow}\top = {\downarrow}z$ and $y \in Y_i$, we obtain $y \leq \top = z$, a contradiction. 
If $z = b$, then $Y_i = Y_3$ because $x \in Y_i$, $x \leq z = b$, and $Y_1 \cup Y_2 \subseteq ({\downarrow}b)^c$. As a consequence, $y \in Y_i = Y_3$. Since $Y_3 \subseteq {\downarrow}b = {\downarrow}z$, this implies that $y \leq z$, a contradiction.

Next we consider the case where $\bot < x, y \leq d$, and $c < z$. Since $c < z$, Condition (\ref{Eq:order-4}) implies that $z = \top$ or $z = a$. If $z = \top$, then $z \geq d \geq y$, a contradiction. Suppose that $z = a$. Recall that $x, y \in Y_i$ for some $i \leq 5$ because $x$ and $y$ are distinct and related by $S$. As $\bot < x\leq d$, Condition~(\ref{Eq:KG-bot-maximal}) implies that $Y_i$ is $Y_4$ or $Y_5$. But, as $x \in Y_i$, $x \leq z = a$, and $Y_4 \subseteq ({\downarrow}a)^c$, we must have $Y_i = Y_5$. Consequently, $y \in Y_i = Y_5 \subseteq {\downarrow}a$. This yields that $y \leq a = z$, a contradiction. Hence, we conclude that $S$ is an E-partition of $X$.
\end{proof}

\begin{Claim}\label{Eq:does-not-alter-T2}
For every $x, y \in \{ \bot, a, b, c, d, \top \}$,
\[
x \leq y \Longleftrightarrow [x] \leq [y]. 
\]
\end{Claim}

\begin{proof}[Proof of the Claim.]
The implication from left to right is obvious.\ To prove the other implication, suppose that $[x] \leq [y]$. The definition of $S$ guarantees that
\[
[\bot] = \{ \bot \} \qquad [\top] = \{ \top \} \qquad [a] = \{ a \} \qquad [b] = \{ b \} \qquad [d] = (\bot, d] \smallsetminus {\downarrow}a
\]
and
\[
[c] = \left\{ \begin{array}{ll}
 (\bot, c] \smallsetminus {\downarrow}b & \text{if $\top$ does not exist or it exists and $c \leq \top$}\\
 (\bot, c] \smallsetminus {\downarrow}\top & \text{if $\top$ exists and $c \nleq \top$.}
  \end{array} \right. 
 \]
Consequently, $y$ is the maximum of the equivalence class $[y]$. Therefore, the assumption that $[x] \leq [y]$ guarantees the existence of some $x' \in [x]$ such that $x' \leq y$. If the equivalence class $[x]$ is a singleton, then $x \leq y$ and we are done. 
Otherwise, in view of the above displays, $x$ is either $c$ or $d$. As $\bot$ is maximal in ${\downarrow}c \cap {\downarrow}d$, the above displays guarantee that $[c]$ and $[d]$ are incomparable. Then we may assume that $x \in \{ c, d \}$ and  $y \in \{ \bot, \top, a, b \}$. 

We begin by the case where $x = d$. If $y \in \{ b, \top \}$, then clearly $x \leq y$. Therefore, we consider the case where $y \in \{ \bot, a \}$. As the set $[x] = [d] = (\bot, d] \smallsetminus {\downarrow}a$ does not contain any element below $\bot$ or $a$, this case never happens and we are done.

Then we turn our attention to the case where $x = c$. If $y = a$, then clearly $x  \leq  y$. Moreover, $[x] = [c] \subseteq (\bot, c]$, and hence $[x]$ does not contain any element below $\bot$ or $b$. Together with $[\bot] = \{ \bot \}$ and $[b] = \{ b \}$, this yields that $y \notin \{ b, \bot \}$. It only remains to consider the case where $y = \top$. But the above display guarantees that if $[c]$ contains an element below $\top$, then $c \leq \top$, and hence $x \leq y$ as desired. 
\end{proof}

Together with the fact that $X / S$ validates all the formulas valid in $X$, Claim \ref{Eq:does-not-alter-T2} allows us to assume that $S$ is the identity relation (otherwise we replace $X$ by $X / S$). 

\begin{Claim}\label{Claim:Conditions:M:part:1}
One of the following conditions holds.
\begin{enumerate}
\renewcommand{\labelenumi}{(\theenumi)}
\renewcommand{\theenumi}{\textup{M\arabic{enumi}}}
\item\label{C:M1} $(\bot, c) = \varnothing$;
\item\label{C:M2} $(\bot, c) = \{ x \}$ for some $x$ such that ${\uparrow}x = \{x, a, c\} \cup {\uparrow}b$;
\item\label{C:M3} $\top$ exists, $c \nleq \top$, and $(\bot, c) = \{ x \}$ for some $x$ such that ${\uparrow}x = \{x, a, c, \top\}$;
\item\label{C:M4} $\top$ exists, $c \nleq \top$, and $(\bot, c) = \{ x, y \}$ for $x$ and $y$ such that ${\uparrow} y = \{ y, a, c, \top \}$ and ${\uparrow}x = \{x, y, a, b, c, \top\}$.
\end{enumerate}
\end{Claim}

\begin{proof}[Proof of the Claim.]
First, if $(\bot, c)$ is empty, Condition (\ref{C:M1}) holds. 
Suppose $ (\bot, c) \ne \varnothing$. The definition of the various $Y_i$ guarantees that $(\bot, c] \subseteq Y_1 \cup Y_2 \cup Y_3$. Bearing in mind that each $Y_i$ is either empty or a singleton (since $S$ is the identity relation), this implies that $(\bot, c]$ has at most three elements, which in turn means that $(\bot, c)$ has at most two. Furthermore,
\begin{equation}\label{Eq:counting-the-Yi}
\{ \{ x \} : x \in (\bot, c] \} \subseteq \{ Y_1, Y_2, Y_3 \}.
\end{equation} 

First suppose that $(\bot, c)$ has precisely two elements $x$ and $y$. By Condition~(\ref{Eq:KG-bot-maximal}), $x$ and $y$ are incomparable with $d$. As $X$ is a rooted poset of width $\leq 2$, this yields that both $x$ and $y$ must be comparable. Without loss of generality we may assume that $\bot < x < y < c$. From Condition~(\ref{Eq:counting-the-Yi}) it follows that
\[
\{ \{ x\}, \{ y\},\{ c \} \} = \{ Y_1, Y_2, Y_3 \}.
\]
Together with $x < y < c$ and the definition of the various $Y_i$, this implies that
\[
Y_1 = \{ c \} \qquad Y_2 = \{ y \} \qquad Y_3 = \{ x \}.
\]
This, in turn, guarantees that $\top$ exists and that
\[
x \leq b \qquad y \leq \top \qquad y \nleq b \qquad c \nleq \top.
\]
Together with Condition (\ref{Eq:order-4}) and the facts that $(\bot, c ) = \{ x, y \}$, $\bot < x < y <  c$, and  $x, y,  c \nleq d$, this implies that 
\[
{\uparrow}x = \{ x, y, a, b, c, \top \} \qquad {\uparrow}y = \{ y, a, c, \top \}.
\]
Therefore, Condition~(\ref{C:M4}) holds. 

It only remains to consider the case where 
$(\bot, c) = \{ x \}$ for some $x \in X$. First suppose that $x \leq b$. As before, Condition~(\ref{Eq:KG-bot-maximal}) implies that $x \nleq d$. Together with Condition~(\ref{Eq:order-4}) and the assumption that $(\bot, c) = \{ x \}$, this implies that ${\uparrow}x = \{x, a, c \} \cup {\uparrow}b$. Hence, Condition~(\ref{C:M2}) holds. 
Next suppose that $x \nleq b$. Since $x, c \nleq b$, both $[c]$ and $[x]$ are different from $Y_3$. By Condition~(\ref{Eq:counting-the-Yi}) and the fact that $c \ne x$, this implies that $\{\{c\}, \{x\}\} = \{ Y_1, Y_2 \}$. As $x \leq c$, the definition of the various $Y_i$ guarantees that $c \in Y_1$ and $x \in Y_2$. Consequently, $\top$ exists and 
\[
c \nleq \top \qquad x \leq \top \qquad x \nleq b.
\]
Bearing in mind that $x \nleq d$ by Condition (\ref{Eq:KG-bot-maximal}), we conclude that ${\uparrow}x = \{ x, a, c, \top \}$. Therefore, Condition (\ref{C:M3}) holds as desired.
\end{proof}

A similar (but shorter) argument 
yields the following:\footnote{The statement of Claim \ref{Claim:Conditions:N:part:1} is simpler than that of Claim \ref{Claim:Conditions:M:part:1} because of the asymmetric behavior of the elements $d$ and $c$ (see Figure \ref{Fig:who-is-P2}).  
In particular, the definition of the sets $Y_1, \dots, Y_5$ and the fact that the relation $S$ in Claim \ref{Claim : collapsing Yi : correct partition} is assumed to be the identity relation ensure that the interval $(\bot, d]$ has at most two elements, while $(\bot, c]$ may have three.}

\begin{Claim}\label{Claim:Conditions:N:part:1}
One of the following conditions holds.
\begin{enumerate}
\renewcommand{\labelenumi}{(\theenumi)}
\renewcommand{\theenumi}{\textup{N\arabic{enumi}}}
\item\label{C:N1} $(\bot, d) = \varnothing$;
\item\label{C:N2} $(\bot, d) = \{ x \}$ for some $x$ such that ${\uparrow}x = \{x, a \} \cup {\uparrow}d$.
\end{enumerate}
\end{Claim}

At last, we are ready to give a more concrete description of the poset $X$. First, the order structure of the subposet $\{ \bot, a, b, c, d\}$ of $X$ is that of Figure~\ref{Fig:who-is-P2-in-X} 
 and, if $\top$ exists, the subposet $\{ \bot, a, b, c, d, \top \}$ is one of those depicted in Figure \ref{Fig:who-is-P2}. By Condition~(\ref{Eq:order-4}), the elements of $X$ other than $a, b, c, d, \top$, and  $\bot$ lie in $(\bot, c) \cup (\bot, d)$. But recall from  Claim \ref{Claim:Conditions:M:part:1} that one of Conditions~(\ref{C:M1})--(\ref{C:M4}) holds and that each of them gives a complete description of the interval $(\bot, c)$. Similarly, one of (\ref{C:N1}) or (\ref{C:N2}) holds by Claim \ref{Claim:Conditions:N:part:1} and each of them gives a complete description of the interval $(\bot, d)$.

As a consequence, we obtain that $X$ is a subposet of one of the rooted posets in Figure~\ref{Fig:the-K-poset-final-part}. Moreover, $X$ contains $\bot, a, b, c, d$ plus $\top$ if $\top$ appears in the corresponding picture.
\begin{figure}
\begin{tabular}{ccccccc}
\begin{tikzpicture}
    \tikzstyle{point} = [shape=circle, thick, draw=black, fill=black , scale=0.35]
\node[label=above:{$Z_1$}]  at (-1.5,-1.5)  {};
\node[label=below:{$\bot$}] (d1) at (-1.5,0) [point] {};
\node[label=left:{$x$}] (d2) at (-2.5,1) [point] {};
\node[label=right:{$y$}] (d3) at (-0.5,1) [point] {};
\node[label=left:{$c$}] (d4) at (-2.5,2) [point] {};
\node[label=right:{$d$}] (d5) at (-0.5,2) [point] {};
\node[label=left:{$a$}] (d6) at (-2.5,3) [point] {};
\node[label=right:{$b$}] (d7) at (-0.5,3) [point] {};

\node[label=above:{$Z_2$}]  at (2,-1.5)  {};
\node[label=below:{$\bot$}] (a1) at (2,0) [point] {};
\node[label=left:{$x$}] (a2) at (1,1) [point] {};
\node[label=right:{$y$}] (a3) at (3,1) [point] {};
\node[label=left:{$c$}] (a4) at (1,2) [point] {};
\node[label=right:{$d$}] (a5) at (3,2) [point] {};
\node[label=left:{$a$}] (a6) at (1,3) [point] {};
\node[label=right:{$b$}] (a7) at (3,3) [point] {};
\node[label=above:{$\top$}] (a8) at (2,4) [point] {};

\node[label=above:{$Z_3$}]  at (5.5,-1.5)  {};
\node[label=below:{$\bot$}] (b1) at (5.5,0) [point] {};
\node[label=left:{$y$}] (b2) at (4.5,1) [point] {};
\node[label=right:{$z$}] (b3) at (6.5,1) [point] {};
\node[label=left:{$x$}] (b4) at (4.5,2) [point] {};
\node[label=right:{$d$}] (b5) at (6.5,2) [point] {};
\node[label=left:{$a$}] (b6) at (4.5,4) [point] {};
\node[label=left:{$c$}] (b9) at (4.5,3) [point] {};
\node[label=right:{$b$}] (b7) at (6.5,3) [point] {};
\node[label=right:{$\top$}] (b8) at (6.5,4) [point] {};

\node[label=above:{$Z_4$}]  at (9,-1.5)  {};
\node[label=below:{$\bot$}] (c1) at (9,0) [point] {};
\node[label=left:{$x$}] (c2) at (8,1) [point] {};
\node[label=right:{$y$}] (c3) at (10,1) [point] {};
\node[label=left:{$c$}] (c4) at (8,2) [point] {};
\node[label=right:{$d$}] (c5) at (10,2) [point] {};
\node[label=right:{$b$}] (c6) at (10,3) [point] {};
\node[label=left:{$a$}] (c7) at (8,3) [point] {};
\node[label=right:{$\top$}] (c8) at (10,4) [point] {};

  \draw  (a6) -- (a4) -- (a2) -- (a1) -- (a3) -- (a5) -- (a7) -- (a8) -- (a6)  (a2) -- (a7) (a3) -- (a6);
    \draw  (b6) -- (b4) -- (b2) -- (b1) -- (b3) -- (b5) -- (b8) (b2) -- (b7) (b4) -- (b8) (b3) -- (b6);
    \draw (c6) -- (c8) -- (c4) (c3) -- (c7)  -- (c4) -- (c2) -- (c1) -- (c3) -- (c5) -- (c6)-- (c2);
  \draw  (d6) -- (d4) -- (d2) -- (d1) -- (d3) -- (d5) -- (d7)  (d2) -- (d7) (d3) -- (d6);
       
\end{tikzpicture}
\end{tabular}
\caption{The posets $Z_1, Z_2, Z_3$, and $Z_4$.}
\label{Fig:the-K-poset-final-part}
\end{figure}
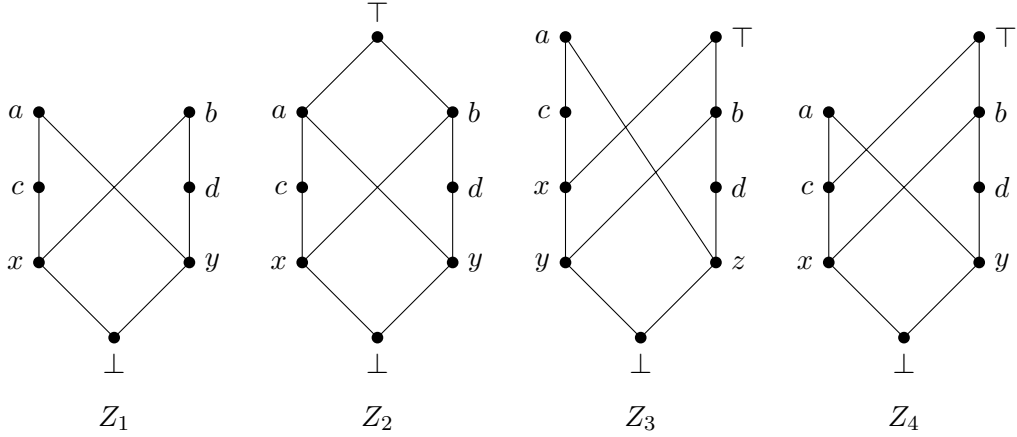

Now, recall that $X$ validates the Jankov formulas of $K_1, \dots, K_7$. Therefore, none of these posets is a p-morphic image of an upset of $X$. Bearing this in mind, we begin by considering the case where $X$ is a subposet of $Z_1$ containing $\bot, a, b, c$, and $d$. Notice that $X \ne \{ \bot, a, b, c, d \}$, otherwise we would obtain $X \cong K_1$, a contradiction. Then $X$ contains $x$ or $y$. If $X$ does not contain both, it is isomorphic to $K_2$ against the assumptions. Therefore, we conclude that $X$ contains both $x$ and $y$. But this is also impossible as in this case $K_3$ is a p-morphic image of $X$.

The case where $X$ is a subposet of $Z_2$ containing $\bot, a, b, c, d$, and $\top$ leads to a contradiction in a similar way (where $K_3$ takes the role of $K_1$ and $K_4$ that of $K_2$).

Next we consider the case where $X$ is a subposet of $Z_3$ containing $\bot, a, b, c, d$, and $\top$. Notice that $X$ contains one of $x, y$, and $z$ (otherwise $K_1$ is a p-morphic of $X$, which is impossible). If $x$ and $z$ or $y$ and $z$ belong to $X$, then $K_3$ is a p-morphic image of $X$, which is also false. More precisely, 
when $x, z \in X$, we collapse $\{ a, b, \top \}$ if $y \notin X$ and we collapse  $\{ a, b, c, \top\}$ if $y \in X$. In both cases, we obtain a p-morphic image of $X$ isomorphic to $K_3$. On the other hand, 
when $y, z \in X$ but $x \notin X$, we collapse $\{ a, b, \top \}$, thus obtaining a p-morphic image of $X$ isomorphic to $K_3$. Therefore, we may assume that the universe of $X$ is the union of $A \coloneqq \{ \bot, a, b, c, d, \top \}$ with $\{ x\}$ or $\{ y \}$ or $\{ x, y \}$ or $\{ z \}$. We will show that each of these cases leads to a contradiction. 

If $X = A \cup \{ x \}$ or $X = A \cup \{ y \}$, then $K_2$ is a p-morphic image of $X$ obtained by collapsing $\{ b, d, \top \}$. Moreover, if $X = A \cup \{ x, y \}$, then $K_6$ is a p-morphic image of $X$ obtained by collapsing $\{ a, c \}$. Lastly, if $X = A \cup \{ z \}$, then $K_2$ is a p-morphic image of $X$ obtained by collapsing $\{ b, \top \}$.

It remains to consider the case where $X$ is a subposet of $Z_4$ containing $\bot, a, b, c, d$, and $\top$. Observe that $X \ne \{ \bot, a, b, c, d, \top \}$ (otherwise $X \cong K_5$, which is false). Therefore, $x$ or $y$ belong to $X$. If both $x$ and $y$ belong to $X$, then $K_3$ is a p-morphic image of $X$ obtained by collapsing $\{ a, b, \top \}$, against the assumptions. Thus, $X$ is 
$A$ together with $x$ or $y$. If $X = A \cup \{ x \}$, then $X \cong K_6$, which is false. On the other hand, if $X = A \cup \{ y \}$, then $X \cong K_7$, which is also false. Hence, we reach the desired contradiction. 
\end{proof}

\begin{proof}[Proof of Lemma~\ref{Lem:jankov-axioms-KG-2}.]
Let $X$ be a rooted Noetherian poset of width $\leq 2$ validating $\beta(P_2)$. By Theorem~\ref{Thm:subframe-formulas}(\ref{subframe-formulas 2}), from $X \vDash \beta(P_2)$ it follows that $P_2$ is not a p-morphic image of any subposet of $X$. This fact will be used repeatedly in the proof.

First suppose that $X\vDash \beta(P_3)$ and consider a nonnegative integer $i \leq 6$. By the Fine Lemma, 
to show that  $X$ validates $\mathcal{J}(G_i)$, it suffices to prove that $G_i$ is not a p-morphic image of any upset of $X$. Suppose the contrary. Then there exist an upset $U$ of $X$ and a surjective p-morphism $\alpha \colon U \to G_i$. As a consequence, $G_i$ validates all the formulas valid in $X$ and, in particular, $\beta(P_3)$. But in view of Theorem~\ref{Thm:subframe-formulas}(\ref{subframe-formulas 2}) this is false because $P_3$ is isomorphic to a subposet of $G_i$, as it can be checked by inspecting the posets in Figure~\ref{Fig:posets-G}.

To prove the converse, assume that $X$ validates the Jankov formulas of the posets $G_1, \dots, G_6$ in Figure~\ref{Fig:posets-G} and suppose, with a view to contradiction, that $X\nvDash \beta(P_3)$. By Theorem~\ref{Thm:subframe-formulas}(\ref{subframe-formulas 2}) this implies that $P_3$ is isomorphic to a p-morphic image of a subposet of $X$. The definition of a p-morphism and the structure of $P_3$ imply that actually $P_3$ is isomorphic to a subposet of $X$. We name the elements of this subposet as in Figure~\ref{Fig:posets-G-in-X}.
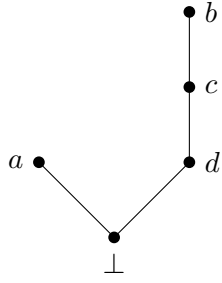
\begin{figure}
\begin{tabular}{ccccccc}
\begin{tikzpicture}
    \tikzstyle{point} = [shape=circle, thick, draw=black, fill=black , scale=0.35]
    
\node[label=below:{$\bot$}] (a1) at (2,0) [point] {};
\node[label=left:{$a$}] (a2) at (1,1) [point] {};
\node[label=right:{$d$}] (a3) at (3,1) [point] {};
\node[label=right:{$c$}] (a5) at (3,2) [point] {};
\node[label=right:{$b$}] (a6) at (3,3) [point] {};

  \draw (a2) -- (a1) -- (a3) -- (a5)-- (a6) ;
         
\end{tikzpicture}
\end{tabular}
\caption{The poset $P_3$ viewed as a subposet of $X$.}
\label{Fig:posets-G-in-X}
\end{figure}

We may assume that $\bot$ is the minimum of $X$ (otherwise we replace $X$ by ${\uparrow}\bot$). In addition, since $X$ is Noetherian, we may assume that 
\begin{equation}\label{Eq:G-part-first-assumption}
\text{$\bot$ is maximal in ${\downarrow}a \cap {\downarrow}d$ and $b$ (resp.\ $c$) is an immediate successor of $c$ (resp.\ $d$).}
\end{equation}

\begin{Claim}\label{Claim:wlog:X:has:some:structure}
We may assume, without loss of generality, that either $X = {\downarrow} \{ a, b \}$ or $X = \{ \top \} \cup {\downarrow} \{ a, b \}$ for some $\top \in X$ such that the subposet of $X$ with the universe $\{ \bot, a, b, c, d, \top \}$ is one of the posets $X_i$ depicted in Figure \ref{Fig:who-is-P3}. 
\end{Claim}

\begin{figure}[t]
\begin{tabular}{ccccccc}
\begin{tikzpicture}
    \tikzstyle{point} = [shape=circle, thick, draw=black, fill=black , scale=0.35]
\node[label=below:{$X_1$}] at (2,-0.5) {};

\node[label=below:{$\bot$}] (a1) at (2,0) [point] {};
\node[label=left:{$a$}] (a2) at (1,1) [point] {};
\node[label=right:{$d$}] (a3) at (3,1) [point] {};
\node[label=right:{$c$}] (a5) at (3,2) [point] {};
\node[label=right:{$b$}] (a6) at (3,3) [point] {};
\node[label=above:{$\top$}] (a7) at (2,4) [point] {};

\node[label=below:{$X_2$}] at (5.5,-0.5) {};

\node[label=below:{$\bot$}] (b1) at (5.5,0) [point] {};
\node[label=left:{$a$}] (b2) at (4.5,1) [point] {};
\node[label=right:{$d$}] (b3) at (6.5,1) [point] {};
\node[label=right:{$c$}] (b5) at (6.5,2) [point] {};
\node[label=right:{$b$}] (b6) at (6.5,3) [point] {};
\node[label=above:{$\top$}] (b7) at (4.5,3) [point] {};

  \draw (a2) -- (a1) -- (a3) -- (a5)-- (a6) -- (a7) -- (a2) (b5) -- (b7) -- (b2) ;
  
      \draw (b2) -- (b1) -- (b3) -- (b5)-- (b6)  ;
   
\end{tikzpicture}
\end{tabular}
\caption{The posets $X_1$ and $X_2$.}
\label{Fig:who-is-P3}
\end{figure}
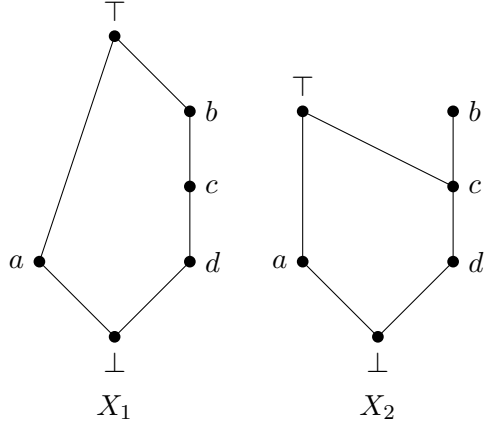

\begin{proof}[Proof of the Claim.]
If $X = {\downarrow} \{ a, b \}$, we are done. 
Suppose $X \ne {\downarrow} \{ a, b \}$. Since $({\downarrow}\{a, b \})^c$ is an upset, the following relation is an E-partition of $X$:
\[
R = \{ \langle x, y \rangle \in X \times X : x = y \text{ or } x, y \notin {\downarrow}\{ a, b\} \}.
\]
Notice that the subposet of $X / R$ with the universe $\{ [a], [b], [c], [d], [\bot] \}$ is isomorphic to the subposet of $X$ with the universe $\{ a, b, c, d, \bot \}$. Therefore, since $X / R$ validates all the formulas valid in $X$, we may assume that $R$ is the identity relation (otherwise we replace $X$ by $X / R$). Consequently, being nonempty by assumption, the set $({\downarrow} \{ a, b \})^c$ is an equivalence class of the identity relation $R$. Therefore, $({\downarrow} \{ a, b \})^c = \{ \top \}$ for some $\top \in X$. This implies that $X = \{ \top \} \cup {\downarrow} \{ a, b \}$.

It only remains to prove that the subposet of $X$ with the universe $\{ \bot, a, b, c, d, \top \}$ is one of the posets $X_i$. First, recall that $\top \nleq a, b$. Since $X$ is a rooted poset of width $\leq 2$, this yields that either $a < \top$ or $b < \top$. If $a, b < \top$, then the subposet $\{ \bot, a, b, c, d, \top \}$ is isomorphic to $X_1$ and we are done. 
Suppose $a \nleq \top$ or $b \nleq \top$. Since $\top \nleq a, b$, we have that
\[
(a \text{ and }\top \text{ are incomparable and }b < \top) \text{ or }(b \text{ and }\top \text{ are incomparable and }a < \top).
\]

First suppose that $a$ and $\top$ are incomparable and $b < \top$. Since ${\uparrow}b$ is an upset, the following relation is an E-partition of $X$:
\[
S = \{ \langle x, y \rangle \in X \times X : x = y \text{ or } b \leq x, y \}.
\]
Notice that $S$ does not alter the order relation between $\bot, a, b, c, d$. Together with the fact that $X / S$ validates all the formulas valid in $X$, this means that $X / S$ is a poset of width $\leq 2$ that validates $\beta(P_2)$ and the various $\mathcal J(G_i)$. Moreover, $X/S$ contains a subposet isomorphic to $P_3$, namely $\{ [\bot], [a], [b], [c], [d] \}$.  Therefore, in our proof we may replace $X$ by $X/ S$ and each element $\bot, a, b,c,d$ by its equivalence class. Furthermore, since $b < \top$, the definition of $S$ ensures that $[\top] = [b]$. Bearing in mind that $X = \{ \top \} \cup {\downarrow} \{ a, b \}$, this means that $X / S = {\downarrow} \{ [a], [b] \}$. Because of this, by replacing $X$ by $X / S$, we may assume that $X = {\downarrow} \{ a, b \}$ as desired.

Therefore, it only remains to consider the case where $b$ and $\top$ are incomparable and $a < \top$. We will prove that $c \leq \top$. Suppose the contrary. 
Then $c$ and $\top$ are incomparable because by assumption $c \leq b$ and $\top \nleq b$. Therefore, both $c$ and $b$ are incomparable with $\top$. By assumption, they are also incomparable with $a$. Together with the fact that $\bot < a < \top$ and $\bot<  c < b$, this implies that $\{ \bot, a, b, c, \top \}$ is a subposet of $X$ isomorphic to $P_2$. But this contradicts the assumption that $X$ validates $\beta(P_2)$. Hence, we conclude that $c \leq \top$ as desired. Bearing in mind that $b$ and $\top$ are incomparable, that $a < \top$, and that the structure of the poset $\{ \bot, a, b, c, d \}$ is as in Figure \ref{Fig:posets-G-in-X}, this implies that the subposet of $X$ with universe $\{ \bot, a, b, c, d, \top \}$ is $X_2$.
\end{proof}

\begin{Claim}\label{Eq:order-G4}
We have that 
\begin{align*}
    \begin{split}
X = \left\{ \begin{array}{ll}
 \{ \bot, a, b, c, d\} \cup ({\downarrow}c \cap (\bot, a)) \cup (\bot, d) & \text{if $\top$ does not exist}\\
 \{ \bot, a, b, c, d, \top\} \cup ({\downarrow}c \cap (\bot, a)) \cup (\bot, d) & \text{if $\top$ exists.}\\
  \end{array} \right.
\end{split}
\end{align*} 
\end{Claim}

\begin{proof}[Proof of the Claim.]
The inclusion from right to left is obvious. To prove the other inclusion, consider some $x \in X \smallsetminus \{ \bot, a, b, c, d \}$ other than $\top$. In view of the Claim \ref{Claim:wlog:X:has:some:structure}, either $X = {\downarrow} \{ a, b \}$ or $X  = \{ \top \} \cup {\downarrow} \{ a, b \}$. Since $x$ is different from $a, b, \top$ and from the minimum $\bot$, we have two cases: $x \in (\bot, a)$ or $x \in (\bot, b)$.

First suppose that $x \in (\bot, a)$. To prove that $x$ belongs to the set in the right hand side of the statement, it suffices to show that $x \leq c$. Suppose the contrary, with a view to contradiction. Since $x \leq a$ and $c \nleq a$, this means that $x$ and $c$ are incomparable. Furthermore, since $\bot$ is maximal in ${\downarrow}a \cap {\downarrow}d$ by Condition (\ref{Eq:G-part-first-assumption}) and $\bot < x \leq a$, we obtain $x \nleq d$. In addition, $d \nleq x$ because $x \leq a$ and $d \nleq a$. Thus, $x$ is also incomparable with $d$. Therefore, $\bot < x < a$ and $\bot < d < c$ and $x, a$ are incomparable with $d, c$. Consequently, $\{ \bot, a, b, c, d \}$ is a subposet of $X$ isomorphic to $P_2$. But this contradicts the assumption that $X\vDash \beta(P_2)$.

Next suppose that $x \in (\bot, b)$. We may assume that $x \nless a$ (otherwise $x \in (\bot, a)$ and we repeat the argument of the previous case). Consequently, in order to prove that $a$ and $x$ are incomparable, it suffices to show that $a \nleq x$. But this is clear because by assumption $x < b$ and $a \nleq b$. Furthermore, by assumption, $a$ is incomparable with $c$ and $d$. Together with the facts that $a$ is incomparable with $x$ and 
that $X$ has width $\leq 2$, this implies that $x$ is comparable with both $c$ and $d$. If $x < d$, then $x \in (\bot, d)$, and hence $x$ belongs to the right hand side of the statement. 
Suppose $x \nless d$. Since $x$ and $d$ are comparable and distinct, this means that $d < x$. Together with the assumption in Condition~(\ref{Eq:G-part-first-assumption}) that $c$ is an immediate successor of $d$ and the fact that $x \ne c$, this implies that $x \nleq c$. Since $x$ and $c$ are comparable, we obtain that $c < x$. But then we have $c < x < b$, a contradiction to the assumption that $b$ is an immediate successor of $c$ (see Condition~(\ref{Eq:G-part-first-assumption})). 
\end{proof}

Now, we consider the relation
\[
T = \{ \langle x, y \rangle \in X \times X : x = y \text{ or } x, y \in (\bot, a) \cap {\downarrow}c \text{ or } x, y \in (\bot, d] \},
\]
where $(\bot, d]$ stands for $\{ x \in X : \bot < x \leq d \}$.

\begin{Claim}
The relation $T$ is an E-partition of $X$.
\end{Claim}

\begin{proof}[Proof of the Claim.]
Since $\bot$ maximal in ${\downarrow}a \cap {\downarrow} d$ by Condition~(\ref{Eq:G-part-first-assumption}), the sets $(\bot, a) \cap {\downarrow}c$ and $(\bot, d]$ are disjoint, and hence $T$ is an equivalence relation on $X$. We will prove that it is also an E-partition.

To this end, it suffices to show that there are no distinct $x, y \in X$ such that $\langle x, y \rangle \in T$ and for which there exists an element $z \in X$ such that $x \leq z$ and $y \nleq z$ and $\langle x, z \rangle \notin T$. Suppose the contrary, with a view to contradiction. Since $x$ and $y$ are distinct and related by $T$, we have that
\[
x, y \in (\bot, a) \cap {\downarrow}c \text{ or } x, y  \in (\bot, d].
\]

First suppose that $x, y \in (\bot, a) \cap {\downarrow}c$. From Claim \ref{Eq:order-G4} it follows that $z \in \{ \bot, a, b, c, d, \top \} \cup ((\bot, a) \cap {\downarrow}c) \cup (\bot, d)$. Clearly, $z \notin (\bot, a) \cap {\downarrow}c$ (otherwise $\langle x, z \rangle \notin T$ contradicting the assumption). Moreover, $z \notin (\bot, d]$ because otherwise $\bot < x \leq z \leq d$ and by assumption $\bot < x \leq a$, contradicting the maximality of $\bot$ in ${\downarrow}a \cap {\downarrow}d$ (see Condition~(\ref{Eq:G-part-first-assumption})). Therefore, $z \in \{ \bot, a, b, c, \top \}$. Since by assumption $y \nleq z$ and $y \in (\bot, a) \cap {\downarrow}c$, we obtain that $z \notin \{ a, c, b, \top \}$. Consequently, $z = \bot$. But this contradicts the assumption that $\bot < x \leq z$. 

Next suppose that $x, y  \in (\bot, d]$. Since $y \nleq z$, this implies that $d \nleq z$. Furthermore, $z \nleq d$ (otherwise $\bot < x \leq z \leq d$, and hence $\langle x, z \rangle \in T$, a contradiction). Therefore, $z$ and $d$ are incomparable. Since $a$ and $d$ are also incomparable and $X$ is a rooted poset of width $\leq 2$, we conclude that $z$ and $a$ are comparable. As $\bot$ is maximal in ${\downarrow}a \cap {\downarrow}d$ by Condition~(\ref{Eq:G-part-first-assumption}) and $\bot < x \leq d$, we obtain that $x \nleq a$. Together with $x \leq z$, this implies that $z \nleq a$. Thus, since $z$ and $a$ are comparable, we must have $a < z$. As $\{ \bot, a, b, c, d, \top \}$ is one of the posets depicted in Figure~\ref{Fig:who-is-P3}, we conclude that $z = \top$. But since $y \leq d \leq \top$, this implies that $y \leq z$, a contradiction. 
\end{proof}

Lastly, we will make use of the following.

\begin{Claim}\label{Eq:does-not-alter-2nd-case}
For every $x, y \in \{ \bot, a, b, c, d, \top \}$,
\[
x \leq y \Longleftrightarrow [x] \leq [y].
\]
\end{Claim}

\begin{proof}[Proof of the Claim.]
The implication from left to right is obvious.\ To prove the other one, suppose that $[x] \leq [y]$. The definition of $T$ guarantees that
\[
[\bot] = \{ \bot \} \qquad [\top] = \{ \top \} \qquad [a] = \{ a \} \qquad [b] = \{ b \} \qquad [c] = \{ c \} \qquad [d] = (\bot, d].
\]
In view of the above display, if $x, y \ne d$, then $[x] = \{ x \}$ and $[y] = \{ y \}$, whence $[x] \leq [y]$ implies $x \leq y$ as desired. Therefore, we consider the case where either $x = d$ or $y = d$. First suppose that $x = d$. If $y \in \{ b, c, d, \top\}$, then $x = d \leq y$ and we are done. Thus, it suffices to show that $y \notin \{ \bot, a \}$. Since there is no element in $[y] = [d] = (\bot, d]$ below $\bot$ or $a$ (the latter, by the maximality of $\bot$ in ${\downarrow} a \cap {\downarrow}d$; see Condition~(\ref{Eq:G-part-first-assumption})), the fact that $[\bot] = \{ \bot \}$ and $[a] = \{ a \}$ implies that $[y] \nleq [\bot], [a]$, thus preventing $y$ from being $\bot$ or $a$ as desired. 
Next suppose that $y = d$. If $x = d$, we are done. Therefore, we suppose that $x \in \{ \bot, a, b, c, \top \}$. In this case, $[x] = \{ x \}$, thus the assumption that $\{ x \} = [x] \leq [y] = [d] = (\bot, d]$ implies that $x \leq d = y$. 
\end{proof}

Together with the fact that $X / T$ validates all the formulas valid in $X$, this allows us to assume that $T$ is the identity relation (otherwise we replace $X$ by $X / T$). Consequently, Claim \ref{Eq:order-G4} specializes to the following:
\begin{align}\label{Eq:order-G4-b}
    \begin{split}
X = \left\{ \begin{array}{ll}
 \{ \bot, a, b, c, d\} \cup ({\downarrow}c \cap (\bot, a)) & \text{if $\top$ does not exist}\\
 \{ \bot, a, b, c, d, \top\} \cup ({\downarrow}c \cap (\bot, a)) & \text{if $\top$ exists,}\\
  \end{array} \right.
\end{split}
\end{align}
where ${\downarrow}c \cap (\bot, a)$ is either empty or a singleton. Bearing in mind that if $\top$ exists, then the subposet of $X$ with the universe $\{ \bot, a, b, c, d, \top \}$ is one of the posets depicted in Figure \ref{Fig:who-is-P3}, we conclude that $X$ is a subposet of one of the posets depicted in Figure \ref{Fig:the-final-shape-of-X-2nd-case} containing $\bot, a, b, c,$ and $d$. Furthermore, when we identify $X$ with a subposet of $Z_2$ or $Z_3$ we assume that it contains $\top$, otherwise we identify it with a subposet of $Z_1$.

\begin{figure}
\begin{tabular}{ccccccc}
\begin{tikzpicture}
    \tikzstyle{point} = [shape=circle, thick, draw=black, fill=black , scale=0.35]

\node[label=below:{$Z_1$}] at (-1.5,-0.5) {};

\node[label=below:{$\bot$}] (c1) at (-1.5,0) [point] {};
\node[label=left:{$x$}] (c2) at (-2.5,1) [point] {};
\node[label=right:{$d$}] (c3) at (-0.5,1) [point] {};
\node[label=right:{$c$}] (c5) at (-0.5,2) [point] {};
\node[label=right:{$b$}] (c6) at (-0.5,3) [point] {};
\node[label=left:{$a$}] (c7) at (-2.5,3) [point] {};


\node[label=below:{$Z_2$}] at (2,-0.5) {};

\node[label=below:{$\bot$}] (a1) at (2,0) [point] {};
\node[label=left:{$x$}] (a2) at (1,1) [point] {};
\node[label=right:{$d$}] (a3) at (3,1) [point] {};
\node[label=right:{$c$}] (a5) at (3,2) [point] {};
\node[label=right:{$b$}] (a6) at (3,3) [point] {};
\node[label=above:{$\top$}] (a7) at (2,4) [point] {};
\node[label=left:{$a$}] (a8) at (1,2) [point] {};

\node[label=below:{$Z_3$}] at (5.5,-0.5) {};

\node[label=below:{$\bot$}] (b1) at (5.5,0) [point] {};
\node[label=left:{$x$}] (b2) at (4.5,1) [point] {};
\node[label=left:{$a$}]  at (4.5,2) [point] {};

\node[label=right:{$d$}] (b3) at (6.5,1) [point] {};
\node[label=right:{$c$}] (b5) at (6.5,2) [point] {};
\node[label=right:{$b$}] (b6) at (6.5,3) [point] {};
\node[label=above:{$\top$}] (b7) at (4.5,3) [point] {};

  \draw (c5) -- (c2) -- (c1) -- (c3) -- (c5)-- (c6)  (c7) -- (c2)  ;

  \draw (a2) -- (a1) -- (a3) -- (a5)-- (a6) -- (a7) (a5) -- (a2) --(a8) -- (a7) (b5) -- (b7) -- (b2) ;
  
      \draw (b5) -- (b2) -- (b1) -- (b3) -- (b5)-- (b6)  ;
   
\end{tikzpicture}
\end{tabular}
\caption{The posets $Z_1$, $Z_2$, and $Z_3$.}
\label{Fig:the-final-shape-of-X-2nd-case}
\end{figure}
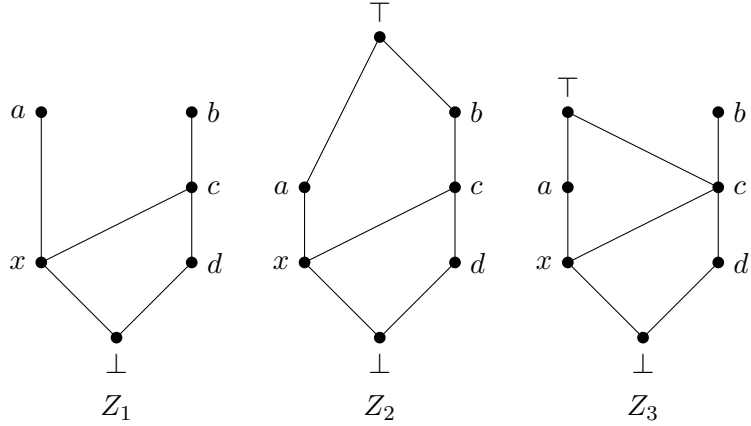

Now, recall that $X$ validates the Jankov formulas of $G_1, \dots, G_6$, and hence none of these posets is a p-morphic image of an upset of $X$ by the Fine Lemma. Bearing this in mind, we begin by considering the case where $X$ is a subposet of $Z_1$ containing $\bot, a, b, c,$ and $d$. In this case, $X$ is isomorphic to either $G_1$ or $G_2$, a contradiction. Next we consider the case where $X$ is a subposet of $Z_2$ (resp.\ $Z_3$) containing $\bot, a, b, c, d,$ and $\top$. In this case, $X$ is isomorphic to either $G_3$ or $G_4$ (resp.\ $G_5$ or $G_6$), which is also false. Hence, we reach the desired contradiction. 
\end{proof}

We are now ready to prove that $\mathsf{KG}$ can be axiomatized by Jankov formulas.

\begin{proof}[Proof of Theorem~\ref{Cor:KG-Jankov-axiom}]
Let $\Sigma$ be the union of the set of Jankov formulas that axiomatize ${\sf BW}_2$ 
and the set of Jankov formulas of the posets in Figures \ref{Fig:posets-K} and \ref{Fig:posets-G}. We will prove that $\Sigma$ axiomatizes $\mathsf{KG}$. 

Suppose the contrary. Since sums of one-generated Heyting algebras have width $\leq 2$, we have that ${\sf BW}_2 \subseteq \mathsf{KG}$, and hence the Jankov formulas axiomatizing ${\sf BW}_2$  belong to $\mathsf{KG}$. Furthermore, observe that the posets in Figures~\ref{Fig:posets-K} and~\ref{Fig:posets-G} are not models of $\mathsf{KG}$ (because each of them contains one of the posets $P_1, P_2, P_3$ in Figure~\ref{Fig:KG-axiom} as a subposet and $\mathsf{KG}$ is axiomatized by $\beta(P_1), \beta(P_2), \beta(P_3)$ by Theorem \ref{Thm:axioms-subframe-KG}). Therefore, in view of the Dual Jankov Lemma, the Jankov formulas of these posets belong to $\mathsf{KG}$. As a consequence, we obtain that $\Sigma \subseteq \mathsf{KG}$. Since by assumption $\Sigma$ does not axiomatize $\mathsf{KG}$, this yields that the si-logic $\mathsf{L}$ axiomatized by $\Sigma$ is strictly contained in $\mathsf{KG}$. \color{black}

From Theorem \ref{Thm:width2-kracht} it follows that ${\sf BW}_2\subseteq\mathsf{L}$. 
By Fine Completeness Theorem, there is a class $K$ of rooted Noetherian posets of width $\leq 2$ such that $\mathsf{L} = \textup{Log}(K)$.
Since $\mathsf{KG} \nsubseteq \mathsf{L}$, by Theorem \ref{Thm:axioms-subframe-KG} there is a poset $X \in K$ refuting $\beta(P_i)$ for some $i \leq 3$. Because $X$ has width $\leq 2$, we have that $X \vDash \beta(P_1)$. Therefore, either $X \nvDash \beta(P_2)$ or $X \nvDash \beta(P_3)$. By Lemmas \ref{Lem:jankov-axioms-KG} and \ref{Lem:jankov-axioms-KG-2}, there is a poset $Y$ in Figure \ref{Fig:posets-K} or \ref{Fig:posets-G} such that $X \nvDash \mathcal{J}(Y)$. Since $X \in K$ and $\mathsf{L} = \textup{Log}(K)$, we obtain
that $\mathcal{J}(Y)$ does not belong to $\mathsf{L}$. The obtained contradiction proves that ${\sf KG}={\sf IPC}+\Sigma$.
\end{proof}

\end{document}